\documentclass[11pt]{article}
\usepackage{pgfplots}
\usepackage{amsthm}
\usepackage{array}
\usepackage{multirow}
\usepackage{longtable}
\usepackage{tablefootnote}
\usepackage{helvet}
\usepackage{amsmath}
\usepackage{amsfonts}
\usepackage{amssymb}
\usepackage{gensymb}
\usepackage{enumitem}
\usepackage{fancyvrb}
\usepackage{tikz}
\usetikzlibrary{fadings,patterns}
\usepackage{caption}
\usepackage{comment}
\usepackage[a4paper, portrait, margin=1.5in]{geometry}
\usepackage{graphicx}

\newcommand{\JE}{\hfill J.\,E.}

\newcommand{\OP}{\operatorname}

\newcommand{\RR}{\mathbb{R}}
\newcommand{\ZZ}{\mathbb{Z}}

\renewcommand{\sin}{\operatorname{sin.}}
\renewcommand{\cos}{\operatorname{cos.}}
\newcounter{eulercounter}

\newenvironment{LemmaUnnum}[1][]{\refstepcounter{eulercounter}\par\medskip
  \begin{center}\Large Lemma.\end{center}\medskip \qquad\S. \theeulercounter.\,\,\,\,}{}

\newenvironment{Problem}[1][]{\refstepcounter{eulercounter}\par\medskip
  \begin{center}\Large Problem #1.\end{center}\medskip \qquad\S. \theeulercounter.\,\,\,\,}{}

\newenvironment{Scholium}[1][]{\refstepcounter{eulercounter}\par\medskip
  \begin{center}\large Scholium. #1\end{center}\medskip \qquad\S. \theeulercounter.\,\,\,\,}{}

\newenvironment{SolutionUnnum}[1][]{\par\medskip
  \begin{center}\large Solution.\end{center}\medskip\par}{}
\newenvironment{Corollary}[1][]{\refstepcounter{eulercounter}\par\medskip
  \begin{center}\large Corollary #1.\end{center}\medskip \qquad\S. \theeulercounter.\,\,\,\,}{}
\newenvironment{CorollaryUnnum}[1][]{\refstepcounter{eulercounter}\par\medskip
  \begin{center}\large Corollary.\end{center}\medskip \qquad\S. \theeulercounter.\,\,\,\,}{}

\newenvironment{Example}[1][]{\refstepcounter{eulercounter}\par\medskip
  \begin{center}\large Example. #1\end{center}\medskip \qquad\S. \theeulercounter.\,\,\,\,}{}

\newcommand{\geom}[1]{\mathrm{#1}}

\newlength{\leftbarwidth}
\newlength{\leftbarsep}
\usepackage{xcolor}
\usepackage{framed}
\DeclareGraphicsExtensions{.eps}
\begin{document}
\pagestyle{empty}
\VerbatimFootnotes

\begin{center}
  {\Huge Euler's explorations of\\ extremal ellipses:

    \bigskip
    
    \Large Translations of and commentary on E563, E691 and E692
  }

  \bigskip
  
  Jonathan David Evans,\\ School of Mathematical Sciences,
  Lancaster University \verb|j.d.evans@lancaster.ac.uk|\\

\end{center}
\bigskip

\VerbatimFootnotes

\begin{abstract}
  In the 1770s, Euler wrote a series of papers (E563, E691 and
  E692) about finding the ellipse with minimal area or perimeter
  in the family of all ellipses passing through a fixed set of
  points. This is a translation of all three papers from the
  original Latin, together with a commentary which discusses
  Euler's results and an appendix which addresses a question
  from E691 which Euler left for others to consider.
\end{abstract}

\section{Introduction}

In the 1770s, Euler worked on a series of extremal problems for
families of ellipses. Namely, he fixes a finite set of points
and looks at the family of ellipses passing through these
points; he asks which amongst these ellipses has the smallest
area or perimeter. The main strategy is of course to find a
convenient expression for the area or perimeter of the ellipse
in terms of the parameter or parameters of the family,
differentiate with respect to those parameters and set the
result equal to zero, then solve for the parameters. In each
case, this strategy poses difficulties which make the problem
interesting. The results were published in the following papers:
\begin{itemize}[leftmargin=30.5pt]
\item[E563] {\em De ellipsi minima dato parallelogrammo
    rectangulo circumscribenda.}\\
  (On the minimal ellipse circumscribing a given right-angled parallelogram.)

  \medskip

  This paper appeared in the Acta of the St Petersburg Academy
  in 1780, but was written in 1773. It deals with the case when
  the ellipses are required to pass through four points at the
  vertices of a rectangle; there is a one-parameter family of
  such ellipses. The ellipse of minimal area is quickly found,
  and the bulk of the paper concerns the much harder problem of
  finding the ellipse of minimal perimeter. The principal
  difficulty is that the formula for the perimeter of the
  ellipse is only known as a power series; in the end Euler only
  gives the solution in reverse: that is, by specifying the
  ellipse and finding the inscribed rectangle for which the
  ellipse minimises area amongst all circumscribed
  ellipses. Moreover, the rectangle is specified by expressing
  its aspect ratio in terms of a power series, so in finite time
  one finds only approximate solutions. But Euler points out
  that if one tabulated the values of the solution sufficiently
  densely, one could work backwards to solve the original
  problem, in the same way one could figure out a number from its
  logarithm using a log table.
\item[E691] {\em Problema geometricum quo inter omnes ellipses
    quae per data quatuor puncta traduci possunt ea quaeritur
    quae habet aream minimam.}\\
  (A geometric problem, in which
  amongst all ellipses which can be drawn through four given
  points, that one is sought which has minimal area.)

  \medskip

  This paper appeared (posthumously) in the Nova Acta of the St
  Petersburg Academy in 1795, but was written in
  1777.\footnote{Publication delays seem to have been a huge
    issue for the St Peterburg Academy throughout the 18th
    century. Already in a 1733 letter to Euler (E864, \cite{Bernoulli}), Daniel Bernoulli
    says {\em Man kann hierdurch sehen wie
      pr\"{a}judierlich es unsern Commentariis ist, so langsam
      gedruckt zu werden, indem wir allzeit als die alte
      Fastnacht nach den andern kommen werden} (``One can see
    from this how prejudicial it is to our Commentarii to be
    printed so slowly, in that, like alte Fastnacht, we always
    come after the others.'' -- alte Fastnacht is a lent
    carnival which starts after the end of the usual lent
    carnival; it is famously celebrated in Bernoulli and Euler's
    hometown of Basel in place of the more traditionally-timed
    Fastnacht).} It deals with the problem of minimal area in
  the case when the ellipses are constrained to pass through
  four points. The four points are arbitrary provided that at
  least one ellipse can be drawn through them. The main
  innovation is a general formula for the area of an ellipse
  given its quadratic equation (in arbitrary coordinates); this
  lets Euler reduce the problem of finding the minimal area
  ellipse to solving a cubic equation. He works out in detail
  the case when the four points are at the vertices of a
  parallelogram, since in this case the cubic simplifies and can
  be fully factorised. He leaves for others the investigation of
  what happens when there are three real roots of the cubic (are
  there three critical ellipses?). We address this problem in
  Appendix \ref{app:cubic}: it turns out there are always three roots and,
  except in certain degenerate cases, they can be interpreted as
  hyperbolas for which an explicit area-like invariant has a
  local maximum.
\item[E692] {\em Solutio problematis maxime curiosi quo inter
    omnes ellipses quae circa datum triangulum circumscribi
    possunt ea quaeritur cuius area sit omnium minima.}\\
  (Solution
  to a problem of greatest curiosity in which, amongst all
  ellipses which can be circumscribed around a given triangle,
  that one is sought, whose area is smallest of all.)

  \medskip
  
  This paper appeared in the same volume of the Nova Acta,
  immediately after E691. It deals with the problem of minimal
  area in the case when the ellipses are constrained to pass
  through three (non-collinear) points. The same formula for
  area discovered in E691 is used to give an ``exceedingly
  simple and easy construction'' of the minimal area ellipse
  passing through the three points. This minimal area ellipse
  enjoys the ``remarkable property'' that its centre is at the
  centre of gravity of the three points and the ratio of the
  area of the ellipse to that of the triangle is
  \(\frac{4\pi}{3\surd 3}\). This result is nowadays
  attributed\footnote{See for example {\cite[\S 99]{Dorrie}},
    {\cite[p.59]{Sommerville}}.} to J. Steiner (1796--1863) who
  proved it in an 1829 article \cite{Steiner2} in Gergonne's
  journal.
\end{itemize}
Of course, area and perimeter are not the only quantities one
could consider minimising, for example one could try to minimise
the eccentricity (as proposed in 1827 in Gergonne's journal
\cite{Gergonne} and solved by Steiner in \cite{Steiner1}; see
D\"{o}rrie's book {\cite[\S 54]{Dorrie}} for a nice
exposition). One could also weaken the constraint by asking for
the ellipse to simply {\em contain} the points (either in its
interior or boundary) -- this usually leads to a smaller
ellipse. This culminates in the theory of the outer
L\"{o}wner--John ellipsoid \cite{John} which is the ellipse of
smallest volume containing a given convex body (in any
dimension). In this commentary, however, I will confine myself
to discussing and expanding on Euler's papers, which seem to be
largely forgotten.

\section{The geometric setting}

I will start by explaining the general geometric setting for all
three papers; I will use modern terminology here because I think
it helps. If we fix a collection of (distinct) points and ask
for the family of all ellipses passing through these points, we
find ourselves in one of several scenarios:
\begin{itemize}
\item If we have fixed five or more points then there is at most
  one ellipse: two distinct ellipses can intersect in at most
  four points.
\item If we have fixed four points then there are either no
  ellipses passing through them (if any of the points lies in
  the convex hull of the other three) or an infinite family
  depending on a single parameter. This is part of a bigger
  family of conics: what we would now call the {\em pencil}
  whose {\em base locus} is the given collection of points; this
  pencil will include hyperbolas, possibly parabolas, and some
  singular conics (unions of straight lines). Since Euler is
  interested in bounded areas and lengths, he considers only the
  ellipses in the pencil. We will loosely refer to this as a
  {\em pencil of ellipses}.
\item If we have fixed three (non-collinear) points then there
  is a two-parameter family of ellipses passing through
  them. This is part of a {\em net} of conics.
\end{itemize}
When talking about the base locus of our pencil or net, it will
be convenient to use the word {\em conelliptic} to describe the
case when there is at least one ellipse through the points. This
guarantees that the points are ``in general position'' in the
sense that no three are collinear.

The general equation of a plane conic (in Euler's notation) is
\[Axx+2Bxy+Cyy+2Dx+2Ey+F=0.\]
The condition that it be an ellipse is
\[B^2<AC.\] If two quadratic forms like this vanish at the same
point then any linear combination of them also vanishes at that
point, so imposing point constraints always yields a {\em
  linear} space of equations, but we only care about the
equation up to scale, so it gives a {\em projective space} of
curves. We call such a family of curves a {\em linear system} of
dimension \(d\) if the projective space has dimension \(d\). A
pencil is a 1-dimensional linear system and a net is a
2-dimensional linear system.

\begin{figure}[htb]
  \centering
  (a)
  $\vcenter{\hbox{\includegraphics[width=5cm]{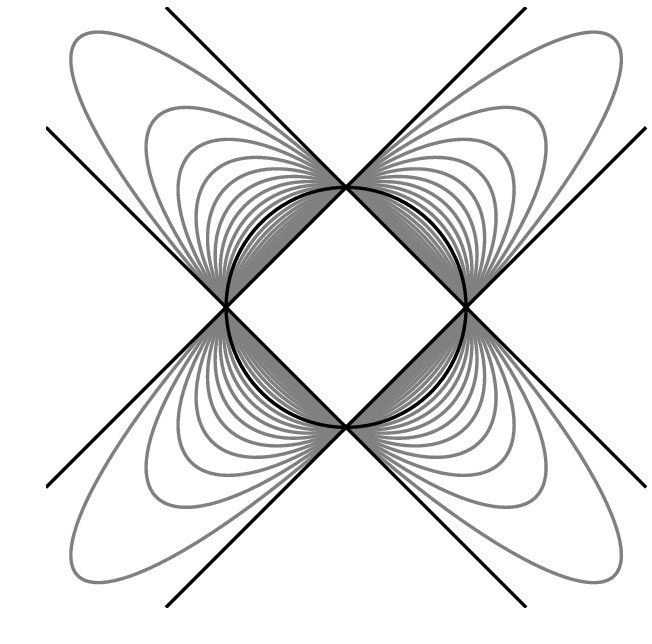}}}$\quad
  (b)
  $\vcenter{\hbox{\includegraphics[width=5cm]{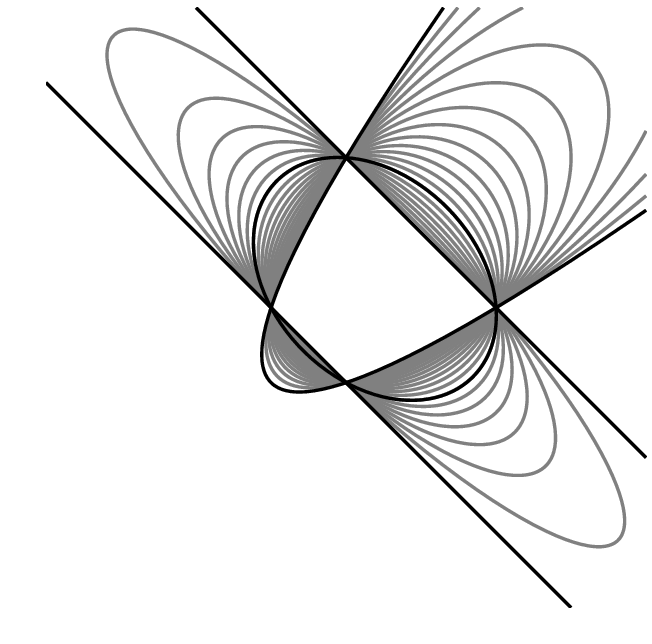}}}$\\
  (c)
  $\vcenter{\hbox{\includegraphics[width=5cm]{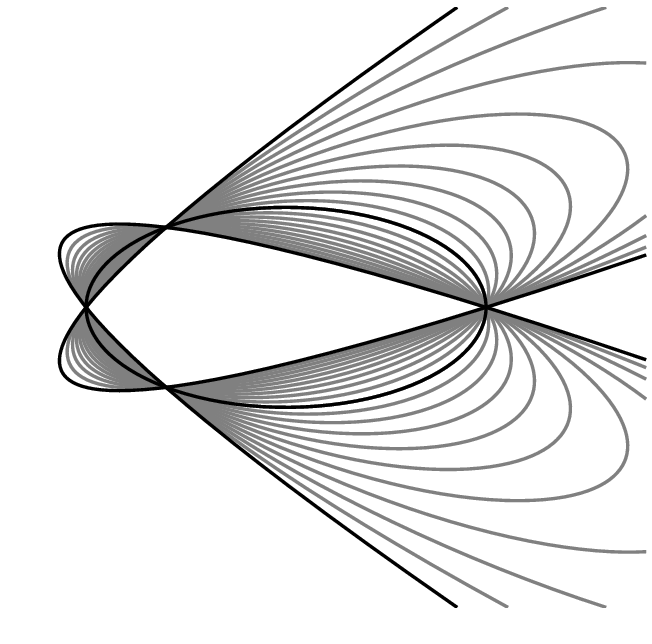}}}$\quad
  (d)
  $\vcenter{\hbox{\includegraphics[width=5cm]{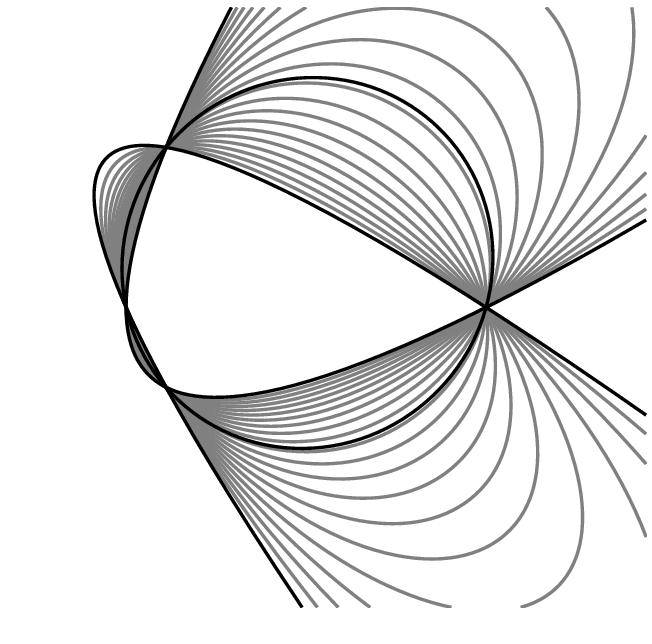}}}$
  \caption{Four pencils of ellipses through the points
    \(\geom{A}=(a,0)\), \(\geom{B}=(b,0)\), \(\geom{C}=(0,c)\)
    and \(\geom{D}=(0,d)\). In each case, both the minimal
    ellipse and the limiting conics are shown in black. The
    quadrilateral (base locus) is: \\ \hspace*{1ex} (a) a square
    (\(a=-b=c=-d=1\)); \\ \hspace*{1ex} (b) a trapezium
    (\(a=c=2\), \(b=d=-1\)); \\ \hspace*{1ex} (c) a kite
    (\(a=4\), \(b=d=-1\), \(c=1\)); \\ \hspace*{1ex} (d)
    irregular (\(a=2\), \(b=-1\), \(c=4\), \(d=-1/2\)).}
  \label{fig:pencils}
\end{figure}

Figure \ref{fig:pencils} shows some pencils of ellipses through
four points. The parameter \(B\) in these pencils varies from
\(-\sqrt{AC}\) to \(\sqrt{AC}\). When it takes on one of these
limiting values the conic we obtain is either a parabola or a
pair of lines; I have highlighted these limiting conics in
black. I have also highlighted in black the minimal area ellipse
in each pencil (found using Euler's method from E691).

%\newpage

\section{The minimal area problem}

\subsection{E563}

In E563, the four points sit at the vertices of a
rectangle, and Euler chooses coordinates to make these points
sit at \((\pm f,\pm g)\) for some \(f,g\). In this case, the
equation of the ellipse passing through the four points
simplifies to:
\[\frac{x^2}{a^2}+\frac{y^2}{b^2}=1.\] 
Euler points out that \(a\) and \(b\) are not independent: they
must be chosen so that
\[\frac{f^2}{a^2}+\frac{g^2}{b^2}=1\] in order to satisfy the
point constraint. This leaves the single degree of freedom I
alluded to. The formula for the area of the ellipse is very
simple: it is just \(\pi ab\). Euler treats this as a
constrained optimisation problem.

I will explain his calculation in a little more detail because
it is the model for his later calculation on minimising
perimeter. The vanishing of the differential of area gives
\[a\,db+b\,da = 0.\] Since the constraint is always satisfied,
its differential of the constraint must vanish:
\[-\frac{2f^2}{a^3}\,da-\frac{2g^2}{b^3}\,db=0,\] so that
\(da\,:\,db=a\,:-b\), and ``with the proportional \(a\) and
  \(-b\) written in place of \(da\) and \(db\)'' we get
  \[-\frac{2f^2}{a^2}+\frac{2g^2}{b^2}=0,\] or
  \(\frac{g^2}{b^2}=\frac{f^2}{a^2}\). Substituting this into
  the constraint tells us \(2f^2/a^2=2g^2/b^2=1\), so the result
  is that \(a=f\surd 2\) and \(b=g\surd 2\), which means
  geometrically that the tangents to the ellipse at the vertices
  of the rectangle are parallel to the lines connecting the
  midpoints of the sides of the rectangle (as Euler explains in
  \S.2).

\subsection{E691 and E692: the equations of the ellipses}
\label{sct:equations}

The issue with three or four {\em general} points is that the
equation for the ellipse will be more complicated, and one also
needs a formula for its area. This put Euler off for several
years:

\begin{quote}
  ```...but the general problem I did not venture to attempt at the
  time on account of the vast number of quantities which needed
  to be introduced in calculations, whence thoroughly
  inextricable analytic formulas arose.''
\end{quote}

In E691 Euler shows that, if we have four point constraints then
we can choose coordinates (and scaling of the equation) so that
the coefficients \(A\), \(C\), \(D\), \(E\) and \(F\) of the
conic are all fixed and only \(B\) varies. Namely, Euler writes
\(\geom{A}\), \(\geom{B}\), \(\geom{C}\), and \(\geom{D}\) for
the points\footnote{Euler confusingly uses the same letters for
  points and for variables. In the translations, I have tried to
  ameliorate this by using \(\geom{A}\) for a point and \(A\)
  for a variable. I use the same convention here for
  consistency.}; we can assume (maybe after relabelling) that
the lines \(\geom{A}\,\geom{B}\) and \(\geom{C}\,\geom{D}\)
intersect at some point \(\geom{O}\) and use these lines as
coordinate axes. The points are now at \((a,0)\), \((b,0)\),
\((0,c)\) and \((0,d)\), so substituting these into the equation
for the conic Euler derived
\begin{align*}
  Aaa+2Da+F&=0,& Abb+2Db+F&=0,\\
  Ccc+2Ec+F&=0,&Cdd+2Ed+F&=0,
\end{align*}
which we can think of simultaneous equations
expressing \(A,C,D,E,F\) in terms of a single free variable,
which we can scale to be equal to \(1\). This gives
\[A=cd,\quad C=ab,\quad D=-\frac{1}{2}cd(a+b),\quad
  E=-\frac{1}{2}ab(c+d),\quad F=abcd,\] and \(B\) is free. 

In E692, with three points \(\geom{A}\), \(\geom{B}\), \(\geom{C}\),
Euler uses the coordinate axes \(\geom{B}\,\geom{A}\) and
\(\geom{B}\,\geom{C}\), with origin at \(\geom{B}\). The three points
are then at \((0,0)\), \((a,0)\) and \((0,c)\), and the
coefficients need to satisfy:
\[F=0,\quad Aaa+2Da=0,\quad Ccc+2Ec=0,\]
so the equation is
\[Axx+2Bxy+Cyy-Aax-Ccy,\] which has three free variables, but
only two up to scale.

From a modern point of view, Euler is exploiting the affine
symmetry of the problem. Affine transformations send ellipses to
ellipses and act in a predictable way on areas (scaling
them). In fact, affine transformations act transitively on
triangles, so in E692 Euler could have reduced without loss of
generality to the case of an equilateral triangle. From this
point of view, his very elegant result becomes intuitively very
clear: the minimal area ellipse circumscribing an equilateral
triangle should be a circle centred at the centre of gravity of
the triangle. The area ratio between this circle and the
triangle is \(\frac{4\pi}{3\surd 3}\). Since affine
transformations send centres of gravity to centres of gravity,
and preserve area ratios, it is clear that we should have
expected Euler's remarkable property purely from symmetry
considerations.

\subsection{E691 and E692: the area formula}

The other key insight from E691 and E692 was the discovery of an

\begin{quote}
  ``expression...most worthy of note, because with its help the
total area of all ellipses can be sufficiently quickly assigned
from only the equation between its coordinates, whether they are
rectangular or obliquangular.''
\end{quote}

\noindent This formula for the area in terms of the equation for
the conic is:
\begin{equation}
  \label{eq:area_formula}
  \pi\sin\omega\left(\frac{CDD+AEE-2BDE}{(AC-BB)^{\frac{3}{2}}}-\frac{F}{\surd(AC-BB)}\right),
\end{equation}
where \(\omega\) is the angle \(\geom{A\,O\,C}\) and one should
choose the square root in the denominator to make the overall
sign positive. The denominator gives rise to poles at
\(B=\pm\sqrt{AC}\) and the fractional power in the denominator
means that the function is only real-valued for
\(B\in (-\sqrt{AC},\sqrt{AC})\); let us plot the graphs of this
function against \(B\) for the pencils shown in Figure
\ref{fig:pencils}. In these plots we actually extend the domain
to \(B^2>AC\) by replacing the denominator with
\(\sqrt{B^2-AC}\); we will revisit both Euler's area formula and
what it means for hyperbolas in Appendix \ref{app:cubic}.

\begin{figure}[htb]
  \centering
  (a)
  $\vcenter{\hbox{\includegraphics[width=5cm]{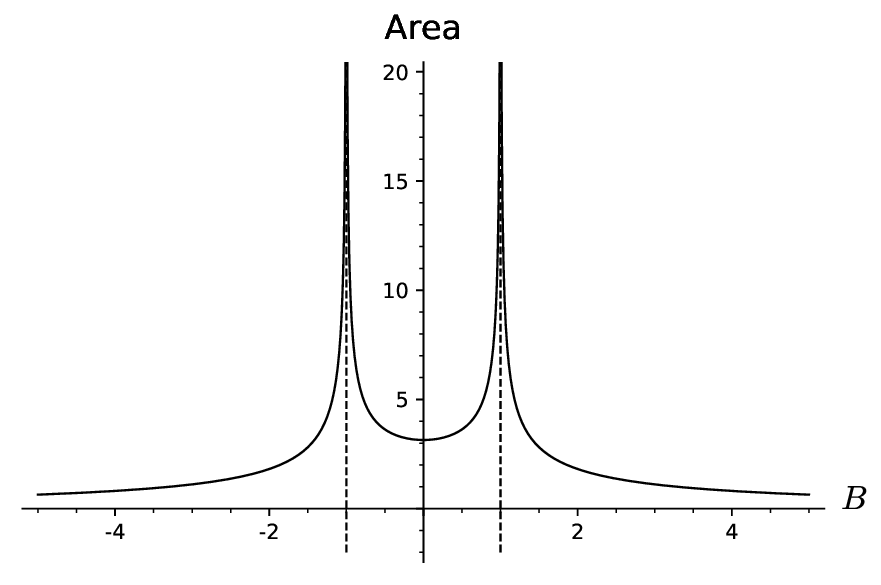}}}$\quad
  (b)
  $\vcenter{\hbox{\includegraphics[width=5cm]{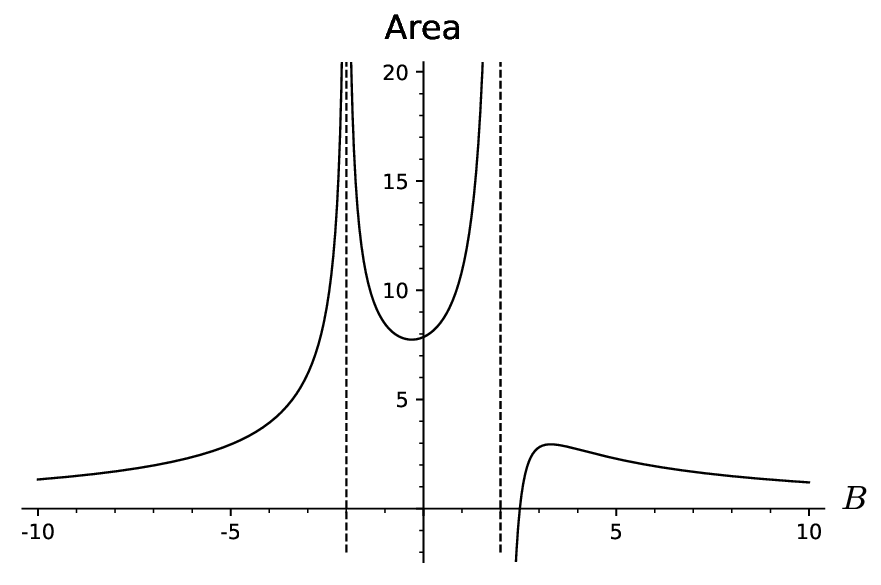}}}$\\
  (c)
  $\vcenter{\hbox{\includegraphics[width=5cm]{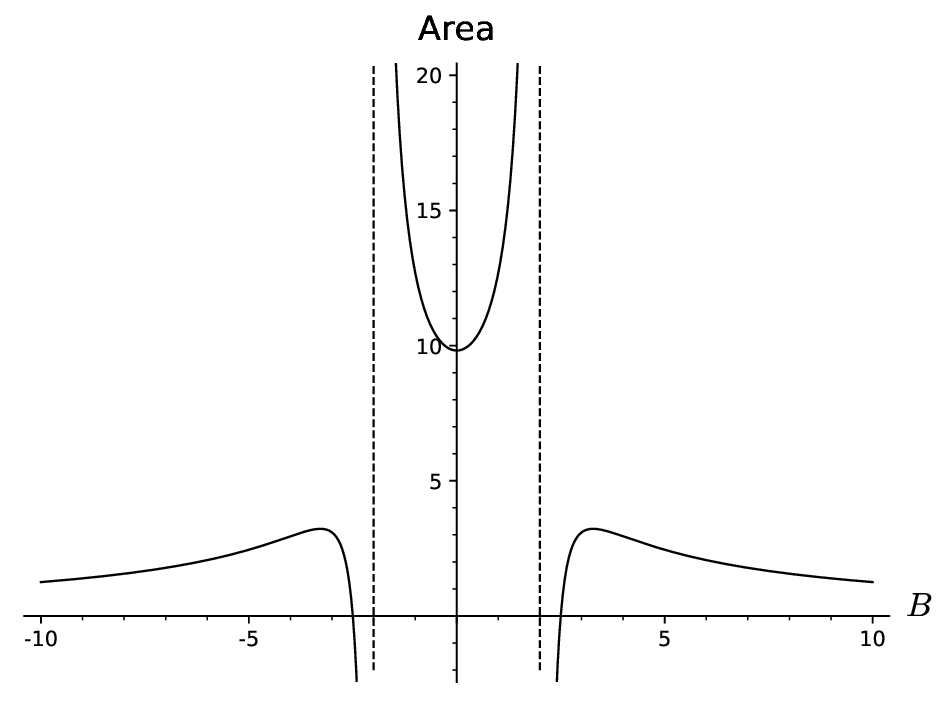}}}$\quad
  (d)
  $\vcenter{\hbox{\includegraphics[width=5cm]{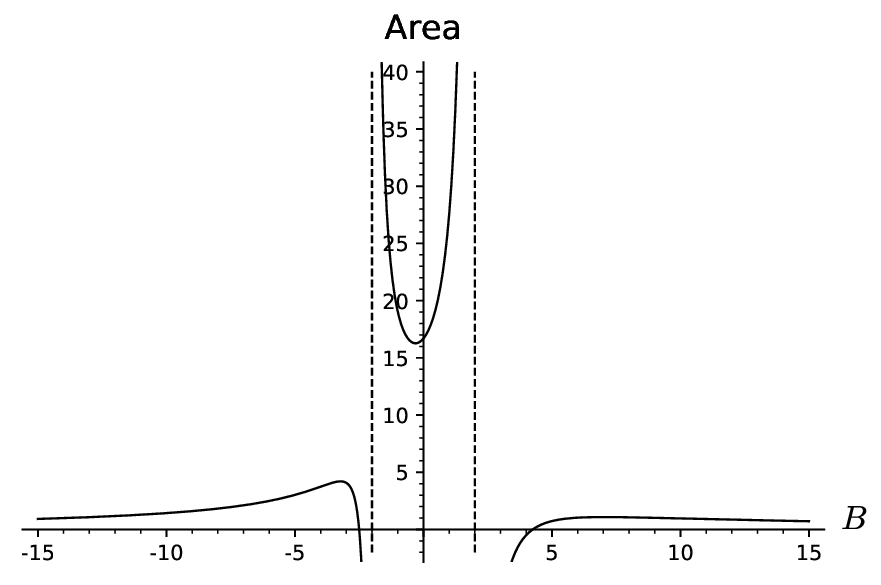}}}$
  \caption{The graphs of area against \(B\) for the four pencils
    in Figure \ref{fig:pencils}. There are poles at
    \(B=\pm\sqrt{AC}\), corresponding to the limiting conics of
    infinite area in the pencil. In each case you can see a
    clear unique minimum of the area between
    \(\pm\sqrt{AC}\). We have continued the function beyond that
    by using the denominator \((B^2-AC)^{3/2}\) instead of
    \((AC-B^2)^{3/2}\).}
  \label{fig:pencils_areas}
\end{figure}
\newpage

\subsection{E691: the minimal area ellipse}

In E691, Euler differentiates Equation \eqref{eq:area_formula}
with respect to \(B\), and sets the result equal to zero:
\[\frac{2DE}{(AC-BB)^{\frac{3}{2}}} -
  \frac{BF}{(AC-BB)^{\frac{3}{2}}} +
  \frac{3B(\Delta-2BDE)}{(AC-BB)^{\frac{5}{2}}} = 0,\]
where \(\Delta=CDD+AEE\), which yields the cubic equation
\[FB^3-4DEB^2+(3CD^2+3AE^2-ACF)B-2ACDE=0.\] Euler now argues
that since

\begin{quote}
  ``this always has real roots, it is certain howsoever the
four points be distributed, one ellipse can always be assigned
passing through these four points, whose area is smallest of
all.''
\end{quote}

\noindent This argument is not quite sufficient, because to give
a minimal {\em ellipse} (rather than hyperbola) we need
\(B^2<AC\), and there is no {\em a priori} guarantee that this
condition will be satisfied for any of the roots. Euler also
says that

\begin{quote}
  ``if perchance it were to happen that this cubic were
to admit three real roots, just as many solutions will take
place, the nature of which, however, I leave the examining to
others.''
\end{quote}

\noindent We will take up the Euler's thread in Appendix \ref{app:cubic},
and show that the cubic always has three roots: one satisfying
\(B\leq -\sqrt{AC}\), one in the interval
\((-\sqrt{AC},\sqrt{AC}\), and one in the interval
\(B\geq \sqrt{AC}\). We will also give a geometric
interpretation of the quantity in Equation
\eqref{eq:area_formula} for the hyperbolas in the pencil.

In \S\S.13--23, Euler works out in detail the situation when the
four points are at the vertices of a parallelogram. Recall that
in Euler's coordinate system, this means that the vertices are
at the points \((\pm a,0)\) and \((0,\pm c)\). The minimal
ellipse corresponds to the element \(B=0\) in the pencil, that
is
\[c^2x^2+a^2y^2=a^2c^2,\] or
\[\frac{x^2}{a^2}+\frac{y^2}{c^2}=1.\] Note, however, that
Euler's coordinates are not orthogonal: they meet at an angle
Euler calls variously \(\omega\) or \(\theta\). The area of this
ellipse is \(\pi ac\sin\theta\), whilst the area of the
parallelogram is \(2 ac \sin\theta\). Euler spends \S\S.18--20
converting to an orthogonal coordinate system to determine the
principal semiaxes \(f\) and \(g\):
\[f^2=\frac{a^2c^2\OP{sin}^2\theta}{c^2\OP{sin}^2(\theta-\phi)+a^2\OP{sin}^2\phi}
  \quad\text{and}\quad
  g^2=\frac{a^2c^2\OP{sin}^2\theta}{c^2\OP{cos}^2(\theta-\phi)+a^2\OP{cos}^2\phi},\]
where \(\phi\) is the angle the major axis makes with the line
\(\geom{O\,A}\), which Euler determines via the equation
\[\cot 2\phi=\frac{a^2+c^2\OP{cos} 2\theta}{c^2\OP{sin} 2\theta}.\] There
is a typo in the original: Euler writes \(\cos 2\phi\) instead
of \(\cot 2\phi\).

In \S\S.21--22 derives two well-known properties expressed by
the relations
\[f^2g^2=a^2c^2\OP{sin}^2\theta\quad\text{and}\quad
  a^2+b^2=f^2+g^2.\] The first relation implies that the

\begin{quote}
``parallelogram described around two conjugate diameters is
ascribed area equal to the parallelogram described around the
principal axes,''
\end{quote}

\noindent the second that

\begin{quote}
  ``in every ellipse the sum of
the squares of two [conjugate] diameters is always equal to the
sum of the squares of the principal axes.''
\end{quote}

\noindent Recall that two diameters \(D_1\) and \(D_2\)
(i.e. chords passing through the centre) of an ellipse are said
to be {\em conjugate} if \(D_1\) intersects all chords parallel
to \(D_2\) and vice versa, which happens if and only if the
tangents of the ellipse at the endpoints of \(D_1\) are parallel
to \(D_2\) and vice versa. The fact that \(\geom{A\,B}\) and
\(\geom{C\,D}\) are conjugate diameters follows from the
equation \(c^2x^2+a^2y^2=a^2c^2\) of the ellipse: these
diameters are just the \(x\)- and \(y\)-axes, and the tangents
at the \(x\)- (respectively \(y\)-) intercepts are vertical
(respectively horizontal).

Finally, in \S.23, Euler specialises to the case of a
rectangular parallelogram, and observes that he recovers the
minimal area result from E563. Unfortunately, to compare the two
papers, one must switch the notation: in E563, he uses \(a\) and
\(b\) for the semi-axes of the ellipse and \(f\) and \(g\) for
the half-sides of the rectangle, whilst for E691 he uses \(f\)
and \(g\) for the semiaxes of the ellipse and
\(a\OP{cos}(\theta/2)\) and \(b\OP{sin}(\theta/2)\) for the
half-sides of the rectangle.

\subsection{E692: the minimal area ellipse}

Having reviewed his material from E691 on the equations for
ellipses and the area formula, Euler proceeds (\S.4) to adapt it
to the case of ellipses through three points. In his chosen
coordinates, with origin at \(\geom{B}\), he gets the equation
\[Axx+2Bxy+Cyy-Aax-Ccy\] for the general ellipse in his net. He
introduces variables \(\phi\) and \(s\) such that
\(B=\OP{cos}\phi\surd AC\) (note that \(B^2<AC\) means that
such a \(\phi\) exists) and \(C=As^2\), so the area formula
becomes
\[\frac{\pi}{4}\sin\omega\left(\frac{a^2+c^2s^2 -
      2acs\OP{cos}\phi}{2\OP{sin}^3\phi}\right),\] where
\(\omega\) is the angle \(\geom{A\,B\,C}\). By this change of
variables, he has now eliminated the freedom of scaling the
equation, and so we just need to take partial derivatives with
respect to \(s\) and \(\phi\) and set them equal to zero. In
\S.7, he shows that the vanishing of the \(s\)-derivative
happens if \(s=\frac{c}{a}\), and if we take \(A=c^2\) (by
scaling the equation) then we get \(C=a^2\). In \S.8 he shows
that the \(\phi\) derivative vanishes if and only if
\(\phi=60\degree\). This completely solves the problem, giving
the minimal area to be
\(\frac{2\pi ac\OP{sin}\omega}{3\surd 3}\), whose ratio to the
area \(\frac{1}{2}ac\OP{sin}\omega\) of the triangle is
\(\frac{4\pi}{3\surd 3}\approx 2.41840\), and in \S.9 he gives
the first few convergents of the continued fraction for this
number. In \S\S.10--11, Euler works from the equation
\[c^2x^2+acxy+a^2y^2=ac^2x+a^2cy\] of the minimal ellipse to
deduce that its centre is at the centre of gravity of the
triangle \(\geom{A\,B\,C}\) and that the tangents of the ellipse
at the vertices of this triangle are parallel to the opposite
sides of the triangle. In \S\S.12--14, Euler works out the cases
of an equilateral triangle (where he shows the minimal area
ellipse is a circle) and an isosceles triangle.

\section{The minimal perimeter problem}

The problem of minimising perimeter is substantially harder than
minimising area, and Euler doesn't fully solve it to his
satisfaction: E563 is rather an ``attempt at solving it''. The
issue with his solution is that, rather than giving a
closed-form expression for the minimal perimeter ellipse in
terms of the data of the rectangle, he is able to express the
data for the rectangle in terms of the data for the minimal
ellipse, and the expression is an infinite power series. Of
course, for the purpose of getting an {\em approximate} answer,
this is sufficient, as Euler works out in several examples
towards the end of the paper. Below, we will illustrate that,
nowadays with a computer, Euler's method allows us to plot the
dependence of the ellipse on the rectangle as accurately as we
would ever want.

The reason that this problem is so much harder than the problem
of minimising area is that the formula for the perimeter of an
ellipse is in the form of a power series. Indeed, the
rectification of the ellipse (i.e. determining its perimeter) is
a problem which absorbed Euler for much of his mathematical
life; indeed Volumes 20--21 of the Opera Omnia are largely
dedicated to papers about the integrals that arose from his
study of arc-length of conic sections. In \S.3, he says:

\begin{quote}
  ``yet one may not even enter upon the other unless the perimeter
of every ellipse were first expressed neatly in general as an
infinite series, which converges in all cases as far as
possible. However, although several such series can already be
found, nevertheless hardly would any be discovered which would
snatch the palm [of victory] from the following, which I am
about to give.''
\end{quote}

\noindent In \S.4, he then derives the series
\[
  \frac{\pi c}{2\surd 2}\left(1-\frac{1.\,1}{4.\,4}\cdot n^2-\frac{1.\,1.\,3.\,5}{4.\,4.\,8.\,8}\cdot n^4-\frac{1.\,1.\,3.\,5.\,7.\,9}{4.\,4.\,8.\,8.\,12.\,12}\cdot n^6-\text{etc.}\right)
\]
for a quarter of the perimeter, which he first derived about a
year beforehand in E448, in a slightly different way. This appears
to be the breakthrough which allowed him to write E563. He was
clearly very excited about this series; for example in E448, he wrote:

\begin{quote}
  ``After I had formerly been much occupied investigating several
infinite series by which the perimeter of every ellipse would be
expressed, I had hardly expected to be able to dig up such a
series yet simpler and more adapted to calculation than I gave
here and there in Comment. Petrop. or in Actis Berolin.''
\end{quote}

\noindent Adolf Krazer, the editor of the Opera Omnia edition,
remedies Euler's casual approach to citation by suggesting E52
\cite{E52} in Comment. Petrop. and E154 \cite{E154} in Acta
Berolin. Using this power series, Euler digresses in \S.5 and
\S.6. In \S.5 he points out that when \(n=1\) the ellipse
becomes infinitely thin, so that this power series is just
computing the length of an interval; this yields the beautiful
infinite product formula
\[\frac{\pi}{2\surd 2}=\frac{4.\,4}{3.\,5}\cdot\frac{8.\,8}{7.\,9}
  \cdot\frac{12.\,12}{11.\,13}\,\,\text{etc.}\] In \S.6, he
derives a differential equation for \(s\), namely:
\begin{equation}\label{eq:diff_eq_s}
  4n\frac{d^2s}{dn^2}+4\frac{ds}{dn}+\frac{ns}{1-n^2}=0.
\end{equation}
This equation will be used in \S.10 and in \S.12 to derive some
further differential equations, which ultimately do not help
Euler to solve the problem.

With all this in place, the meat of the argument appears in \S
7. Euler proceeds in the same way as before: differentiating
both the formula for perimeter and the constraint, setting both
differentials equal to zero, and attempting to solve for \(a\)
and \(b\). Instead of differentiating with respect to \(a\) and
\(b\), he uses instead the more convenient parameters
\[c=\sqrt{a^2+b^2}\quad\text{and}\quad
  n=\frac{a^2-b^2}{a^2+b^2}.\] He also introduces
analogous quantities for the rectangle
\[h^2=f^2+g^2\quad\text{and}\quad i=\frac{f^2-g^2}{f^2+g^2},\]
where recall the coordinates of the vertices are at
\((\pm f,\pm g)\).

The quantity Euler is trying to minimise is \(cs\) (ignoring the
factor of \(\pi/(2\surd 2)\)). The vanishing of the differential
of \(cs\) gives
\[s\,dc=\frac{ct}{n}\,dn,\]
where
\[t=-n\frac{ds}{dn}=2\cdot\frac{1.\,1}{4.\,4}\cdot
  n^2+4\cdot\frac{1.\,1.\,3.\,5}{4.\,4\,.8.\,8}\cdot\,n^4
  +6\cdot\frac{1.\,1.\,3.\,5.\,7.\,9}{4.\,4\,.8.\,8.\,12.\,12}\cdot\,n^6.\]
The vanishing of the differential of the (logarithm of) the
constraint becomes
\[\left(\frac{-i}{1-in}+\frac{2n}{1-n^2}\right)\,dn=\frac{2}{c}\,dc\]
and together these equations imply
\[\left(2n^2-in-in^3\right)s=2(1-in)(1-n^2)t.\] Here, the
quantity \(i\) is fixed by the rectangle, and the quantity \(n\)
is sought; both \(i\) and \(n\) take values strictly between
\(0\) and \(1\). Unfortunately, both \(s\) and \(t\) are
infinite power series in \(n\), so it seems hopeless to look for
a solution. Instead, Euler is content to reverse the problem:
fix \(n\) and solve for \(i\),

\begin{quote} ``because if it is
  established for many cases, it will easily allow one to
  determine which value of \(n\) corresponds to any given value
  of \(i\).''
\end{quote}

\noindent If we truncate the power series, we get a good
approximation and Euler gives several examples of such
approximations:
\begin{itemize}
\item If you drop terms of order \(n^4\) and above, you get
  \(i=\frac{7n/4}{1+11n^2/16}\approx \frac{7}{4}n\) valid for
  sufficiently small \(n\), so \[n\approx \frac{4}{7}\cdot i.\]
\item If you drop terms of order \(n^5\) and above, you get
  \[i=\frac{7n/4+n^3/128}{1+11n^2/16}\approx\frac{7}{4}\cdot
    n-\frac{153}{128}\cdot n^3,\]
  which yields
  \[n\approx \frac{4}{7}\cdot i-\frac{306}{2401}\cdot i^3.\]
\item Neither of these approximations are valid for values of \(n\)
that are close to the upper limit \(n=1\). For example, we know
that when \(i\to 1\), so that \(g\to 0\), the rectangle becomes
less and less tall, and so the ellipse should become less and
less tall, so that \(b\to 0\) and hence \(n\to 1\). Euler
suggests the approximation
\[n\approx \frac{4i}{7-3i^2}\] as a compromise:
\end{itemize}

\begin{quote}
  ``Indeed as
  long as \(i\) is an exceedingly small fraction, it will be
  that \(n=\frac{4}{7}i\); and if it were considerably bigger it
  will be that \(n=\frac{4}{7}i+\frac{12}{49}i^3\), which value
  is a little bigger than
  \(\frac{4}{7}\cdot i-\frac{306}{2401}\cdot i^3\); but when
  \(i=1\) it again yields \(n=1\).''
\end{quote}

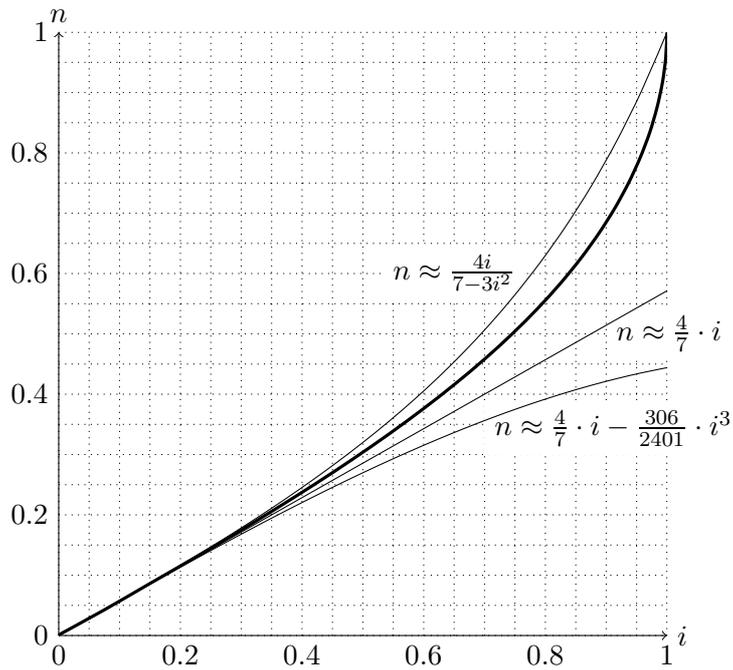
\begin{figure}[htb]
  \centering
  \begin{tikzpicture}[scale=8]
    \draw[->] (0,0) -- (0,1) node [above] {\(n\)};
    \draw[->] (0,0) -- (1,0) node [right] {\(i\)};
    \draw[step=0.05,black,dotted] (0,0) grid (1,1);
    \node[fill=white] at (1-0.1,4/7-0.07) [right] {\(n\approx \frac{4}{7}\cdot i\)};
    \node[fill=white] at (1-0.3,4/7-306/2401-0.1) [right] {\(n\approx \frac{4}{7}\cdot i-\frac{306}{2401}\cdot i^3\)};
    \node[fill=white] at (0.65,0.6) {\(n\approx\frac{4i}{7-3i^2}\)};
    \node at (0,0) [below] {\(0\)};
    \node at (0,0) [left] {\(0\)};
    \node at (1,0) [below] {\(1\)};
    \node at (0,1) [left] {\(1\)};
    \node at (0.2,0) [below] {\(0.2\)};
    \node at (0.4,0) [below] {\(0.4\)};
    \node at (0.6,0) [below] {\(0.6\)};
    \node at (0.8,0) [below] {\(0.8\)};
    \node at (0,0.2) [left] {\(0.2\)};
    \node at (0,0.4) [left] {\(0.4\)};
    \node at (0,0.6) [left] {\(0.6\)};
    \node at (0,0.8) [left] {\(0.8\)};
    \begin{scope}
      \clip (0,0) -- (1,0) -- (1,1) -- (0,1) -- cycle;
      \draw[domain=0:1, smooth, variable=\n] plot ({7*\n/4},{\n});
      \draw[domain=0:1, smooth, variable=\i] plot ({\i},{4*\i/7-306*\i*\i*\i/2401});
      \draw[domain=0:1, smooth, variable=\i] plot ({\i},{4*\i/(7-3*\i*\i)});
      \draw[very thick,domain=0.001:1, smooth, variable=\n, black] plot ({(2*\n*\n*(1-(\n)^2*0.0625-(\n)^4*0.0146484375-(\n)^6*0.0064087-(\n)^8*0.0035799-(\n)^(10)*0.0022821-(\n)^(12)*0.0015809)-2*(1-\n*\n)*(2*(\n)^2*0.0625+4*(\n)^4*0.0146484375+6*(\n)^6*0.0064087+8*(\n)^8*0.0035799+10*(\n)^(10)*0.0022821+12*(\n)^(12)*0.0015809))/(\n*(1+\n*\n)*(1-(\n)^2*0.0625-(\n)^4*0.0146484375-(\n)^6*0.0064087-(\n)^8*0.0035799-(\n)^(10)*0.0022821-(\n)^(12)*0.0015809)-2*\n*(1-\n*\n)*(2*(\n)^2*0.0625+4*(\n)^4*0.0146484375+6*(\n)^6*0.0064087+8*(\n)^8*0.0035799+10*(\n)^(10)*0.0022821+12*(\n)^(12)*0.0015809))},{\n});  
    \end{scope}
  \end{tikzpicture}
  \caption{Euler's approximations for the dependence of the
    quantity \(n\) (encoding the ellipse) in terms of the
    quantity \(i\) (encoding the rectangle). The thick line is
    the approximation
    \(i=\frac{2n^2(s+t)-2t}{n(1+n^2)s-2n(1-n^2)t}\) where \(s\)
    and \(t\) are both truncated by ignoring terms of order
    \(n^{13}\) and higher. Keeping more terms does not affect
    the picture visibly.}
  \label{fig:graph_comparison}
\end{figure}

\noindent With the benefit of modern computers, we can now plot Euler's
approximations, along with higher approximations, to see how the
graphs compare (Figure \ref{fig:graph_comparison}). As Euler
says:

\begin{quote}
  ``in general we will not deviate too far from the
  truth if we were to set \(n=\frac{4i}{7-3ii}\).''
\end{quote}

\noindent As you can see from the graphs, the biggest deviation is
around \(0.1\) when \(n\approx 0.9\), and when \(n<0.3\) the
three approximations are almost indistinguishable from the true
value.

\begin{figure}[htb]
  \centering
  \begin{tikzpicture}[scale=2]
    \draw (0,0) circle [x radius = 2, y radius = 1];
    \begin{scope}[scale=0.51]
      \node (a) at (3.372108,1) {\(\cdot\)};
      \node (b)at (-3.372108,1) {\(\cdot\)};
      \node (d) at (3.372108,-1) {\(\cdot\)};
      \node (c) at (-3.372108,-1) {\(\cdot\)};
      \draw (a.center) -- (b.center) -- (c.center) -- (d.center) -- cycle;
    \end{scope}
  \end{tikzpicture}
  \caption{The ellipse with aspect ratio \(2\,:\,1\)
    circumscribing a rectangle with aspect ratio
    \(3.372108\,:\,1\). In {\cite[\S.11]{E563}}, Euler shows
    that this ellipse minimises perimeter amongst all
    circumscribed ellipses for this rectangle.}
  \label{fig:ellipse_2_1}
\end{figure}
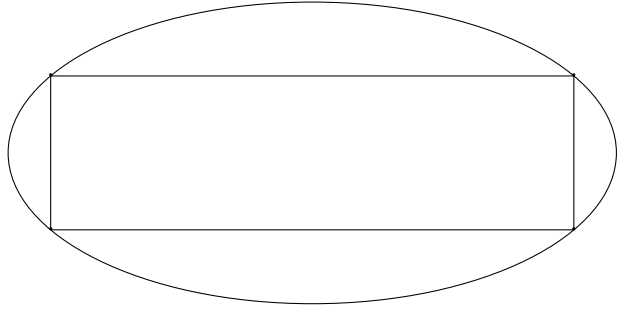

Euler works out an example in \S.11, showing that the ellipse
with \(n=0.6\) (i.e. aspect ratio \(a\,:\,b=2\,:\,1\))
corresponds to the rectangle with \(i\approx 0.838333\)
(i.e. aspect ratio \(f\,:\,g\approx 3.372108\,:\,1\)). See
Figure \ref{fig:ellipse_2_1}.

In \S.9, Euler uses the (second order) differential equation
\eqref{eq:diff_eq_s} satisfied by the power series \(s\) to
deduce a {\em first order} differential equation of Riccati type
satisfied by \(z=s/t=-\frac{1}{ns}\frac{ds}{dn}\), namely:
\[\frac{dz}{dn}=\frac{n}{4(1-n^2)}+\frac{z^2}{n}.\] This is then
used in \S.12 to establish a first order differential equation
relating \(i\) and \(n\), namely:
\[-2n(1-n^2)\frac{di}{dn}=-7n+3n^2+2i(1+3n^2)+i^2n(1-5n^2).\] At
this point, Euler is unable to make any further progress in
solving this differential equation, despite which he concedes
that this has
\begin{quote}
  ``certainly turned out much simpler than one
  had been permitted to hope for at the outset''
\end{quote}

\noindent and avows that
finding a solution of this equation in either exact or integral
form

\begin{quote}
  ``would be judged to have added an extraordinary
  development to analysis.''
\end{quote}

\section{Notes on translations}

\begin{itemize}
\item Footnotes are comments by the translator (\JE).
\item In the original, figures were placed in a table at the
  end; whilst I have moved the figures inline, the original
  figure numbering has been retained.
\item Where the papers say ``Presented to the meeting'', this is
  of course a meeting of the Imperial Academy of Sciences in
  St. Petersburg.
\item In all three papers, Euler uses some letters to denote
  multiple entities. Most commonly, \(A\), \(B\), \(C\),
  etc. are used to denote both points and coefficients. I have
  taken the liberty of distinguishing these typographically
  by using \(\geom{A}\) for a point and \(A\) for a
  coefficient.
\item Euler uses {\em applicata} (``applicate'') instead of
  ``ordinate'' to mean the \(y\)-axis. I have translated this as
  ``applicate''.
\item Euler writes \(\sin\phi^n\) where (for better or for
  worse) we would now usually write \(\OP{sin}^n\phi\). I have
  retained his original notation.
\item I have not always translated the subjunctive
  faithfully. For example I have translated {\em ita ut recta
    \(\geom{A\,F\,E}\) latus \(\geom{B\,C}\) bisecet} as ``so
  that the line \(\geom{A\,F\,E}\) bisects the side
  \(\geom{B\,C}\)'' rather than ``bisect'', which jars to my
  ear. I have also omitted some instances of {\em sit} (``it
  be'') like {\em ita ut sit
    \(\geom{B\,O}=\frac{2}{3}\geom{B\,G}\)} (``so that
  \(\geom{B\,O}=\frac{2}{3}\geom{B\,G}\)'').
\end{itemize}

\section{Acknowledgements}

The author would like to thank the referee for their careful
reading of the manuscript, Yank{\i} Lekili for helpful
discussions about Euler's area formula for ellipses, and Ciara
and Eoin for letting Lewis, Short and Euler come on holiday with
us.

\newpage
\setcounter{eulercounter}{0}

\begin{center}
  {\sc \LARGE On the\\ \Huge minimal ellipse \\ \LARGE circumscribing a given\\ \medskip right-angled parallelogram
  }\\

  \medskip
  \normalsize
  
  An English translation of
  \Large\medskip
  
  {\sc De ellipsi minima dato parallelogrammo rectangvlo circvmscribenda} (E563)
  
  \bigskip

  by Leonhard Euler\\

  \bigskip

  \small
  
  Translated by Jonathan David Evans, School of Mathematical Sciences,\\
  Lancaster University \verb|j.d.evans@lancaster.ac.uk|\\

\end{center}
\bigskip

\noindent {\Huge S}ince infinitely many ellipses may
circumscribe every rectangle, the questions seem not a little
intriguing which ask to find amongst these ellipses either that
of which the area is to be minimal, or that which is to have the
smallest perimeter. In fact, the former works with absolutely no
difficulty; yet nevertheless the solution seems worthy of
attention; but at the same time, the other question of perimeter
is especially hard, so that a perfect solution of it is hardly
to be expected. Therefore to that end it will be all the more
useful to make public an attempt at solving it.

\begin{figure}[htb]
  \centering
  \begin{tikzpicture}
    \draw (3.74,2.2) node [above right] {\(\geom{M}\)} -- (-3.74,2.2) node [above left] {\(m\)} -- (-3.74,-2.2) node [below left] {\(n\)} -- (3.74,-2.2) node [below] {\(\geom{N}\)} -- cycle;
    \draw (0,0) circle [x radius={3.74*sqrt(2)}, y radius = {2.2*sqrt(2)}];
    \draw (0,{-2.2*sqrt(2)}) -- (0,{2.2*sqrt(2)});
    \draw ({-3.74*sqrt(2)},0) -- ({3.74*sqrt(2)},0);
    \node at (0,0) [above right] {\(\geom{C}\)};
    \node at (3.74,0) [below right] {\(\geom{F}\)};
    \node at (-3.74,0) [above left] {\(\geom{f}\)};
    \node at (0,-2.2) [above right] {\(\geom{g}\)};
    \node at (0,2.2) [above left] {\(\geom{G}\)};
    \node at ({3.74*sqrt(2)},0) [right] {\(\geom{A}\)};
    \node at ({-3.74*sqrt(2)},0) [left] {\(\geom{a}\)};
    \node at (0,{2.2*sqrt(2)}) [above] {\(\geom{B}\)};
    \node at (0,{-2.2*sqrt(2)}) [below] {\(\geom{b}\)};
    \draw (0,2.2) -- (3.74,0);
    \draw (0,{2.2*sqrt(2)}) -- ({3.74*sqrt(2)},0);
    \node at (-2.5,3.7) {\textit{Tab. I.}};
    \node at (-2.5,3.2) {\textit{Fig. 1.}};
  \end{tikzpicture}
\end{figure}

\newpage

\begin{Problem}[I.]\marginpar{Tab. I.\\Fig. 1.}
  {\it Around the given rectangle \(\geom{M\,m\,N\,n}\), to describe
    the ellipse \(\geom{A\,a\,B\,b}\) whose area is smallest of all.}
\end{Problem}

\begin{SolutionUnnum}
  With the centre of both the rectangle and the ellipse placed
  at the point \(\geom{C}\), let the halves of the sides of the
  rectangle be called \(\geom{C\,F}=f\) and \(\geom{C\,G}=g\);
  indeed let the semiaxes of the ellipse be \(\geom{C\,A}=a\)
  and \(\geom{C\,B}=b\); and since the point \(\geom{N}\) is
  situated on the ellipse itself, it is well-known from the
  basics to be that \(\frac{f^2}{a^2}+\frac{g^2}{b^2}=1\). Now
  since the area of the ellipse is \(=\pi ab\), it is necessary
  that the quantities \(a\) and \(b\) be determined so that the
  differential of the area vanishes, naturally by having taken
  the semiaxes \(a\) and \(b\) variable; whence we get
  \(a\,db+b\,da=0\), or \(da\,:\,db=a\,:\,-b\), wherefore that
  same equation \(\frac{f^2}{a^2}+\frac{g^2}{b^2}=1\), if
  differentiated and with the proportional \(a\) and \(-b\)
  written in place of \(da\) and \(db\), will yield
  \(\frac{g^2}{b^2}=\frac{f^2}{a^2}\), and hence we deduce
  \(\frac{2f^2}{a^2}=1\), and so \(a=f\surd 2\) and
  \(b=g.\surd 2\). Indeed by having determined the semiaxes this
  ellipse is most easily described.
\end{SolutionUnnum}

\begin{CorollaryUnnum}
  Hence, since the semiaxes will have turned out to be
  proportional to the sides of the rectangle, and
  \(a\,:\,b=f\,:\,g\), it is evident that if the straight lines
  \(\geom{A\,B}\) and \(\geom{F\,G}\) were drawn, they would be
  parallel to one another, by which condition the ellipse is
  then determined.
\end{CorollaryUnnum}

\begin{Scholium}
  With the former question therefore cleared up, yet one may not
  even enter upon the other unless the perimeter of every
  ellipse were first expressed neatly in general as an infinite
  series, which converges in all cases as far as
  possible. However, although several such series can already be
  found,\footnote{Euler discovered the series he is about to
    give in his paper E448 \cite{E448} which he had written a
    year earlier. \JE} nevertheless hardly would any be discovered
  which would snatch the palm\footnote{i.e. the palm of
    victory. \JE} from the following, which I am about to give.
\end{Scholium}

\newpage

\begin{figure}[htb]
  \centering
  \begin{tikzpicture}
    \draw ({-3*sqrt(2)},0) node [below left] {\(\geom{A}\)} -- (0,0) node [below right] {\(\geom{C}\)} -- (0,{3*sqrt(2)}) node [above right] {\(\geom{D}\)} arc [radius={3*sqrt(2)},start angle=90, end angle=180];
    \draw (0,{2*sqrt(2)}) node [right] {\(\geom{B}\)} arc[x radius={3*sqrt(2)},y radius={2*sqrt(2)},start angle=90, end angle = 180];
    \node (L) at (140:{3*sqrt(2)}) {};
    \node at (L) [above left] {\(\geom{L}\)};
    \node (P) at ({3*sqrt(2)*cos(140)},0) {};
    \node at (P) [below] {\(\geom{P}\)};
    \draw (L.center) -- (P.center);
    \draw (L.center) -- (0,0);
    \node at ({3*sqrt(2)*cos(140)},{2*sqrt(2)*sin(140)}) [below right] {\(\geom{M}\)};
    \node at ([shift=({0,1})]L) {\textit{Fig. 2.}};
    \node at ([shift=({0,1.5})]L) {\textit{Tab. I.}};
  \end{tikzpicture}
\end{figure}

\begin{LemmaUnnum}
  To find a maximally convergent series which produces the
  perimeter of a given ellipse.\footnote{The same series is
    derived in E448 \cite{E448}, using almost the same
    method. the only difference is that Euler performs the
    integral with respect to the variable
    \(z=\frac{x^2}{a^2}-\frac{y^2}{b^2}\) instead of the angle
    \(\phi\). \JE}

  Let \(\geom{A\,M\,B}\) be a quadrant of the proposed ellipse,
  whose centre is at \(\geom{C}\) and semiaxes are
  \(\geom{C\,A}=a\) and \(\geom{C\,B}=b\).  From the centre
  \(\geom{C}\), let the quadrant \(\geom{A\,L\,D}\) of a circle
  be constructed on the semimajor
  axis\marginpar{Tab. I. \\Fig. 2.}  \(\geom{C\,A}\), whose
  radius is therefore \(\geom{A\,C}=a\), and let any radius
  \(\geom{C\,L}\) whatsoever be drawn; then indeed let the
  applicate \(\geom{L\,P}\) be drawn, intersecting
  the ellipse at the point \(\geom{M}\), for which the
  coordinates are called \(\geom{C\,P}=x\) and
  \(\geom{P\,M}=y\); and by putting the angle
  \(\geom{A\,C\,L}=\phi\) we infer \(\geom{C\,P}=x=a.\cos\phi\)
  and \(\geom{P\,L}=a.\sin\phi\). Because indeed\footnote{The
    ellipse is obtained by rescaling the vertical direction by a
    factor of \(b/a\), so this is also the ratio
    \(\geom{P\,M}/\geom{P\,L}\). \JE}
  \(\geom{C\,D}\,:\,\geom{C\,B}=\geom{P\,L}\,:\,\geom{P\,M}=a\,:\,b\),
  we will have \(y=b.\sin\phi\), whence
  \(dx=-a\,d\phi.\sin\phi\) and \(dy=b\,d\phi.\cos\phi\), and
  hence the elliptical line element is
  \[d\phi.\surd\left(a^2\sin\phi^2+b^2\cos\phi^2\right);\]
  wherefore, by integrating, the elliptical arc \(\geom{A\,M}\)
  will be
  \[=\int\!\!\!.\,d\phi.\surd\left(a^2\sin\phi^2+b^2\cos\phi^2\right),\]
  with the integral taken so that it would vanish by putting
  \(\phi=0\). Hence indeed the elliptical quadrant
  \(\geom{A\,M\,B}\) itself will be found by setting
  \(\phi=90\). And in this way, the whole business is reduced to
  a suitable integration of the formula
  \[d\phi\surd\left(a^2\sin\phi^2+b^2\cos\phi^2\right).\]
  But since
  \[\sin\phi^2 = \frac{1-\cos 2\phi}{2}\quad
    \text{and} \quad\cos\phi^2=\frac{1+\cos 2\phi}{2}\] our
  formula to be integrated will transform into this:
  \[\int\!\!\!.\ d\phi.\sqrt{\vphantom{\frac{1}{1}}}
    \left(\frac{a^2+b^2}{2} - \frac{\left(a^2-b^2\right)}{2}
      \cos 2\phi\right)\]
  which we will render more elegantly by putting
  \[a^2+b^2=c^2 \quad\text{and}\quad \frac{a^2-b^2}{a^2+b^2}=n;\]
  then indeed the arc becomes
  \[\geom{A\,M}=\frac{a}{\surd 2}\int\!\!\!.\,d\phi\surd\left(1-n\cos
      2\phi\right)\] on which simple formula the rectification
  of the ellipse depends.  Therefore let us convert this
  irrational formula into a series, which will be
  \begin{align*}
    \surd\left(1\right.&\left.-n.\cos 2\phi\right)=1-\frac{1}{2}.n.\cos 2\phi\\
                       &-\frac{1.\,1}{2.\,4}\,n^2\cos 2\phi^2-\frac{1.\,1.\,3}{2.\,4.\,6}\,n^3\cos 2\phi^3\\
                       &-\frac{1.\,1.\,3.\,5}{2.\,4.\,6.\,8}\,n^4\cos 2\phi^4-\text{etc.}
  \end{align*}
  As to carrying out these integrations, let us call for help
  upon this reduction
  \begin{align*}
    \int\!\!\, d\phi.&\cos 2\phi^{\lambda+2}=\frac{\lambda+1}{\lambda+2}\int\!\!\,d\phi.\cos 2\phi^\lambda \\
                     &+\frac{1}{2\lambda+4}\sin 2\phi.\cos 2\phi^{\lambda+1},
  \end{align*}
  where the final term vanishes by itself in the case \(\phi=0\)
  so that no constant needs adding.\footnote{Recall Euler's
    convention stated earlier that the integrals are definite
    and taken to vanish when \(\phi=0\). \JE} But now for our
  purposes let us extend this integration up to \(\phi=90\), so
  that \(2\phi=180\degree\), and once again that unbound
  term\footnote{{\em terminus ille absolutus} -- Euler means the
    term
    \(\frac{1}{2\lambda+4} \sin 2\phi.\cos
    2\phi^{\lambda+1}\). I have translated {\em absolutus}
    (literally ``finished'' or ``unrestricted'') as ``unbound''
    in the sense that it is no longer inside an integral
    sign. \JE} vanishes; wherefore our reduction becomes
  \[\smallint d\phi.\cos
    2\phi^{\lambda+2}=\frac{\lambda+1}{\lambda+2}\smallint\!\!.\,d\phi.\cos
    2\phi^{\lambda}\] with the help of which, from the first two
  terms of our series, all the following will be most easily
  integrated, namely from the limit \(\phi=0\) up to
  \(\phi=90\degree=\frac{\pi}{2}\), as the following table
  shows.
  \[\arraycolsep=1.4pt\def\arraystretch{1.4}
    \begin{array}{lcc}
    \smallint\!\!.\,d\phi=\phi&=&\frac{\pi}{2}\\
    \smallint\!\!.\,d\phi \cos 2\phi&=&0\\
    \smallint\!\!.\,d\phi. \cos 2\phi^2&=&\frac{1}{2}\cdot\frac{\pi}{2}\\
    \smallint d\phi. \cos 2\phi^3&=&0\\
    \smallint d\phi. \cos 2\phi^4&=&\frac{1}{2}\cdot\frac{3}{4}\cdot\frac{\pi}{2}\\
    \smallint d\phi. \cos 2\phi^5&=&0\\
    \smallint d\phi. \cos 2\phi^6&=&\frac{1.\,3.\,5}{2.\,4.\,6}\cdot\frac{\pi}{2}\\
    \smallint d\phi. \cos 2\phi^7&=&0\\
    \smallint d\phi. \cos 2\phi^8&=&\frac{1.\,3.\,5.\,7}{2.\,4.\,6.\,8}\cdot\frac{\pi}{2}\\
      &\text{etc.}&
    \end{array}
  \]
  Therefore with these values inserted our elliptical quadrant
  will turn out to be
  \begin{align*}
    \geom{A\,M\,B}&=\frac{\pi c}{2\surd 2}\left(1-\frac{1.\,1}{2.\,4}\cdot\frac{1}{2}\,n^2-\frac{1.\,1.\,3.\,5}{2.\,4.\,6.\,8}\cdot\frac{1.\,3}{2.\,4}\,n^4\right.\\
    &\left.\,\,\,\,-\frac{1.\,1.\,3.\,5.\,7.\,9}{2.\,4.\,6.\,8.\,10.\,12}\cdot\frac{1.\,3.\,5}{2.\,4.\,6}\,n^6-\text{etc.}\right).
  \end{align*}
  Because here each numerical coefficient contains in
  itself the preceding, let us write this series in the
  following form
  \[\geom{A\,M\,B}=\frac{\pi c}{2\surd
      2}\left(1-\alpha.n^2-\alpha\beta.n^4-\alpha\beta\gamma.n^6-\alpha\beta\gamma\delta.n^8-\text{etc.}\right)\]
  the values of which letters proceed thus:\footnote{\(\epsilon\) is a typo and should be \(\delta\). \JE}
  \[\arraycolsep=1.4pt\def\arraystretch{1.4}
    \begin{array}{lcccc}
      \alpha&=&\frac{1.\,1.}{2.\,4.}\frac{1}{2}&=&\frac{1.\,1}{4.\,4}\\
      \beta&=&\frac{3.\,5}{6.\,8}\cdot\frac{3}{4}&=&\frac{3.\,5}{8.\,8}\\
      \gamma&=&\frac{7.\,9}{10.\,12}\cdot\frac{5}{6}&=&\frac{7.\,9}{12.\,12}\\
      \varepsilon&=&\frac{11.\,13}{14.\,16}\cdot\frac{7}{8}&=&\frac{11.\,13}{16.\,16}.
    \end{array}\]
  Therefore, for an ellipse whose semiaxes are \(a\) and \(b\), by putting
  \[a^2+b^2=c^2\quad\text{and}\quad \frac{a^2-b^2}{a^2+b^2}=n,\] for the sake of brevity, a quarter of the perimeter will be expressed by the following series 
  \[
    \frac{\pi c}{2\surd 2}\left(1-\frac{1.\,1}{4.\,4}\cdot n^2-\frac{1.\,1.\,3.\,5}{4.\,4.\,8.\,8}\cdot n^4-\frac{1.\,1.\,3.\,5.\,7.\,9}{4.\,4.\,8.\,8.\,12.\,12}\cdot n^6-\text{etc.}\right)
  \]
  which series always converges to a limit, however different
  the axes of the ellipse will have been from one another,
  because \(n\) is always less than unity, and moreover the
  numerical coefficients are strongly decreasing.\footnote{For
    a discussion of convergence, see Appendix \ref{app:convergence}. \JE}
\end{LemmaUnnum}

\begin{CorollaryUnnum}
  The case in which this series converges slowest is where
  \(n=1\), which happens when \(b=0\), or when the conjugate
  axis vanishes; but then it is clear the elliptical quadrant
  will be equal to the semiaxis \(\geom{A\,C}=a\) itself; whence
  because also \(c=a\), we will have for this case\footnote{This
    series for \(\frac{2\surd 2}{\pi}\) was derived in the same
    way by Euler in his earlier paper E448 \cite{E448}. Though
    the infinite product was not derived there, similar series
    and infinite products coming from rectification of an
    ellipse can be found in E154, \S.21 \cite{E154}. This is a
    relation of Wallis's product formula \cite{Wallis},
    which Euler cites in E154. \JE}
  \[\frac{2\surd 2}{\pi}=1 - \frac{1.\,1}{4.\,4} -
    \frac{1.\,1.\,3.\,5}{4.\,4.\,8.\,8} -
    \frac{1.\,1.\,3.\,5.\,7.\,9}{4.\,4.\,8.\,8.\,12.\,12} -
    \text{etc.}\] which series seems in any case worthy of
  attention, and this all the more so because by successively
  collecting terms\footnote{In a little more detail, Euler is
    repeatedly using the fact
    that
    \[1-\frac{1}{(4n)^2} =
      \frac{4n.\,4n-1}{4n.\,4n}=\frac{(4n-1).\,(4n+1)}{4n.\,4n}.\]
    Namely,
    \begin{gather*}
      1-\frac{1.\,1}{4.\,4}=\frac{3.\,5}{4.\,4}\\
      1-\frac{1.\,1}{4.\,4}-\frac{1.\,1.\,3.\,5}{4.\,4.\,8.\,8}
      =\frac{3.\,5}{4.\,4}\left(1-\frac{1}{8.\,8}\right)
      =\frac{3.\,5}{4.\,4}\cdot\frac{7.\,9}{8.\,8},\text{
        etc.}\end{gather*} \JE} it would yield the following
  extraordinary equality:
  \[\frac{3.\,5}{4.\,4}\cdot\frac{7.\,9}{8.\,8}\cdot \frac{11.\,13}{12.\,12}\cdot
    \frac{15.\,17}{16.\,16}\,\text{etc.}
    = \frac{2\surd 2}{\pi}\] and so
  \[\frac{\pi}{2 \surd 2}=\frac{4.\,4}{3.\,5}\cdot\frac{8.\,8}{7.\,9}
    \cdot\frac{12.\,12}{11.\,13}\,\text{etc.}\] the truth of
  which is readily apparent from that which I have published in
  the past on infinite products.
\end{CorollaryUnnum}

\begin{Corollary}[2]
  If we were to set the series obtained earlier as
  \[1-\frac{1.\,1}{4.\,4}\,n^2-\frac{1.\,1.\,3.\,5}{4.\,4.\,8.\,8}\,n^4-\text{etc.}
    = s\] then by that which I have previously
  demonstrated\footnote{In the Opera Omnia edition, the editor
    Andreas Speiser tells us to compare with E52 \cite{E52}. In
    that paper, Euler certainly derives second order
    differential equations for the arc-length of an ellipse (for
    example in \S.17) but the variables used are slightly
    different, and Euler does not discuss the power series
    expansion. I think it is more likely that Euler is referring
    to E448 \cite{E448}; this was written only a year before the
    current paper (E563), and Euler (a) first derived the power
    series \(s\), (b) used it to derive the differential equation
    stated here. The only difference in notation is that Euler
    uses the variable \(v=n/2\). \JE} about the summation of such series, this
  second order differential equation will be obtained:
  \[\frac{4n\,dds}{dn^2}+\frac{4\,ds}{dn}+\frac{ns}{1-nn}=0\]
  where \(dn\) is taken constant.\footnote{Euler essentially
    means ``we differentiate with respect to \(n\)'': if you
    think of \(ds\) as the infinitesimal increment in \(s\)
    caused by the infinitesimal increment \(dn\) in \(n\) then
    you also need to say how much of an increment \(dn\) is. If
    instead you think of \(ds/dn\) as a rate of change then the
    corresponding rate of change of \(n\) is \(dn/dn=1\), a
    constant. The informative Stackexchange answer
    \cite{Bachtold} by Michael B\"{a}chtold discusses this way
    of thinking about differentials, and points to other
    instances of its use by Euler and Leibniz. \JE} Hence
  indeed,
  since \[\frac{4ndds}{dn^2}(1-nn)+\frac{4ds}{dn}(1-nn)+ns=0,\]
  if we were to imagine
  \[s=1+Ann+Bn^4+Cn^6+\text{etc.}\] it will be as
  follows:\footnote{The left-hand side of the first equation
    should be \(\frac{4n\,dds}{dn^2}(1-nn)\). The term
    \(-4.\,2.\,A.\,n^2\) in the expansion of \(4ds/dn\) should be
    \(-4.\,2.\,A.\,n^3\). These typos are corrected in the
    Opera Omnia edition. \JE}
  \[\arraycolsep=1.4pt
    \begin{array}{rcrcrcrcrc}
      \frac{4n\,dds}{dn^2}&=&4.\,2.\,1.\,A.\,n   &+&4.\,4.\,3.\,B.\,n^3&+&4.\,6.\,5.\,C.\,n^5&+&4.\,8.\,7.\,D.\,n^7\,&\text{etc.}\\
                          & &            &-&4.\,2.\,1.\,A.\,n^3&-&4.\,4.\,3.\,B.\,n^5&-&4.\,6.\,5.\,C.\,n^7\,&\text{etc.}\medskip\\
      \frac{4ds}{dn}(1-nn)&=&4.\,2.\,A\,n    &+&4.\,4.\,B.\,n^3  &+&4.\,6.\,C\,n^5 &+&4.\,8.\,D.\,n^7\,&\text{etc.}\\
                          & &            &-&4.\,2.\,A.\,n^2  &-&4.\,4.\,B.\,n^5  &-&4.\,6.\,C.\,n^7\,&\text{etc.}\medskip\\
      ns                  &=&n           &+&A.\,n^3      &+&B\,n^5     &+&C.\,n^7\,    &\text{etc.}
    \end{array}
  \]
  whence arise the following conclusions:
  \begin{align*}
    0&=4.\,4.\,A+1;\text{ therefore }A=\frac{-1}{4.\,4};\\
    0&=4.\,4.\,4.\,B-3.\,5.\,A;\text{ therefore }B=\frac{3.\,5}{8.\,8}.\,A;\\
    0&=4.\,6.\,6.\,C-7.\,9\,.B;\text{ therefore }C=\frac{7.\,9}{12.\,12}.\,B;\\
    0&=4.\,8.\,8.\,D-11.\,13\,.C;\text{ therefore }D=\frac{11.\,13}{16.\,16}.\,C;
  \end{align*}
  whence results the same series which we found above.
\end{Corollary}

\newpage 
\begin{Problem}[2]
  Around a given rectangle \(\geom{M\,m\,N\,n}\), to describe that
  ellipse whose perimeter is minimal.
\end{Problem}
\begin{SolutionUnnum}
  As before, by putting \(\geom{C\,F}=f\), \(\geom{C\,G}=g\) for
  the\marginpar{Tab. I.\\ Fig. 1.}  half-sides of the given rectangle, and
  \(\geom{C\,A}=a\) and \(\geom{C\,B}=b\) for the semiaxes of the desired
  ellipse, we have firstly
  \(\frac{f^2}{a^2}+\frac{g^2}{b^2}=1\). Then indeed, if we were
  to put
  \[a^2+b^2=c^2\quad\text{and}\quad \frac{a^2-b^2}{a^2+b^2}=n,\] we
  find a quarter of the perimeter in the same way as
  before
  \[=\frac{\pi.\,c}{2\surd 2}\left(1-\alpha n^2-\alpha\beta
      n^4-\alpha\beta\gamma n^6-\alpha\beta\gamma\delta
      n^8-\text{etc.}\right),\]
  with it being that
  \[\alpha=\frac{1.\,1}{4.\,4};\, \beta=\frac{3.\,5}{8.\,8};\,
    \gamma=\frac{7.\,9}{12.\,12};\,
    \delta=\frac{11.\,13}{16.\,16};\,\text{etc.}\] which quantity,
  that it might be made minimal, needs its differential to be
  made equal to zero, provided of course the letters \(c\) and
  \(n\) are treated as variables, whence upon dividing by
  \(\frac{\pi}{2\surd 2}\) the following equation will arise:
  \[\begin{array}{r}
      dc\left(1-\alpha\, n\, n-\alpha\beta n^4-\alpha\beta\gamma n^6-\alpha\beta\gamma\delta n^8-\text{etc.}\right)\\
      -c\,dn\left(2\alpha n + 4\alpha\beta n^3+6\alpha\beta\gamma n^5+8\alpha\beta\gamma\delta n^7+\text{etc.}\right)
    \end{array}\biggr\}=0\]
  which may be exhibited in this manner:
  \begin{align*}
    &\frac{dc}{c}\left(1-\alpha n^2-\alpha\beta n^4-\alpha\beta\gamma. n^6-\alpha\beta\gamma\delta n^8-\text{etc.}\right)\\
    =&\frac{dn}{n}\left(2\alpha n^2+4\alpha\beta n^4+6.\alpha\beta\gamma n^6+8.\alpha\beta\gamma\delta.n^8-\text{etc.}\right)
  \end{align*}
  But here the differentials \(dc\) and \(dn\) must maintain a
  fixed relation between them, which it necessary to derive from
  the fundamental equation
  \(\frac{f^2}{a^2}+\frac{g^2}{b^2}=1\). Therefore, first of
  all, it is convenient here for the letters \(c\) and \(n\) to
  be introduced in place of \(a\) and \(b\). Indeed, since
  \[
    \begin{array}{lcrr}
    a^2+b^2=c^2&\quad \text{and}\quad &a^2-b^2=c^2.n,& \text{it will be that}\medskip\\
    a^2=\frac{1}{2}c^2(1+n)&\quad \text{and}\quad&b^2=\frac{1}{2}c^2(1-n);&
    \end{array}
  \]
  whence our equation will transform into this:
  \[\frac{2f^2}{1+n}+\frac{2g^2}{1-n}=c^2.\]
  Because \(f\) and \(g\) are given, let us put for sake of
  brevity\footnote{\(f^2-g^3=h^2.i\) should be
    \(f^2-g^2=h^2.i\). This is corrected in the Opera
    Omnia. \JE}
  \[f^2+g^2=h^2\quad\text{and}\quad
    \frac{f^2-g^2}{f^2+g^2}=i,\quad \text{or}\quad
    f^2-g^3=h^2.\,i;\] having done which, our equation will be
  \(\frac{2h^2(1-in)}{1-nn}=c^2\), of which let us take
  logarithms,\footnote{Euler writes \(l\,x\) where we would write
    \(\log x\). \JE} and it becomes
  \[l\,2h^2+l(1-in)-l(1-nn)=2.\,l\,c\]
  which equation differentiated gives
  \[\frac{-i\,dn}{1-in}+\frac{2n\,dn}{1-nn}=\frac{2\,dc}{c},\]
  whence
  \[\frac{dc}{c}=\frac{n\,dn}{1-nn}-\frac{i\,dn}{2(1-in)}=\frac{2n-i-inn}{2(1-in)(1-nn)}\cdot dn\]
  by substituting which value our equation will be
  \begin{align*}
    &\left(2n^2-in-in^3\right)\left(1-\alpha n^2-\alpha\beta n^4-\alpha\beta\gamma.n^6-\alpha\beta\gamma\delta.n^8\,\text{etc.}\right)\\
    &=2(1-in)(1-nn)\left(2\alpha n^2+4\alpha\beta n^4+6\alpha\beta\gamma n^6\,\,\,\text{etc.}\right)
  \end{align*}
  in which equation only the two quantities \(i\) and \(n\) are
  involved, of which the former is given by the rectangle, but
  the latter \(n\) must be determined by it; therefore that
  cannot otherwise be done except by a solution of the
  infinite equation. It will be convenient therefore to regard
  the quantity \(n\) as known and thence in turn determine
  \(i\), because if it is established for many cases, it will
  easily allow one to determine which value of \(n\) corresponds
  to any given value of \(i\).

  \bigskip
  
  Let us first consider the quantity \(n\) to be minimal, since
  to this will also correspond the minimal value of \(i\);
  therefore by having rejected terms containing \(\beta\),
  \(\gamma\), \(\delta\), etc., we will have
  \[\left(2n^2-in-in^3\right)\left(1-\alpha n^2\right) =
    2\left(1-in\right)\left(1-nn\right).2\alpha n^2\]
  which equation, by rejecting powers of \(n\) greater than the
  third, transforms into this:
  \[2n^2-4\alpha n^2-in-in^3+5\alpha in^3=0,\]
  whence it turns out that
  \[i=\frac{2n(1-2\alpha)}{1+n^2-5\alpha.n^2}.\] Wherefore, if
  \(n\) were an exceedingly small fraction, it becomes very
  nearly
  \[i=2n(1-2\alpha)=\frac{7}{4}n\] whence if, in turn, \(i\)
  were as small a fraction one may conclude it to be that
  \(n=\frac{4}{7}i\).

  \bigskip

  Now, therefore, let us come somewhat closer to the
  truth by rejecting the powers of \(n\) higher than the
  fifth\footnote{The editor of Opera Omnia points out that this
    should be ``the fourth''. \JE} and it will be that
  \[2n^2-4\alpha n^2+2\alpha n^4-8.\alpha\beta n^4-in-in^3+5\alpha in^3=0,\]
  from which is elicited
  \[i=\frac{2n(1-2\alpha)+2\alpha n^3(1-4\beta)}{1+n^2-5\alpha
      n^2},\] and, by dividing\footnote{A quick way to see this
    is to use
    \(\frac{1}{1+n^2-5\alpha n^2}=1-n^2(1-5\alpha)-\cdots\),
    multiply through and drop higher order
    terms. \JE}
  \begin{align*}
    i=2n&(1-2\alpha)+2\alpha
          n^3(1-4\beta)\\&-2n^3(1-2\alpha)(1-5\alpha),
  \end{align*}
  which value, because \(\alpha=\frac{1}{16}\) and
  \(\beta=\frac{15}{64}\), yields
  \[i=\frac{7}{4}\cdot n-\frac{153}{128}\cdot n^3,\] whence in
  turn for a given \(i\) we conclude\footnote{Rather than
    solving the cubic here, Euler is simply giving an
    approximation; namely if we assume that \(n\) can be
    expressed as a power series in \(i\), say
    \(n=Pi+Qi^3+\cdots\) then we can substitute this into
    \(i=\frac{7}{4}\cdot n-\frac{153}{128}\cdot
    n^3=\frac{7}{4}\cdot Pi+\left(\frac{7}{4}\cdot
      Q-\frac{153}{128}\cdot P^3\right) i^3\) and compare terms
    of the same order, giving \(\frac{7}{4}\cdot P=1\) and
    \(\frac{7}{4}\cdot Q=\frac{153}{128}\cdot P^3\), which gives
    the coefficients Euler
    states. \JE} \[n=\frac{4}{7}\,i+\frac{306}{2401}\cdot i^3\] and
  thus one will be able to come closer to the truth.

  \bigskip
  
  Let us additionally develop the case in which \(n\) acquires a
  maximal value, which is equal to unity and for sake of brevity
  let us put
  \begin{align*}
    &1-\alpha-\alpha\beta-\alpha\beta\gamma\,.\,.\,.\,.=s\text{ and}\\
    &2\alpha+4\alpha\beta+6\alpha\beta\gamma\,.\,.\,.\,.=t
  \end{align*}
  so that we would have this equation:
  \begin{align*}
    &2(1-i)s=2(1-i)(1-n\,n)t,\text{ or}\\
    &2(1-i)(s-(1-n\,n)t)=0
  \end{align*}
  whence the first factor manifestly gives \(i=1\), that which the
  nature of the thing demands; indeed if the width of the
  ellipse were to vanish, that is if \(b=0\), then also the
  width of the rectangle must vanish. Therefore, since we will
  have developed both the case in which the letters \(i\) and
  \(n\) are as small as possible, and also that where they
  are allotted the maximum value \(=1\), for the remaining cases,
  judgement will hardly be difficult. Indeed, since for small
  values we have \(n=\frac{4}{7}\,i\), but for the biggest
  \(n=i\); in general we will not deviate too far from the truth
  if we were to set
  \[n=\frac{\frac{4}{7}\,i}{1-\frac{3}{7}\,i\,i} =
    \frac{4i}{7-3.i\,i}.\] Indeed as long as \(i\) is an
  exceedingly small fraction, it will be that
  \(n=\frac{4}{7}i\); and if it were considerably bigger, it
  will be that \(n=\frac{4}{7}\,i+\frac{12}{49}\,i^3\), which
  value is a little bigger than
  \(n=\frac{4}{7}\,i+\frac{306}{2401}\,i^3\); but when \(i=1\)
  it again yields \(n=1\).
\end{SolutionUnnum}

\begin{CorollaryUnnum}
  In order that we may discern more thoroughly the true nature of
  this solution, let us put in general:
  \begin{align*}
    & 1-\alpha\,n^2-\alpha\beta\,n^4-\alpha\beta\gamma\,n^6\,.\,.\,.\,.=s\\
    &2\alpha\,n^2+4\alpha\beta\,n^4+6\alpha\beta\gamma\,n^6\,.\,.\,.\,.=t
  \end{align*}
  so that our equation becomes
  \[\left(2n^2-in-in^3\right)s=2\left(1-in-n^2+in^3\right)t,\]
  whence it follows
  \[i=\frac{2n^2s-2(1-n^2)t}{n(1+n^2)s-2n(1-n^2)t},\] whence, in
  the case \(n=1\), manifestly \(i=1\). But indeed in any other
  case, whatsoever is assumed for \(n\), the true value of \(i\)
  will be able to be determined sufficiently exactly because one
  may sum both series \(s\) and \(t\) readily enough.
\end{CorollaryUnnum}

\begin{Corollary}[2]
  Both series which enter into our solution manifestly depend on
  each other in such a way that \(t=-\frac{n\,ds}{dn}\). But
  now, let us set \(t=sz\), which makes
  \[i=\frac{2n^2-2(1-n^2)z}{n(1+n^2)-2n(1-n^2)z},\] so that here
  a unique value of \(z\) must be found, which will be possible
  to achieve in the following way. Since
  \(t=-\frac{n\,ds}{dn}\), it will now be that
  \(sz=-\frac{n\,ds}{dn}\), and hence
  \(\frac{ds}{s}=-\frac{z\,dn}{n}\), and by
  differentiating\footnote{This differentiation will look more
    familiar to readers with modern sensibilities if we divide
    everything by \(dn^2\):
    \[\frac{1}{s}\frac{d^2s}{dn^2} - \frac{1}{s}\left(\frac{ds}{dn}\right)^2
      = \frac{z}{n^2} - \frac{1}{n}\frac{dz}{dn},\]
    which is the result of differentiating
    \(-\frac{1}{s}\frac{ds}{dn}=\frac{z}{n}\) with respect to
    \(n\) using the quotient rule. \JE}
  \begin{align*}
    &\frac{dds}{s}-\frac{ds^2}{s^2}=+\frac{z\,dn^2}{n^2}-\frac{dz.dn}{n}.\text{
      Let}\\
    &\frac{ds^2}{s^2}=\frac{z^2\,dn^2}{n^2}\text{ be added and
      it makes}\\
    &\frac{dds}{s}=\frac{z\,dn^2}{n^2}+\frac{z^2\,dn^2}{n^2}-\frac{dz\,dn}{n}.
  \end{align*}
  But earlier we gave this equation between \(s\) and \(n\):
  \[\frac{4n.dds}{dn^2}+\frac{4n\,ds}{dn}+\frac{ns}{1-nn}=0,\]
  which, divided by \(s\), yields
  \[\frac{4n}{dn^2}\cdot\frac{dds}{s}+\frac{4}{dn}\cdot\frac{ds}{s}+\frac{n}{1-nn}=0;\]
  but here the newly-ascertained values, having been substituted, produce
  \begin{align*}
    &\frac{4.z}{n}+\frac{4.z^2}{n}-\frac{4dz}{dn}-\frac{4z}{n}+\frac{n}{1-nn}=0,\text{ whence}\\
    &dz=\frac{n\,dn}{4(1-nn)}+\frac{z^2\,dn}{n},
  \end{align*}
  which is a differential equation of only first order, of a
  form similar to that which now for some time is customarily
  referred to as a Riccati equation. But it is always right to
  undertake the solution in such a way that the letter \(n\) be
  considered as known, and, from it, the due value of \(z\) be
  extracted; then indeed from this the value of the letter \(i\)
  is determined. Thus indeed, in turn, we will be able to affirm
  that the assumed value of \(n\) agrees with this same value of
  \(i\); by having found which, since
  \(c^2=\frac{2h^2(1-in)}{1-nn}\), the quantity \(c\) is also
  known; but then finally from \(c\) and \(n\) the semiaxes
  \(a\) and \(b\) are determined.
\end{Corollary}

\begin{Corollary}[3]
  In order that all these calculations may easily be carried
  out, let us present the values of the coefficients \(\alpha\),
  \(\beta\), \(\gamma\), \(\delta\), as decimal fractions, each
  with its logarithm\footnote{In the right-hand column, the
    value \(0.0015808\) appears to have been rounded down from
    \(0.0015808695\) instead of up. In the left-hand column, the
    logarithms are taken with base \(10\) and are to be read as,
    for example, \(8.7958800-10=-1.20412=\log(0.0625)\). There
    appear to be some further errors in the final one or two
    digits of the logarithms, but it is unclear whether these
    are rounding errors, errors inherited from log tables, or
    errors arising from the accuracy of interpolation in
    computing logarithms. Interestingly, the logarithms appear
    to have been computed for the actual rational numbers
    \(\alpha\beta\gamma\cdots\) rather than from their decimal
    expansions; for example \(\log(0.0015808)\) has mantissa
    \(0.198876927\) whilst
    \(\log(\alpha\beta\gamma\delta\epsilon\zeta)\) has mantissa
    \(0.198896033\), much closer to the stated value. Presumably
    this was achieved by computing the logarithms of \(\alpha\),
    \(\beta\), \(\gamma\), etc. and adding them. \JE}
  \[
    \begin{array}{rcrcrcr}
      l\alpha&=&8,7958800;&\qquad&\alpha&=&0,0625000.\\
      l\alpha\beta&=&8,1657913;&\qquad&\alpha\beta&=&0,0146484.\\
      l\alpha\beta\gamma&=&7,8067693;&\qquad&\alpha\beta\gamma&=&0,0064087.\\
      l\alpha\beta\gamma\delta&=&7,5538653;&\qquad&\alpha\beta\gamma\delta&=&0,0035799.\\
      l\alpha\beta\gamma\delta\epsilon&=&7,3583455;&\qquad&\alpha\beta\gamma\delta\epsilon&=&0,0022821.\\
      l\alpha\beta\gamma\delta\epsilon\zeta&=&7,1988959;&\qquad&\alpha\beta\gamma\delta\epsilon\zeta&=&0,0015808.
    \end{array}
  \]
\end{Corollary}

\begin{Example}
  Let us develop the case where the semiaxes of the ellipse
  \(a\) and \(b\) are in a two-fold ratio, or
  \(a\,:\,b=2\,:\,1\),
  whence
  \[n=\frac{a^2-b^2}{a^2+b^2}=\frac{3}{5},\text{ hence
    }n^2=\frac{9}{25},\] whence, for each of the series \(s\)
  and \(t\), the values of each of the terms will be constituted
  thus:\footnote{The editor of the Opera Omnia edition points
    out that the first \(3\) in \(0,0073934\) should be a
    \(5\). In fact, there are several more errors in the table,
    all in the final two digits. The numbers in the right-hand
    column should be: \(0.0450000\), \(0.0075938\),
    \(0.0017940\), \(0.0004810\), \(0.0001380\), \(0.0000413\),
    \(0.0000127\). The final entry in the left-hand column
    (labelled \(-\,-\,-\,-\)) should be \(0.0000009\). \JE}
  \[
    \begin{array}{rcr|rcr}
      \alpha n^2 &=&0,0225000&2\alpha n^2&=&0,0450000\\
      \alpha\beta n^4&=&0,0018984&4\alpha\beta n^4&=&0,0073934\\
      \alpha\beta\gamma n^6&=&0,0002990&6\alpha\beta\gamma n^6&=&0,0017940\\
      \alpha\beta\gamma\delta n^8&=&0,0000601&8\alpha\beta\gamma\delta n^8&=&0,0004808\\
      \alpha\beta\gamma\delta\epsilon n^{10}&=&0,0000138&10\alpha\beta\gamma\delta\epsilon n^{10}&=&0,0001380\\
      \alpha\beta\gamma\delta\epsilon\zeta n^{12}&=&0,0000034&12\alpha\beta\gamma\delta\epsilon\zeta n^{12}&=&0,0000408\\
      -\,-\,-\,-&&10&-\,-\,-\,-&&215\\\hline
                 &&0,0247757&t&=&0,0550685\\
      \text{therefore }s&=&0,9752242&&&
    \end{array}
  \]
  Hence it is obtained that \(\operatorname{log.}z=8,7517967\). Since now
  \[ni=\frac{3}{5}\,i=\frac{2n^2-2(1-n^2)z}{1+n^2-2(1-n^2)z}=\frac{18-32.z}{34-32.z},\]
  we will obtain
  \[\frac{3}{5}\,i=0,50300,\text{ and hence }i=0,838333,\]
  Hence therefore, in turn, if for the given rectangle
  \[i=\frac{f^2-g^2}{f^2+g^2}=0,838333,\] then for the
  corresponding\footnote{{\em satisfaciente} - literally
    ``satisfying'', as in the circumscribed ellipse satisfying
    the condition of having minimal perimeter. \JE} ellipse
  \[n=\frac{3}{5},\quad\text{or}\quad a\,:\,b=2\,:\,1.\]
\end{Example}

\begin{Problem}[3]
  If a type of rectangle were given to which an ellipse of
  minimal perimeter were required to be circumscribed, to determine
  its type by a finite equation.
\end{Problem}

\begin{SolutionUnnum}
  Since the type of rectangle is comprised of the ratio between
  its sides, our letter \(i\) shows its type, since
  \(i=\frac{f^2-g^2}{f^2+g^2}\); then because the type of the
  ellipse is indicated by the ratio between its axes, this will
  be contained in our letter \(n\), since
  \(n=\frac{a^2-b^2}{a^2+b^2}\). An equation expressed in finite
  terms is therefore sought, which would exhibit the relation
  between these quantities \(i\) and \(n\). But now we have seen
  it to be that
  \[in=\frac{2n^2-2(1-n^2)z}{1+n^2-2(1-n^2)z},\] where \(z\) is
  determined by this differential equation:
  \[dz=\frac{n\,dn}{4(1-nn)}+\frac{z^2.dn}{n}.\] Therefore
  nothing else remains except that this letter \(z\) be
  eliminated. So that this can be done more easily, let us put
  \(2(1-n^2)z=x\), which makes \(in=\frac{2n^2-x}{1+n^2-x}\);
  but because \[z=\frac{x}{2(1-nn)}\quad \text{it will be that}\quad
    dz=\frac{dx}{2(1-nn)}+\frac{nx\,dn}{(1-nn)^2};\] whence it
  will arise
  \[2(1-nn)dx+4nx\,dn=(1-nn)n\,dn+\frac{x^2.dn}{n}.\] Indeed
  from this equation, by putting \(v\) in place of \(in\), is
  derived
  \begin{align*}
    &x=\frac{2n^2-(1+n^2)v}{1-v},\quad\text{and hence}\\
    &dx=-\frac{dv(1-nn)}{(1-v)^2}+\frac{4n\,dn-2nv\,dn}{1-v},
  \end{align*}
  by having substituted which, yields
  \[
    \begin{array}{r}
    -\dfrac{2\,dv(1-nn)^2}{(1-v)^2}+8n\,dn\quad =\quad (1-nn)n\,dn\\
      + \dfrac{4n^4-4n^2(1+n^2)v +
      (1+n^2)^2v^2}{(1-v)^2}\cdot\dfrac{dn}{n},
    \end{array}
  \] and hence with
  fractions cleared, it comes to
  \begin{align*}
    -\frac{2\,dv(1-nn)^2}{dn}=-7n+3n^2&+10.nv+\frac{1}{n}v^2\\
    &-2n^3.v-5n.v^2.
  \end{align*}
  Now finally let us substitute in place of \(v\) its value
  \(in\), and it yields
  \begin{align*}
    -\frac{2n(1-nn)^2.\,di}{dn}&=-7.n+3n^3\\
                                &+2i(1+3nn)\\
    &+i^2.n\left(1-5.n^2\right).
  \end{align*}
  But if therefore a {\em constructio}\footnote{According to
    {\cite[pp.12--14]{CapobiancoEneaFerraro}}, {\em constructio} had a
    specific meaning for Euler, namely an expression of the
    solution as an integral which could not be carried out in
    terms of elementary functions, but which would nevertheless
    be useful in obtaining numerical approximations. This is
    distinct from {\em adeo integrare} (to integrate exactly),
    which would mean a solution of the differential equation in
    terms of elementary functions. \JE} of this equation were
  granted, it would be possible to elicit not only for each
  value of \(n\) the corresponding value for \(i\), but also in
  turn for each value of \(i\) the corresponding \(n\). And
  certainly this equation has turned out much simpler than one
  had been permitted to hope for at the outset. But if it were
  permitted to integrate it exactly, or at least to {\em construct}
  it, it would be judged to have added an extraordinary
  development to analysis.
\end{SolutionUnnum}

\vspace{1.5cm}

\begin{center}
  \begin{tikzpicture}
    \draw[very thick] (0,0) -- (8,0);
    \draw (0,-0.1) -- (8,-0.1);
  \end{tikzpicture}
\end{center}

\newpage

\setcounter{eulercounter}{0}
\begin{center}
  {\sc \Huge A geometric problem,\\ \LARGE

    \medskip

    in which amongst all ellipses which can

    \medskip

    be drawn through four given points, that

    \bigskip

    one is sought which has minimal area.
  }\\

  \bigskip

  \bigskip
  \normalsize
  
  An English translation of
  \Large\medskip
  
  {\sc Problem geometricvm, qvo inter omnes ellipses, qvae per
    data qvatvor pvncta tradvci possvnt, ea qvaeritvr, qvae
    habet aream minimam.} (E691)
  
  \medskip

  by Leonhard Euler\\

  \bigskip

  \small
  
  Translated by Jonathan David Evans, School of Mathematical Sciences,\\
  Lancaster University \verb|j.d.evans@lancaster.ac.uk|\\

\end{center}

\begin{center}
  \begin{tikzpicture}
    \draw[very thick] (0,0) -- (10,0);
  \end{tikzpicture}

  Presented to the meeting on the day 4 Sept. 1777.

  \begin{tikzpicture}
    \draw[very thick] (0,0) -- (10,0);
  \end{tikzpicture}
\end{center}

\begin{center}
  \Large \S.1.
\end{center}

\noindent {\Huge T}o the case of this problem in which the four
points are assumed to have been placed at the corners of a
right-angled parallelogram, I have previously given a
solution;\footnote{See Euler's paper E563 \cite{E563}, presented
  to the St. Petersburg Academy on the 15th February 1773, four
  years earlier. \JE} but the general problem I did not venture
to attempt at the time on account of the vast number of
quantities which needed to be introduced in calculations, whence
thoroughly inextricable analytic formulas arose: for which
reason, I hope it will be by no means unpleasant for Geometers
if I will have related here a sufficiently succinct solution of
that most difficult Problem.

\bigskip

\S.\,\,\,2.\quad First, therefore, the four given points must be so
disposed, that at least one ellipse may be drawn through them,
which occurs when any one of those points lies outside the
triangle formed by the remaining three. And indeed as soon as a
single ellipse can be drawn through these points, it is now
clear enough that infinitely many others may also be drawn,
amongst which therefore our problem bids us to seek the one
whose area is smallest of all.

\bigskip
\newpage
\begin{figure}
  \centering
  \begin{tikzpicture}
    \draw (0,0) -- (4.6,0);
    \draw (0,0) -- (2.8,3.1);
    \draw (1,0) -- ++ (2.8*3.8/4.3,3.1*3.8/4.3) node (y) [above right]
    {\(\geom{Y}\)};
    \node at (0,0) [below] {\(\geom{O}\)};
    \node at (1,0) [below] {\(\geom{X}\)};
    \node at (2.6,0) [below] {\(\geom{A}\)};
    \node at (4.6,0) [below] {\(\geom{B}\)};
    \node at (1.5,3.1*1.5/2.8+0.3) {\(\geom{C}\)};
    \node at (2.8,3.1) [above right] {\(\geom{D}\)};
    \node at (4,2) {\textit{Fig. 2.}};
    \node at (0,3.2) {\textit{Tab. I.}};
  \end{tikzpicture}
\end{figure}

\S.\,\,\,3.\quad Therefore let \(\geom{A}\), \(\geom{B}\),
\(\geom{C}\), \(\geom{D}\), be those four
\marginpar{Tab. I.\\Fig. 2.} points through which ellipses are
required to be drawn. Through any two points \(\geom{A}\) and
\(\geom{B}\), let a straight line \(\geom{O\,A\,B}\) be made,
for having as an axis, which the line drawn through the
remaining two points \(\geom{C}\) and \(\geom{D}\) meets at the
point \(\geom{O}\), at which point let us fix the start of the
abscissa. But as applicates let us fix not the normals to the
axis \(\geom{O\,A\,B}\) following the usual practice, but the
parallels to the other direction \(\geom{O\,C\,D}\); namely, if
the abscissa were called \(\geom{O\,X}=x\), the applicate
corresponding to it \(\geom{X\,Y}=y\) is always to be conceived
parallel to the line \(\geom{O\,C\,D}\). Therefore let that
angle of obliquity be called \(\geom{A\,O\,C}=\omega\), and
since the four points \(\geom{A}\), \(\geom{B}\), \(\geom{C}\),
\(\geom{D}\), have been given, let us call their distances from
the point \(\geom{O}\) as follows: \(\geom{O\,A}=a\);
\(\geom{O\,B}=b\); \(\geom{O\,C}=c\) and \(\geom{O\,D}=d\),
whence immediately we will be able to express both the sides and
the diagonals connecting these four points. Indeed first it will
be that \(\geom{A\,B}=b-a\) and \(\geom{C\,D}=d-c\); then indeed
there will be the lines not drawn in the figure:
\begin{align*}
  \geom{A\,C}&=\surd \left(cc+aa-2ac\cos\omega\right),\\
  \geom{A\,D}&=\surd \left(aa+dd-2ad\cos\omega\right),\\
  \geom{B\,C}&=\surd \left(bb+cc-2bc\cos\omega\right),\\
  \geom{B\,D}&=\surd \left(bb+dd-2bd\cos\omega\right).
\end{align*}
Moreover this is observed to be just the same, whichever of the
given points the axis were drawn through, provided that the
direction of the applicates traverses the two remaining points;
it will help to have noted this when, by chance, the lines
\(\geom{A\,B}\) and \(\geom{C\,D}\) will have been parallel to one another;
then indeed it is fit that the axis be drawn through the points
\(\geom{A}\) and \(\geom{D}\) or \(\geom{B}\) and \(\geom{C}\).

\bigskip

\S.\,\,\,4.\quad Now because it is necessary for the curves
leading through the four points \(\geom{A}\), \(\geom{B}\),
\(\geom{C}\), \(\geom{D}\), to be ellipses, let the equation
between the coordinates \(\geom{O\,X}=x\) and \(\geom{X\,Y}=y\)
be represented in this form:
\[Axx+2Bxy+Cyy+2Dx+2Ey+F=0,\]
in which therefore, first the letters \(A\) and \(C\) must be endowed
with the same sign; moreover indeed their product \(AC\) must
be bigger than \(BB\), because otherwise curves satisfying this
equation will be hyperbolae. That we now adapt this general
equation to the fixed purpose, let us first set \(y=0\), whence
the equation transforms into this form: \(Axx+2Dx+F=0\), which
therefore must furnish two points located on the axis, namely
\(\geom{A}\) and \(\geom{B}\), for the former of which \(x=a\), and for the
latter indeed \(x=b\), which must therefore be roots of this equation:
\(Axx+2Dx+F=0\); for this reason let us set
\[Axx+2Dx+F=m(x-a)(x-b),\] whence it becomes
\[A=m;\quad D=-\frac{m(a+b)}{2}\quad \text{and} \quad F=mab.\]

\bigskip

\S.\,\,\,5.\quad Now in a similar way let us put the abscissa \(x=0\),
whence the equation becomes \(Cyy+2Ey+F=0\), of which the roots
must give the points \(\geom{C}\) and \(\geom{D}\), or its roots must be
\(y=c\) and \(y=d\); for this reason, let it be set
\[Cyy+2Ey+F=n(y-c)(y-d),\] whence
\[C=n;\quad E=-\frac{n(c+d)}{2}\quad\text{and}\quad F=ncd.\] But
indeed before we had found \(F=mab\); in order that these values
agree, it is chosen that \(m=cd\) and \(n=ab\), wherefore the
general equation will encompass the four given points if:
\begin{gather*}
  A=cd;\quad C=ab;\quad 2D=-cd(a+b);\\
  2E=-ab(c+d)\quad\text{and}\quad F=abcd,
\end{gather*}
so that now all the letters, except \(B\), are determined. In
this way, therefore, the letter \(B\) is left indeterminate, and
for different values it gives rise to innumerable ellipses
passing through the same four points \(\geom{A}\), \(\geom{B}\),
\(\geom{C}\), \(\geom{D}\), provided it be assumed that
\(BB<AC\). Indeed if it were chosen that \(BB=AC\), the curve
would be a parabola, or an infinitely extended ellipse, of which
therefore the area would also be infinite, for which reason the
proposed question is restricted to minimal area.\footnote{As
  opposed to maximising area. \JE} But still less, if it were
that \(BB>AC\), the curves would be hyperbolas, and so excluded
from our problem.

\bigskip

\S.\,\,\,6.\quad Let us now seek the applicate \(\geom{X\,Y}\). But manifestly, to
any absissa\footnote{This should say \(\geom{O\,X}=x\). \JE} \(\geom{C\,X}=x\)
must correspond the pair of applicates \(\geom{X\,Y}\) and
\(\geom{X\,Y'}\),\marginpar{Tab. I.\\Fig. 3.}  since indeed that
applicate will cut the curve in two points \(\geom{Y}\) and \(\geom{Y'}\),
which applicates will therefore be roots of our general
equation, whose solution will give\footnote{This should be
  \(y=\frac{-E-Bx\pm\surd\left[(E+Bx)^2-ACxx-2CDx-FC\right]}{C}\). \JE}
\[y=-\frac{E-Bx\pm\surd\left[(E+Bx)^2-ACxx-2CDx-FC\right]}{C},\]
of which two values it will give either the applicate \(\geom{X\,Y}\)
or \(\geom{X\,Y'}\), so that:\footnote{Again, in these equations there
  should be no minus sign outside the fraction and \(-E\)
  inside. \JE}
\begin{align*}
  \geom{X\,Y}&=-\frac{E-Bx\pm\surd\left[(E+Bx)^2-ACxx-2CDx-FC\right]}{C}\text{
        and}\\
  \geom{X\,Y'}&=-\frac{E-Bx\pm\surd\left[(E+Bx)^2-ACxx-2CDx-FC\right]}{C}.
\end{align*}

\begin{figure}
  \centering
  \begin{tikzpicture}
    \draw (0,0) -- (4.6,0);
    \draw (0,0) -- (2.8,3.1);
    \draw (1,0) -- ++ (2.8*3.8/4.3,3.1*3.8/4.3) node (y) {\_};
    \node at (y) [above right] {\(\geom{Y}\)};
    \node at (0,0) [below] {\(\geom{O}\)};
    \node at (1,0) [below] {\(\geom{X}\)};
    \node at (2.6,0) [below] {\(\geom{A}\)};
    \node at (4.6,0) [below] {\(\geom{B}\)};
    \node at (1.5,3.1*1.5/2.8+0.3) {\(\geom{C}\)};
    \node at (2.8,3.1) [above right] {\(\geom{D}\)};
    \node at (1.25,2.5) {\textit{Fig. 3.}};
    \node at (0,3.2) {\textit{Tab. I.}};
    \node at (2.1,1.5) {\(\geom{Y'}\)};
    \node at (2.225,1.375) {\_};
    \draw (1.5,0) -- ++ (2.8*1.525/4,3.1*1.525/4) node (y') {\_}
    -- ++ (2.8*1.7/4,3.1*1.7/4) node (y) {\_};
    \node at (y) [below right] {\(\geom{y}\)};
    \node at (y') [right] {\(\geom{y'}\)};
    \node at (1.5,0) [below] {\(\geom{x}\)};
    \draw (1+2.8*0.08,3.1*0.08) node (v) {} -- (1.5,0);
    \node at (1+2.8*0.08-0.1,3.1*0.08+0.1) {\(\geom{v}\)};
    \begin{scope}[shift={(8,1.75)}]
      \draw (0,0) [radius=1.75] circle;
      \draw (0,0) node [left] {\(\geom{a}\)} -- (1.75,0) node [right]
      {\(\geom{r}\)};
      \draw (0,0) -- (40:1.75) node (Y2) {};
      \draw (0,0) -- (-30:0.9);
      \node at (-30:1.1) {\(\geom{t}\)};
      \node at (0.866*0.6-0.25,-0.5*0.6-0.15) {\(\geom{T}\)};
      \draw (-89:1.75) node [below] {\(\geom{y'}\)} -- (29:1.75) node
      [above right] {\(\geom{y}\)};
      \node (Y'2) at (-100:1.75) {};
      \node at (Y2) [above right] {\(\geom{Y}\)};
      \node at (Y'2) [below] {\(\geom{Y'}\)};
      \draw (Y2.center) -- (Y'2.center);
      \node at (-0.5,2.5-1.75) {\textit{Fig. 4.}};
    \end{scope}
  \end{tikzpicture}
\end{figure}

\bigskip

\S.\,\,\,7.\quad Now because both points \(\geom{Y}\) and \(\geom{Y'}\) are
situated on the ellipse passing through the points \(\geom{A}\), \(\geom{B}\),
\(\geom{C}\), \(\geom{D}\), the interval \(\geom{Y\,Y'}\) will be contained in the
ellipse. Wherefore since \(\geom{Y\,Y'}=\geom{X\,Y'}-\geom{X\,Y}\), that interval
will be:
\[\geom{Y\,Y'}=\frac{2\surd\left[(E+Bx)^2-ACxx-2CDx-FC\right]}{C}.\]
But if now to that applicate were drawn a nearby \(\geom{x\,y\,y'}\),
removed from the former by an interval \(\geom{x\,v}\) (namely a
perpendicular \(\geom{x\,v}\) drawn from \(\geom{x}\) to \(\geom{X\,Y}\), which will
be \(\geom{x\,v}=\partial x\,\sin \omega\) because \(\geom{X\,x}=\partial x\)
and the angle \(\geom{x\,X\,v}=\omega\)), if the interval \(\geom{Y\,Y'}\) were
multiplied by that, it would give rise to the element of area
\(\geom{Y\,Y'\,y\,y'}\), which will therefore be
\[=\frac{2\partial
    x\sin\omega}{C}\surd\left[(E+Bx)^2-ACxx-2CDx-FC\right],\]
the integral of which, therefore, extended across the whole
ellipse, will give the whole area of the ellipse which we are
considering.

\bigskip

\S.\,\,\,8.\quad Because the quadrature of an ellipse depends on the
quadrature of a circle, we will most conveniently find this
integral if we refer the matter to a circle. Let us consider
therefore a circle, the radius of which is \(\geom{a\,r}=r\), and so
its area \(=\pi rr\), in which is chosen the analogous element
\(\geom{Y'\,y'\,y\,Y}\), \marginpar{Tab. I.\\Fig. 4.} to which from the
centre \(\geom{a}\) is drawn a normal \(\geom{a\, T}=t\), and it will give
\(\geom{Y\,Y'}=2\surd(rr-tt)\), and so the element of area
\(\geom{Y\,Y'\,y'\,y}=2\partial t\surd(rr-tt)\).Hence, if the integral
were extended over the whole figure, we learn it to be
\(\int 2\partial t\surd(rr-tt)=\pi rr\), whence if we were to
multiply both sides by \(n\), it will give
\[\int 2\partial t\surd(nnrr-nntt)=\pi nrr,\]
and in the same way
\[\int 2m\partial t\surd(nnrr-nntt)=\pi mnrr.\]

\bigskip

\S.\,\,\,9.\quad Now in order for us to adapt this form to our plan, let us
take \(t=x+f\), and it will give
\[\int 2m\partial x\sqrt{\vphantom{\left[nnrr-nn(x+f)^2\right]}}\left[nnrr-nn(x+f)^2\right]=\pi mnrr,\] and
hence
\[\int 2m\partial x\sin\omega\sqrt{\vphantom{\left[nnrr-nn(x+f)^2\right]}}\left[nnrr-nn(x+f)^2\right]
  = \pi mnrr\sin\omega.\] Therefore all that remains is for us
to adapt that form to our case, which is done by taking
\(m=\frac{1}{C}\), then
indeed \[nnrr-nn(x+f)^2=(E+Bx)^2-ACxx-2CDx-FC,\] which equation
will find itself developed thus:
\begin{align*}
  nnrr&-nnff-2nnfx-nnxx\\
  &=EE-FC+2(BE-CD)x+(BB-AC)xx,
\end{align*}
Hence it is clear it must be that:\quad 1\textsuperscript{st}.)
\(nn=AC-BB\), and so \(n=\surd(AC-BB)\);\quad
2\textsuperscript{nd}.) it must be that \(nnf=CD-BE\), and so
\(f=\frac{CD-BE}{AC-BB}\);\quad 3\textsuperscript{rd}.) it is
certainly necessary that \[rr=ff+\frac{EE-FC}{nn},\] whereupon
if we were to substitute the ascertained values, it will yield
\[rr=\frac{(CD-BE)^2}{(AC-BB)^2}+\frac{EE-FC}{AC-BB},\] or
\[rr=\frac{CCDD-2BCDE+ACEE}{(AC-BB)^2}-\frac{CF}{AC-BB}.\] From
these ascertained values the total area of our ellipse must be
\(=\pi mnrr\sin\omega\), whence by having made a substitution the
following expression will be
obtained:
\[\pi\frac{(CDD-2BDE+AEE)\sin\omega}{(AC-BB)^{\frac{3}{2}}}-\frac{\pi
    F\sin\omega}{\surd(AC-BB)},\] which area can now be
exhibited in this fashion:
\[\pi\sin\omega\left(\frac{CDD+AEE-2BDE}{(AC-BB)^{\frac{3}{2}}} -
    \frac{F}{\surd(AC-BB)}\right).\] This expression is
therefore most worthy of note, because with its help the total
area of all ellipses can be sufficiently quickly assigned from
only the equation between its coordinates, whether they are
rectangular or obliquangular. Thus if it were to have the most
well-known equation for an ellipse: \(ffxx+ggyy=ffgg\), in
rectangular coordinates, it will be first that \(\sin\omega=1\);
then indeed \(A=ff\); \(B=0\); \(C=gg\); \(D=0\); \(E=0\);
\(F=-ffgg\), whence the total area of this ellipse will be
\(=\pi fg\) by our rule.

\bigskip

\S.\,\,\,10.\quad Therefore since in this way the areas of all ellipses
passing through the four given points \(\geom{A}\), \(\geom{B}\), \(\geom{C}\),
\(\geom{D}\), are known, all that remains is that the minimum be sought
amongst all these areas. For this reason since all the remaining
letters except \(B\) are determined by the four given points,
since indeed we have taken it to be that \(A=cd\); \(C=ab\);
\(2D=-cd(a+b)\); \(2E=-ab(c+d)\) and \(F=abcd\): the question
reduces to this, that the value of the letter \(B\) be sought
which render the recently ascertained formula smallest of all,
or that, by putting \(CDD+AEE=\Delta\) for sake of brevity, this
formula:
\[\frac{\Delta-2BDE}{(AC-BB)^{\frac{3}{2}}}-\frac{F}{\surd(AC-BB)},\]
be made minimal.

\bigskip

\S.\,\,\,11.\quad Therefore let the letter \(B\) be treated as a variable,
and let the differential of this expression be set equal to
zero, whence the following equation will arise:\footnote{The
  first term should be
  \(-\frac{2DE}{(AC-BB)^{\frac{3}{2}}}\). This fixes itself on
  the next line. This is corrected in the Opera
  Omnia edition. \JE}
\[\frac{2DE}{(AC-BB)^{\frac{3}{2}}} -
  \frac{BF}{(AC-BB)^{\frac{3}{2}}} +
  \frac{3B(\Delta-2BDE)}{(AC-BB)^{\frac{5}{2}}} = 0,\] whose
product with \((AC-BB)^{\frac{5}{2}}\) will produce this
equation:
\[\arraycolsep=1.4pt\left.\begin{array}{rlcc}
    -2ACDE&+\,3CDDB&-\,4DEBB&+\,FB^3\\
          &+\,3AEEB&      &     \\
          &-\,ACFB&       &
  \end{array}\right\}=0.\]

\bigskip

\S.\,\,\,12.\quad Observe therefore the complete solution to the
proposed problem is brought to the solution of a cubic equation;
since this always has real roots, it is certain\footnote{Note
  that Euler really needs the root to lie between \(-\surd AC\)
  and \(\surd AC\), which is not guaranteed by this argument. We
  give an alternative argument in Appendix \ref{app:cubic}
  (Lemma 3). \JE}, howsoever the four points be distributed, one
ellipse can always be assigned passing through these four
points, whose area is smallest of all, for which the equation
between the coordinates \(x\) and \(y\) can be exhibited, if
only in place of \(B\), a root arising from this cubic be
substituted. But if perchance it were to happen that this cubic
were to admit three real roots, just as many solutions will take
place, the nature of which, however, I leave the examining to
others.\footnote{See Appendix \ref{app:cubic} for an examination
  of this. \JE}

\bigskip

\begin{center}
  {\huge APPLICATION}

  \medskip
  
  \Large of this solution to the case in which a minimal ellipse
  circumscribing a given parallelogram is sought.
\end{center}

\bigskip

\S.\,\,\,13.\quad Since here the opposite sides are parallel to
one another, neither of them is suitable to be taken as an axis;
for this reason, let us take one diagonal as axis, and let us
fix the applicates parallel to the other. Therefore let
\marginpar{Tab. I.\\ Fig. 5.}  \(\geom{A\,D\,B\,C}\) be the
proposed parallelogram, whose diagonals \(\geom{A\,B}\) and
\(\geom{C\,D}\) will intersect one another at \(\geom{O}\), and
let it be called \(\geom{A\,O}=\geom{B\,O}=a\) and
\(\geom{C\,O}=\geom{O\,D}=c\), and indeed the angle
\(\geom{A\,O\,C}=\theta\). Let the abscissa \(\geom{O\,X}=x\) be
placed in any position whatsoever, taken on the diagonal
\(\geom{A\,B}\) from the point \(\geom{O}\), and the applicate
\(\geom{X\,Y}=y\) corresponding to it parallel to the other
diagonal \(\geom{C\,D}\), and let the equation expressing the
relation between \(x\) and \(y\) be:
\[Axx+2Bxy+Cyy+2Dx+2Ey+F=0,\] and moreover in
\S.\,9. above we saw the area of the ellipse to be
\[\pi\sin\theta\left(\frac{CDD+AEE-2BDE}{(AC-BB)^{\frac{3}{2}}}
    - \frac{F}{\surd(AC-BB)}\right).\]

\begin{figure}
  \centering
  \begin{tikzpicture}
    \draw (0,0) node (a) {} -- (2,2.3) node (c) {} -- ++ (4.8,0) node (b) {} -- (4.8,0) node (d) {} -- cycle;
    \node at (a) [below] {\(\geom{A}\)};
    \node at (b) [above] {\(\geom{B}\)};
    \node at (c) [above] {\(\geom{C}\)};
    \node at (d) [below] {\(\geom{D}\)};
    \draw (a.center) -- (b.center);
    \draw (c.center) -- (d.center);
    \node at (3.3,1.15) [below] {\(\geom{O}\)};
    \draw (1.5,2.3*1.5/2) node (y) {} -- ++ (2.8*0.75/2,-2.3*0.75/2) node (x) {};
    \node at (1.5+2.8*0.75/2+0.15,2.3*1.5/2-2.3*0.75/2-0.17) {\(\geom{X}\)};
    \node at (1.5-0.15,2.3*1.5/2+0.15) {\(\geom{Y}\)};
    \node at (3.6,2.6) {\textit{Fig. 5.}};
    \node at (0,2.6) {\textit{Tab. I.}};
  \end{tikzpicture}
\end{figure}

\bigskip

\S.\,\,\,14.\quad Therefore let us adapt this general equation to the
proposed case, and first indeed it is manifest that the
applicate \(y\) must vanish at the points \(A\) and \(B\), for
which \(x=+a\) and \(x=-a\), whence arise these two equations:
\begin{align*}
  Aaa+2Da+F&=0\quad\text{and}\\
  Aaa-2Da+F&=0,
\end{align*}
whence it follows \(F=-Aaa\) and \(D=0\). Next by putting
\(x=0\), it must become either \(y=+c\) or \(y=-c\), whence
arise these two equations:
\begin{align*}
  Ccc+2Ec+F&=0\quad\text{and}\\
  Ccc-2Ec+F&=0,
\end{align*}
and hence \(F=-Ccc\) and \(E=0\). Therefore since it must be
that \(Aaa=Ccc\), it will be fit to take \(A=cc\) and \(C=aa\),
so that \(F=-aacc\), and so the equation for our curve will be
\[ccxx+2Bxy+aayy-aacc=0.\]

\bigskip

\S.\,\,\,15.\quad Hence, therefore, the area of this ellipse will be
expressed in this way:
\(\frac{\pi aacc\sin\theta}{\surd(aacc-BB)}\), which is smallest
of all by having taken \(B=0\). Therefore let \(B=0\), and
moreover for this smallest ellipse of all we will have this
equation:
\[ccxx+aayy-aacc=0,\] whose area will be \(=\pi
ac\sin\theta\). Where let it be noted the area of this
parallelogram is \(=2ac\sin\theta\), so that the area of the
ellipse is to the area of the parallelogram as \(\pi\) to \(2\).

\bigskip

\S.\,\,\,16.\quad Therefore it is clear the centre of this
ellipse falls precisely at the point \(\geom{O}\), and moreover
both diagonals \(\geom{A\,B}\) and \(\geom{C\,D}\) will be
conjugate\footnote{The fact that \(\geom{A\,B}\) and
  \(\geom{C\,D}\) are conjugate diagonals follows from the
  equation \(aaxx+ccyy=aacc\) for the ellipse since these
  diagonals are the \(x\) and \(y\) axes. \JE} diameters of it,
inclined to one another at an angle of obliquity
\(\geom{A\,O\,C}=\theta\); from which it follows the tangents at
the points \(\geom{A}\) and \(\geom{B}\) are parallel to the
diameter \(\geom{C\,D}\), and indeed the tangents at the points
\(\geom{C}\) and \(\geom{D}\) are parallel to the diameter
\(\geom{A\,B}\), whence this curve is easily described. But if
the angle \(\theta\) were a right-angle, the parallelogram
transforms into a rhombus, whose diagonals \(\geom{A\,B}\) and
\(\geom{C\,D}\) will be principal axes of the ellipse.

\bigskip

\S.\,\,\,17.\quad But if both diagonals \(\geom{A\,B}\) and
\(\geom{C\,D}\) were equal to one another, the angle \(\theta\)
staying oblique, our parallelogram becomes a rectangle; and this
case I already contemplated some time ago, and I determined the
minimal ellipse circumscribing such a rectangle, which solution
moreover agrees excellently with the present one.\footnote{Euler
  is referring to his paper E563 {\cite[\S.1]{E563}}. \JE}

\bigskip

\begin{figure}
  \centering
  \begin{tikzpicture}
    \draw (0,0) node [below] {\(\geom{O}\)} -- (4.7,0) node [below] {\(\geom{A}\)};
    \draw (0,0) -- (4.7,2.25) node [right] {\(\geom{F}\)};
    \draw (2.15,0) node [below] {\(\geom{X}\)} -- (2.9,2.3) node [above] {\(\geom{Y}\)};
    \draw (2.8,0) node [below] {\(\geom{t}\)} -- (3.3,1.58) node [below right] {\(\geom{x}\)};
    \draw (2.67,1.58) node [left] {\(\geom{u}\)} -- (3.3,1.58);
    \draw (3.3,1.58) -- (2.9,2.3);
    \node at (1,1.58) {Fig. 1.};
    \node at (0,2.3) {Tab. II.};
  \end{tikzpicture}
\end{figure}

\S.\,\,\,18.\quad Now let us see also in what way the principal
axes of the \marginpar{Tab. II. \\Fig. 1.} ascertained ellipse
are required to be determined in general. By putting coordinates
\(\geom{O\,X}=x\) and \(\geom{X\,Y}=y\) therefore, with the
angle being \(\geom{A\,X\,Y}=\theta\), we find this equation:
\[ccxx+aayy-aacc=0.\] Now let us put \(\geom{O\,F}\) to be the
principal semi-axis of this ellipse,\footnote{Since the
  coordinates are not orthogonal, \(\geom{O\,A}\) is not the
  semi-axis despite the familiar equation. \JE} inclined to the
line \(\geom{O\,A}\) by an angle \(\geom{A\,O\,F}=\phi\), and
let us refer a point \(\geom{Y}\) of the ellipse to this axis
\(\geom{O\,F}\) in orthogonal coordinates \(\geom{O\,x}=X\) and
\(\geom{x\,Y}=Y\).\footnote{Euler's habit of using the same
  letter to denote to different things makes this passage
  particularly confusing. Recall that with our font conventions,
  \(\geom{X}\) is a point and \(X\) a variable, and we have
  \(\geom{O\,X}=x\) and \(\geom{O\,x}=X\). \JE} To which end let
us draw parallels \(\geom{x\,u}\) and \(\geom{x\,t}\) to the
former coordinates from \(\geom{x}\), and in the triangle
\(\geom{O\,x\,t}\) will be the angle
\(\geom{O\,x\,t}=\theta-\phi\); but indeed in the triangle
\(\geom{x\,u\,Y}\), since the angle \(\geom{O\,x\,u}=\phi\), the
angle will be \(\geom{u\,x\,Y}=90\degree-\phi\), and the angle
\(\geom{x\,Y\,u}\) complementary to the angle \(\theta-\phi\),
and then indeed the angle \(\geom{x\,u\,Y}=\theta\).

\bigskip

\S.\,\,\,19.\quad Now the solution of this triangle
yields:\footnote{The equation for \(\geom{t\,x}\) should have a
  factor of \(X\) in the numerator where the original has left
  space for it to be printed; this factor is added in the Opera
  Omnia edition, and fixes itself a few lines later in the
  original. \JE}
\[\geom{O\,t}=\frac{X\sin(\theta-\phi)}{\sin\theta}\quad\text{and}\quad
  \geom{t\,x}=\frac{\,\,\,\,\sin\phi}{\sin\theta};\]
then indeed
\[\geom{x\,u}=\frac{Y\cos(\theta-\phi)}{\sin\theta}\quad\text{and}\quad
  \geom{Y\,u}=\frac{Y\cos\phi}{\sin\theta},\] whence the former
coordinates \(x\) and \(y\) are determined by \(X\) and \(Y\),
so that
\[x=\frac{X\sin(\theta-\phi)-Y\cos(\theta-\phi)}{\sin\theta}
  \quad\text{and}\quad
  y=\frac{X\sin\phi+Y\cos\phi}{\sin\theta},\]
which values, substituted in the equation:
\[ccxx+aayy=aacc,\] produce this equation between \(X\) and
\(Y\):
\[
  \left.
    \begin{array}{r}
      ccXX\sin(\theta-\phi)^2-2ccXY\sin(\theta-\phi)\cos(\theta-\phi)\\
                             +ccYY\cos(\theta-\phi)^2\\
      +aaXX\sin\phi^2-2aaXY\sin\phi\cos\phi\\ +aaYY\cos\phi^2
    \end{array}
  \right\}
  =aacc\sin\theta^2.
\]
Therefore in this equation, since it is referred to the
principal axes, first of all the terms containing \(XY\) must
cancel with one another, whence
\[cc\sin(\theta-\phi)\cos(\theta-\phi)=aa\sin\phi\cos\phi,\]
from which equation the angle \(\phi\) may be extracted. Indeed
since
\[aa\sin 2\phi=cc\sin(2\theta-2\phi)=cc\sin 2\theta\cos
  2\phi-cc\cos 2\theta\sin 2\phi,\]
upon dividing by \(\sin 2\phi\) we will have:
\[aa=cc\sin 2\theta\cos 2\phi-cc\cos 2\theta,\] and hence it
becomes\footnote{There is a missing equality sign after
  \(\cos 2\phi\), fixed in the Opera Omnia edition. There is
  also a more serious mistake in this equation: the left-hand
  side should be \(\cot 2\phi\), as pointed out in the Opera
  Omnia edition. This mistake propagates, but cancels with a
  second mistake later. \JE}
\[\cos 2\phi\,\frac{aa+cc\cos 2\theta}{cc\sin 2\theta},\]
whence a twofold value for the angle \(\phi\) is elicited, for
the position of either principal axis.

\bigskip

\S.\,\,\,20.\quad Now by having removed the \(X\,Y\) term our
equation will be
\[
  \left\{
    \begin{array}{r}
      XX[cc\sin(\theta-\phi)^2+aa\sin\phi^2]\\
      +YY[cc\cos(\theta-\phi)^2+aa\cos\phi^2]
    \end{array}
  \right\} =aacc\sin\theta^2,\] whence both principal semiaxes, which will
be \(f\) and \(g\), are determined in the following way:
\[ff=\frac{aacc\sin\theta^2}{cc\sin(\theta-\phi)^2+aa\sin\phi^2}
  \quad\text{and}\quad
  gg=\frac{aacc\sin\theta^2}{cc\cos(\theta-\phi)^2+aa\cos\phi^2},\]
and thus
\begin{align*}
  \frac{aacc\sin\theta^2}{ff}&=cc\sin(\theta-\phi)^2+aa\sin\phi^2\text{
                               and}\\
  \frac{aacc\sin\theta^2}{gg}&=cc\cos(\theta-\phi)^2+aa\cos\phi^2,
\end{align*}
whence because the angle \(\phi\) is already found, both
principal semiaxes \(f\) and \(g\) will be able to be determined.

\bigskip

\S.\,\,\,21.\quad If the two last inequalities be added, this
equation will arise:
\begin{align*}
  &\frac{aacc\sin\theta^2(ff+gg)}{ffgg}=cc+aa,\text{ or}\\
  &\frac{ff+gg}{ffgg}=\frac{aa+cc}{aacc\sin\theta^2}.
\end{align*}
Then indeed if in the former equation
\[\frac{aacc\sin\theta^2}{ff}=cc\sin(\theta-\phi)^2+aa\sin\phi^2,\]
on the right hand side, in place of \(cc\) were written the
value
\(\frac{aa\sin\phi\cos\phi}{\sin(\theta-\phi)\cos(\theta-\phi)}\),
it will yield this equation:
\begin{align*}
  \frac{cc\sin\theta^2}{ff}&=\frac{\sin\phi\cos\phi\sin(\theta-\phi)}{\cos(\theta-\phi)}+\sin\phi^2\\
                           &=\frac{\sin\phi[\cos\phi\sin(\theta-\phi)+\sin\phi\cos(\theta-\phi)]}{\cos(\theta-\phi)}=\frac{\sin\phi\sin\theta}{\cos(\theta-\phi)},
\end{align*}
and thus
\(\frac{cc\sin\theta}{ff}=\frac{\sin\phi}{\cos(\theta-\phi)}\). Indeed
then if in the other equation
\[\frac{aacc\sin\theta^2}{gg}=cc\cos(\theta-\phi)+aa\cos\phi^2.\]
on the right-hand side, in place of \(aa\) were written the
value
\(\frac{cc\sin(\theta-\phi)\cos(\theta-\phi)}{\sin\phi\cos\phi}\),
it will yield this equation:\footnote{The denominator in the
  final expression should be \(\sin\phi\). \JE}
\begin{align*}
  \frac{aa\sin\theta^2}{gg}&=\cos(\theta-\phi)^2+\frac{\sin(\theta-\phi)\cos(\theta-\phi)\cos\phi}{\sin\phi}\\
                           &=\frac{\cos(\theta-\phi)}{\sin\phi}[\sin\phi\cos(\theta-\phi)+\cos\phi\sin(\theta-\phi)]=\frac{\sin\theta\cos(\theta-\phi)}{\OP{in}.\phi}.
\end{align*}
whence
\(\frac{aa\sin\theta}{gg}=\frac{\cos(\theta-\phi)}{\sin\phi}\).

\bigskip

\S.\,\,\,22.\quad Now the two latter equalities, multiplied together, will
give \(\frac{aacc\sin\theta^2}{ffgg}=1\), and so
\[ffgg=aacc\sin\theta^2,\] consequently \(fg=ac\sin\theta\), in
which equation is contained that notable property by which the
parallelogram described around two conjugate diameters is
ascribed area\footnote{{\em aequale perhibetur} -- Though he
  does not use the word explicitly, Euler is talking about the
  parallelograms having the same {\em area}: the area of the
  parallelogram described around the diameters is
  \(2ac\sin\theta\) whilst the area of the parallelogram
  described around the principal semi-axes is \(2fg\). \JE}
equal to the parallelogram described around the principal
axes. Then since we found above
\(\frac{ff+gg}{ffgg}=\frac{aa+cc}{aacc\sin\theta^2}\), since we
recently saw it to be that \(aa+cc\sin\theta^2=ffgg\), hence the
other principal property follows, that is \(aa+cc=ff+gg\),
namely in every ellipse the sum of the squares of two diameters
is always equal to the sum of the squares of the principal
axes.\footnote{More precisely, the sum of squares of two
  {\em conjugate} diameters of an ellipse equals the sum of
  squares of the principal axes. \JE}

\bigskip

\S.\,\,\,23.\quad Let us apply this to the case of a rectangle, already
formerly worked out, for which \(c=a\), and for the position of
the principal axes it will now have the equation:\footnote{Here
  again, the Opera Omnia points out that the left-hand side
  should be \(\cot 2\phi\). \JE}
\[\cos 2\phi=\frac{1+\cos 2\theta}{\sin 2\theta},\]
whence it is inferred\footnote{The next equation is
  incorrect. Instead, as the Opera Omnia points out, it should
  be
  \[\cot
    2\phi=\frac{1+\OP{cos}2\theta}{\OP{sin}2\theta} =
    \frac{1+\OP{cos}^2\theta -
      \OP{sin}^2\theta}{2\sin\theta\cos\theta} =
    \frac{\cos\theta}{\sin\theta}=\cot\theta.\] This still
  yields \(\phi=\theta/2\) or \(\phi=90\degree+\theta/2\), so
  the mistakes cancel. \JE}
\[\cos 2\phi=\frac{1+\cos 2\theta}{\surd 2(1+\cos
    2\theta)} = \sqrt{\vphantom{\frac{1+\cos 2\theta}{2}}}
  \frac{1+\cos 2\theta}{2}.\]
But it is evident that
\[\sqrt{\vphantom{\frac{1+\cos 2\theta}{2}}}
  \frac{1+\cos 2\theta}{2}=\pm\cos\theta,\]
\newpage
\begin{figure}
  \centering
  \begin{tikzpicture}
    \draw (-2.5,0) node [below] {\(\geom{B}\)} -- (2.6,0) node [below] {\(\geom{A}\)};
    \draw (252:1) node [left] {\(\geom{D}\)} -- (72:2.75) node [above] {\(\geom{C}\)};
    \draw (0,0) -- (36:3.4);
    \node at (36:3.65) {\(\geom{F}\)};
    \draw (0,0) -- (36+90:2.4);
    \node at (36+90:2.65) {\(\geom{G}\)};
    \node at (0.15,-0.2) {\(\geom{O}\)};
    \node at (-0.6,3) {\textit{Tab. II.}};
    \node at (-0.6,2.5) {\textit{Fig. 2.}};
  \end{tikzpicture}
\end{figure}

\noindent whence either \(2\phi=\theta\), and so
\(\phi=\frac{1}{2}\theta\); or \(2\phi=\pi+\theta\), and so
\(\phi=90\degree+\frac{1}{2}\theta\). Hence therefore it is\marginpar{Tab. II\\Fig. 2.}
clear that on the one hand the principal axis \(\geom{O\,F}\)
bisects the angle \(\geom{A\,O\,C}=\theta\), and on the other,
\(\geom{O\,G}\), normal to this, bisects the angle
\(\geom{B\,O\,C}\). Then indeed it will be that
\[fg=aa\sin\theta\text{ and }ff+gg=2aa,\] whence we infer
\[(f+g)^2=2aa(1+\sin\theta) =
  4aa\cos\left(45\degree-\tfrac{1}{2}\theta\right)^2,\] and so
it will be that
\[f+g=2a\cos\left(45\degree-\tfrac{1}{2}\theta\right).\] In a
similar way will be obtained
\[(f-g)^2=2aa(1-\sin\theta) =
  4aa\sin\left(45\degree-\tfrac{1}{2}\theta\right)^2,\]
consequently
\[f-g=2a\sin\left(45\degree-\tfrac{1}{2}\theta\right),\]
wherefore we will have
\[f=a\left[\cos\left(45\degree-\tfrac{1}{2}\theta\right) +
    \sin\left(45\degree-\tfrac{1}{2}\theta\right)\right] =
  a\cos\tfrac{1}{2}\theta.\surd 2,\]
and in a similar way
\[g=a\left[\cos\left(45\degree-\tfrac{1}{2}\theta\right) -
    \sin\left(45\degree-\tfrac{1}{2}\theta\right)\right] =
  a\sin\tfrac{1}{2}\theta.\surd 2,\]
which values manifestly suffice; indeed
\[fg=aa\sin\theta\quad\text{and}\quad ff+gg=2aa,\] and this
solution agrees perfectly with the one I gave
before.\footnote{In his earlier paper \cite{E563}, Euler found
  that the ratio of the semi-axis to the parallel half-side of
  the rectangle is \(\surd 2\). This coincides with what is
  written here, since the half-sides of the rectangle are
  \(a\cos\tfrac{1}{2}\theta\) and
  \(a\sin\tfrac{1}{2}\theta\). Note that Euler has reversed his
  notation from \cite{E563}, where he uses \(f\) and \(g\) for
  the half-sides and \(a\) and \(b\) for the semiaxes. \JE}

\vspace{1.5cm}

\begin{center}
  \begin{tikzpicture}
    \draw[very thick] (0,0) -- (8,0);
    \draw (0,-0.1) -- (8,-0.1);
  \end{tikzpicture}
\end{center}

\newpage
\setcounter{eulercounter}{0}

\begin{center}
  {\sc \Large Solution to a

    \huge \medskip

    problem of greatest curiosity

    \medskip \Large

    in which, amongst all ellipses which can be circumscribed
    around a given triangle, that one is sought, whose area is
    smallest of all

  }

  \bigskip

  \bigskip
  \normalsize
  
  An English translation of
  \Large\medskip
  
  {\sc Solvtio problematis maxime cvriosi, qvo inter omnes
    ellipses, qvae circa datvm triangvlvm circvmscribi possvnt,
    ea qvaeritvr, cvivs area sit omnivm minima.} (E692)
  
  \medskip

  by Leonhard Euler\\

  \bigskip

  \small
  
  Translated by Jonathan David Evans, School of Mathematical Sciences,\\
  Lancaster University \verb|j.d.evans@lancaster.ac.uk|\\

\end{center}

\begin{center}
  \begin{tikzpicture}
    \draw[very thick] (0,0) -- (10,0);
  \end{tikzpicture}

  Presented to the meeting on the day 4 Sept. 1777.

  \begin{tikzpicture}
    \draw[very thick] (0,0) -- (10,0);
  \end{tikzpicture}
\end{center}

\begin{center}
  \Large \S.1.
\end{center}

\noindent {\Huge A}fter I had already attempted this problem
before in many ways in vain, recently at last I hit upon an
truly remarkable method of investigating its solution, which is
rather worthy of note, because it leads to an exceedingly simple
and easy construction. Namely, I used the same method by which I
recently demonstrated how, amongst all ellipses which may be
drawn through four given points, to assign that which has
minimal area, whence I am going to solicit the particular help
of a calculation from that paper.

\bigskip

\S.\,\,\,2.\quad Therefore let \(\geom{A\,B\,C}\) be the
proposed\marginpar{Tab. II\\Fig. 3}
triangle, whose angle at \(\geom{B}\) let us call \(=\omega\), and let
us call the two sides forming this angle \(\geom{B\,A}=a\) and
\(\geom{B\,C}=c\), so that the third side is:
\[\geom{A\,C}=\surd(aa+cc-2ac\cos\omega);\] moreover indeed it will
help to have noted the area of this triangle to be
\(=\frac{1}{2}ac\sin\omega\). But we will see the area of the
minimal ellipse passing through the three points \(\geom{A}\), \(\geom{B}\),
\(\geom{C}\), always holds a fixed ratio to the area of this triangle,
which ratio will be found namely as \(4\pi\) to \(3\surd 3\).

\bigskip
\begin{figure}
  \centering
  \begin{tikzpicture}
    \draw (0,0) node (b) [below] {\(\geom{B}\)} -- (4.6,0) node (a) {};
    \node at (a) [below] {\(\geom{A}\)};
    \draw (0,0) -- (1.1,2.75) node (c) {};
    \node at (c) [above] {\(\geom{C}\)};
    \draw (1.1,0) node (x) [below] {\(\geom{X}\)} -- ++ (1.1*1.1,2.75*1.1) node (y) [right]{\(\geom{Y}\)};
    \draw (a.center) -- (c.center);
    \node at (2.8,2.15) {\textit{Fig. 3.}};
    \node at (-0.5,3.2) {\textit{Tab. II.}};
  \end{tikzpicture}
\end{figure}

\S.\,\,\,3.\quad Let \(\geom{Y}\) be any point whatsoever in the
desired ellipse, whose position let us determine through two
obliquangular coordinates parallel to the two axes
\(\geom{B\,A}\) and \(\geom{B\,C}\); for this reason, a line
\(\geom{X\,Y}\) having been drawn parallel to the side
\(\geom{B\,C}\), let us call these coordinates \(\geom{B\,X}=x\)
and \(\geom{X\,Y}=y\), of which let the relation be expressed as
this general second order equation:
\[Axx+2Bxy+Cyy+2Dx+2Ey+F=0,\] and moreover in the aforementioned
paper where I considered four points, I showed the total area of
this ellipse to be
\[=\pi\sin\omega\left(\frac{CDD+AEE-2BDE}{(AC-BB)^{\frac{3}{2}}}
    - \frac{F}{\surd(AC-BB)}\right),\] where following the
usual practice \(\pi\) denotes the circumference of a circle
whose diameter \(=1\).

\bigskip

\S.\,\,\,4.\quad First of all therefore let us adapt this
general form to the situation of our triangle. And first indeed
since by having taken \(x=0\) for the point \(\geom{B}\) itself
and moreover it be made \(y=0\), it is manifest that \(F=0\)
must be fixed, whence the area above will already be expressed
more simply. Then since \(x=a\) and \(y=0\) for the point
\(\geom{A}\), the general equation will give \(Aaa+2Da=0\),
whence \(2D=-Aa\). Thirdly indeed for the point \(\geom{C}\) it
will be that \(x=0\) and \(y=c\), whence \(Ccc+2Ec=0\), and so
\(2E=-Cc\). Therefore the general equation, having been adapted
to the proposed situation, will be
\[Axx+2Bxy+Cyy-Aax-Ccy=0,\] which therefore completely
encompasses all ellipses which may be drawn through the three
given points \(\geom{A}\), \(\geom{B}\), \(\geom{C}\); in which
equation therefore three indefinite letters are still involved,
namely \(A\), \(B\) and \(C\).

\bigskip

\S.\,\,\,5.\quad Therefore these letters must be determined so
that the area of the ellipse be rendered smallest of
all. Therefore, since \(F=0\); \(D=-\frac{1}{2}Aa\) and
\(E=-\frac{1}{2}Cc\), the area of the ellipse from the general
formula given above will be expressed thus:
\[\frac{1}{4}\pi\sin\omega\left(\frac{AACaa+ACCcc-2ABCac}{(AC -
      BB)^{\frac{3}{2}}}\right),\] where therefore the letters
\(A\), \(B\), \(C\) must be determined so that this equation
become smallest of all, whence it is clear that two
investigations of the minimum must be made.\footnote{In other
  words, since the area depends on the triple \(A,B,C\) of
  coefficients only up to scale, there are really only two free
  variables, and Euler is saying that we need to minimise with
  respect to both. \JE}

\bigskip

\S.\,\,\,6.\quad In order that we adapt this formula more
closely to the calculation, let us first set\footnote{Note that
  in order for the conic to be an ellipse we must have
  \(B^2<AC\), so such a \(\phi\) exists. \JE}
\(B=\cos\phi\surd AC\), so that
\begin{align*}
  &AC-BB=AC\sin\phi^2,\quad \text{and so}\\
  &(AC-BB)^{\frac{3}{2}}=AC\sin\phi^3\surd AC,
\end{align*}
by having introduced which value, the area of the ellipse will become
\[=\frac{1}{4}\pi\sin\omega\left(\frac{Aaa+Ccc-2ac\cos\phi\surd
      AC}{\sin\phi^3\surd AC}\right).\] Now to completely cancel
the irrational\footnote{I believe that Euler is referring to
  \(\surd AC\) as {\em the irrational}. \JE}, let \(C=Ass\), so
that \(\surd AC=As\), and in this way our area will be
\[\frac{1}{4}\pi\sin\omega\left(\frac{aa+ccss-2acs\cos\phi}{s\sin\phi^3}\right),\]
and the question now reduces to how the quantity \(s\) and the
angle \(\phi\) are required to be determined, so that the value
of this expression: \(\frac{aa+ccss-2acs\cos\phi}{s\sin\phi^3}\)
be rendered smallest of all?

\bigskip

\S.\,\,\,7.\quad Let us suppose the angle \(\phi\) is already
given a due value, so that only the quantity \(s\) need be
investigated, by which the minimal value is procured in that
formula; in which investigation therefore the angle \(\phi\)
will be considered as constant, but only \(s\) as variable, and
thus this expression will need to be rendered minimal:
\[\frac{aa+ccss-2acs\cos\phi}{s}=\frac{aa+ccss}{s}-2ac\cos\phi,\]
of which the last part is now constant, whence this formula
alone: \(\frac{aa}{s}+ccs\) must be brought to a minimum, with
the differential of which having been set equal to zero, it
yields this equation: \(ccss-aa=0\), whence it is inferred
\(s=\frac{a}{c}\). But it was that \(ss=\frac{C}{A}\), and so
\(\frac{C}{A}=\frac{aa}{cc}\) must be taken. Since therefore in
our equation only the ratio between the letters \(A\) and \(C\)
is considered, let us take \(A=cc\) and \(C=aa\), and in this
way we have already fulfilled one condition of the minimum.

\bigskip

\S.\,\,\,8.\quad In place of \(A\) and \(C\) let us write those
ascertained values, and the area of the ellipse, already partly
rendered minimal, will be
\[\frac{1}{4}\pi\sin\omega\left(\frac{2ac(1-\cos\phi)}{\sin\phi^3}\right)=\frac{1}{2}\pi
  ac\left(\frac{1-\cos\phi}{\sin\phi^3}\right).\] Therefore it
only remains for the angle \(\phi\) to be determined, so that
the formula \(\frac{1-\cos\phi}{\sin\phi^3}\) comes out
minimal. But since \(\sin\phi^2=1-\cos\phi^2\), that formula
\(\frac{1-\cos\phi}{\sin\phi^3}\) is transmuted into this:
\(\frac{1}{\sin\phi(1+\cos\phi)}\), of which fraction therefore
the denominator ought to be rendered maximal; but its
differential gives this equation
\(\cos\phi+\cos\phi^2-\sin\phi^2=0\), or
\(\cos\phi+2\cos\phi^2-1=0\), which is manifestly resolved into
these factors: \((1+\cos\phi)(2\cos\phi-1)=0\), whence two
solutions follow: on the one hand \(1+\cos\phi=0\) which however
would yield the formula \(\sin\phi(1+\cos\phi)=0\), and so there
will not be a maximum;\footnote{Recall that
  \(\sin\phi(1+\cos\phi)\) is supposed to be the denominator of
  the area, so cannot vanish. \JE} wherefore the other solution
will have its place, which gives \(2\cos\phi-1=0\), whence
\(\cos\phi=\frac{1}{2}\), and hence
\(\sin\phi=\frac{\surd 3}{2}\), namely the angle itself will
be\footnote{Here Euler uses the abbreviation {\em gr.} for {\em
    gradus}, meaning degrees. Elsewhere, for example in E691,
  \S.23, he uses the modern symbol \(\degree\), which came into
  use in the late sixteenth century according to Cajori
  {\cite[\S.514]{CajoriNotation}}. Cajori points out that the
  {\em gr.} notation survived well into the late eighteenth
  century. \JE} \(=60\text{ deg.}\)

\bigskip

\S.\,\,\,9.\quad In view of this therefore we have completely
satisfied the prescribed condition, and now the area of the
ellipse will be expressed in this manner:
\(\frac{2\pi ac\sin\omega}{3\surd 3}\), which is minimal amongst
all ellipses which may be drawn through the three given
points. Therefore since the area of the triangle
\(\geom{A\,B\,C}\) is \(\frac{1}{2}ac\sin\omega\), it is evident
that the area of the desired minimal ellipse is to the area of
the triangle as \(4\pi\,:\,3\surd 3\), exactly as we mentioned
above. But this ratio is close in value to \(2,41840\,:\,1\),
whence the following fractions successively approach the truth
more closely:\footnote{This is just the sequence of convergents
  for the continued fraction expansion of \(4\pi/3\surd
  3\). \JE}
\[\frac{2}{1};\quad \frac{5}{2};\quad \frac{12}{5};\quad
  \frac{17}{7};\quad \frac{29}{12};\quad
  \frac{104}{43};\quad\frac{237}{98}.\]

\bigskip

\S.\,\,\,10.\quad Let us now seek also the very equation of the
ascertained curves, and because we have taken \(A=cc\) and
\(C=aa\), hence we will find the letter \(B\) from setting
\(B=\cos\phi\surd AC\), whence because \(\cos\phi=\frac{1}{2}\)
it will be that \(B=\frac{1}{2}ac\), with which
value substituted, the equation for the smallest ellipse of all will be:
\[ccxx+acxy+aayy-accx-aacy=0,\] whence for whichever abscissa
\(x\) it is possible to determine twin applicates \(y\), indeed
it will be found:\footnote{There is a factor of \(c\) missing:
  it should be \(y=\frac{ac-cx\pm
    c\surd(aa+2ax-3xx)}{2a}\). This reappears on the next
  line. \JE}
\[y=\frac{ac-cx\pm \surd(aa+2ax-3xx)}{2a},\]
which value is elegantly expressed thus:
\[y=\frac{c(a-x)\pm c\surd(a-x)(a+3x)}{2a}.\]\newpage
From this equation
it is firstly clear that the abscissa \(x\) can never become
greater than \(a\), and in the negative direction the abscissa
cannot increase beyond
\(\frac{1}{3}a\). \marginpar{Tab. II.\\Fig. 4.} But by having
taken \(\geom{B\,D}=\frac{1}{3}a\) to the left,\footnote{i.e. to
  the left of the triangle; see the figure. \JE} it will make
the applicate
\(\geom{D\,E}=\frac{2}{3}c=\frac{2}{3}\geom{B\,C}\), and
moreover at this point \(\geom{E}\) the two values of \(y\)
coalesce, and thus the line \(\geom{D\,E}\) touches the curve at
the point \(\geom{E}\). But now if the line \(\geom{A\,E}\) be
drawn intersecting the side \(\geom{B\,C}\) in \(\geom{F}\), by
similar triangles it will be that
\(\geom{A\,D}\,:\,\geom{D\,E}=\geom{A\,B}\,:\,\geom{B\,F}\),
whence it yields\footnote{The similar triangles
  \(\geom{A\,F\,B}\) and \(\geom{A\,E\,D}\) are related by
  scaling by a factor of
  \(\geom{A\,D}\,:\,\geom{A\,B}=\frac{4}{3} a\,:\,a\), so
  \(\geom{B\,F}=\frac{3}{4}\,\geom{D\,E}=\frac{3}{4}\cdot\frac{2}{3}\,c=\frac{1}{2}c\). \JE}
\(\geom{B\,F}=\frac{1}{2}c=\frac{1}{2}\geom{B\,C}\), so that the line
\(\geom{A\,F\,E}\) bisects the side \(\geom{B\,C}\). In turn indeed, since
\(\geom{B\,D}\) is a third part of \(\geom{B\,A}\), it will also be that
\(\geom{E\,F}\) is a third part of \(\geom{A\,F}\), whence the point \(\geom{E}\) is
most easily marked out.

\begin{figure}[htb]
  \centering
  \begin{tikzpicture}
    \node (A) at (-20:2.75) {};
    \node (B) at (206:2.15) {};
    \node (C) at (112.5:2) {};
    \node (D) at (-1.5+-1.93,-0.942) {};
    \node (E) at (160:2.8) {};
    \node (F) at (160:1.425) {};
    \node (G) at (26:1.02) {};
    \node (H) at (292.5:1.02) {};
    \draw (A.center) -- (E.center);
    \draw (A.center) -- (D.center);
    \draw (B.center) -- (G.center);
    \draw (C.center) -- (H.center);
    \node at (0,0) {\(\geom{O}\)};
    \node at (A) [below] {\(\geom{A}\)};
    \node at (B) [below] {\(\geom{B}\)};
    \node at (C) [above] {\(\geom{C}\)};
    \node at (D) [below] {\(\geom{D}\)};
    \node at (E) [above] {\(\geom{E}\)};
    \node at (-1.339-0.07,0.487+0.25) {\(\geom{F}\)};
    \node at (G) [above right] {\(\geom{G}\)};
    \node at (H) [below] {\(\geom{H}\)};
    \draw (A.center) -- (C.center) -- (B.center);
    \draw (E.center) -- (D.center);
    \node[anchor=west] at (0.31,1.7) {\textit{Tab. II}};
    \node[anchor=west] at (0.31,1.3) {\textit{Fig. 4}};
  \end{tikzpicture}
\end{figure}

\bigskip

\S.\,\,\,11.\quad But if in turn we were to make \(x=a\),
because of the vanishing of the square root both values of \(y\)
also coalesce, with both being \(=0\), or with
\(\geom{A}\,\alpha\) drawn parallel to \(\geom{B\,C}\) and
infinitely short, \(\alpha\) will also be a point in the
curve,\footnote{i.e. \(\geom{B\,C}\) is parallel to the tangent
  line at \(\geom{A}\). \JE} and so the line
\(\geom{A}\,\alpha\) touches the curve at \(\geom{A}\), which
since it is parallel to the tangent \(\geom{D\,E}\), it follows
that the line \(\geom{A\,E}\) is a diameter of the ellipse,
whose centre will therefore fall at its midpoint
\(\geom{O}\). Therefore since
\(\geom{F\,E}=\frac{1}{3}\geom{A\,F}\), the centre of the
ellipse will fall at the point \(\geom{O}\), by having taken
\(\geom{A\,O}=\frac{2}{3}\geom{A\,F}\), or
\(\geom{F\,O}=\frac{1}{3}\geom{A\,F}\). Therefore that diameter
will bisect all ordinates\footnote{Euler switches here from
  ``applicate'' to ``ordinate'', but later switches back. \JE}
parallel to the side \(\geom{B\,C}\). Moreover indeed since the
three points may be permuted with one another, in a similar way
the line \(\geom{B\,G}\) bisecting the side \(\geom{A\,C}\), and
indeed also the side \(\geom{C\,H}\), bisecting the side
\(\geom{B\,A}\), will all intersect one another in the same
point \(\geom{O}\). But it is known that the centre of gravity
of the triangle is determined in this way, whence it shows this
remarkable property: \emph{That the centre of the smallest
  ellipse of all passing through three points \(\geom{A}\),
  \(\geom{B}\), \(\geom{C}\), falls at the very centre of
  gravity of the triangle \(\geom{A}\), \(\geom{B}\),
  \(\geom{C}\)}; whence since moreover not only the three points
\(\geom{A}\), \(\geom{B}\), \(\geom{C}\) be given but also the
tangents at these points, namely which are parallel to the
opposite sides, that desired ellipse can be most easily
constructed from known properties of conic sections.

\begin{figure}[htb]
  \centering
  \begin{tikzpicture}
    \node (A) at (-30:2) {};
    \node (B) at (210:2) {};
    \node (C) at (90:2) {};
    \draw (B.center) -- ++ (60:4.1) node (V) {};
    \draw (B.center) -- (A.center) -- (C.center);
    \draw (B.center) -- ++ (30:3.6) node (T) {};
    \node at (T) [right] {\(\geom{T}\)};
    \node (X) at ([xshift=1.05cm]B) {};
    \node at (X) [below] {\(\geom{X}\)};
    \draw (X.center) -- ++ (60:3.075) node [right] {\(\geom{Y}\)};
    \node at (A) [below] {\(\geom{A}\)};
    \node at (B) [below] {\(\geom{B}\)};
    \node at ([shift=({-0.2,0.1})]C) {\(\geom{C}\)};
    \node at (0.1,-0.17) {\(\geom{O}\)};
    \draw (V.center) -- (T.center);
    \node at ([shift=({0.2,0.2})]V) {\(\geom{V}\)};
    \node at (0.866+0.3,0.45) {\(\geom{G}\)};
    \draw (X.center) -- ++ (120:1.05) node (S) {};
    \node at ([shift=({-0.15,0.1})]S) {\(\geom{S}\)};
    \node at (-1.73,-1+2.8) {\textit{Tab. II}};
    \node at (-1.73,-1+2.4) {\textit{Fig. 5}};
  \end{tikzpicture}
\end{figure}

\bigskip

\S.\,\,\,12.\quad Let us develop several cases. And first indeed
\marginpar{Tab. II.\\Fig. 5.} let the triangle
\(\geom{A\,B\,C}\) be equilateral, and so \(c=a\) and the angle
\(\omega=60\degree\), of which the sine \(=\frac{\surd 3}{2}\),
and the equation of the minimal ellipse will be
\(xx+xy+yy-ax-ay=0\), indeed the area itself of this minimal
ellipse will be \(=\frac{1}{3}\pi aa\). But it is easily
discerned that in this case the ellipse will be a circle
circumscribed around the triangle. Indeed from the equation this
can be shown thus: Let \(\geom{Y}\) be a point on the ellipse,
whence let a parallel \(\geom{Y\,X}\) to the side
\(\geom{B\,C}\) be made, so that \(\geom{B\,X}=x\) and
\(\geom{X\,Y}=y\); then indeed let a line \(\geom{B\,G}\) be
extended, bisecting the side \(\geom{A\,C}\) in \(\geom{G}\), to
which from \(\geom{Y}\) let a normal \(\geom{Y\,T}\) be drawn,
and call\footnote{{\em voceturque} -- If we translate literally
  as ``and let it be called'' then ``it'' would sound like it
  referred to \(\geom{Y\,T}\), so I opted for ``call'' even
  though this is neither passive nor subjunctive. \JE}
\(\geom{B\,T}=t\) and \(\geom{T\,Y}=u\). Let \(\geom{T\,Y}\) be
extended until it meet the extended side \(\geom{B\,C}\) in
\(\geom{V}\), then indeed moreover let a line \(\geom{X\,S}\) be
made parallel to \(\geom{A\,C}\), and it will be that
\(\geom{B\,S}=\geom{X\,S}=x=\geom{Y\,V}\). Therefore it will
also be that \(\geom{S\,V}=\geom{X\,Y}=y\), and so
\(\geom{B\,V}=x+y\), and hence because the angle
\(\geom{C\,B\,G}=30\degree\) it will be that
\(\geom{B\,T}=t=\frac{\surd 3}{2}(x+y)\) and
\(\geom{T\,V}=\frac{1}{2}(x+y)\). From this let
\(\geom{Y\,V}=x\) be removed and it leaves
\(\geom{T\,Y}=u=\frac{1}{2}(y-x)\).

\bigskip

\S.\,\,\,13.\quad From the equations here discovered it will be
first that \(x+y=\frac{2t}{\surd 3}\) and \(y-x=2u\), whence it
is inferred \(x=\frac{t}{\surd 3}-u\) and \(y=\frac{t}{\surd
  3}+u\), by substituting which values an equation will arise
between the rectangular coordinates \(t\) and \(u\), which will
be
\[tt+uu-\frac{2at}{\surd 3}=0,\quad\text{or}\quad
  uu=+\frac{2at}{\surd 3}-tt,\] which is manifestly for a circle
whose radius \(=\frac{a}{\surd 3}\). Therefore since the line
\(\geom{B\,G}=\frac{a\surd 3}{2}\), \(\geom{B\,G}\) will be to
that radius as \(3\,:\,2\), and so the centre of the circle
falls at \(\geom{O}\), so that
\(\geom{B\,O}=\frac{2}{3}\geom{B\,G}\).

\begin{figure}[htb]
  \centering
  \begin{tikzpicture}
    \node (A) at (-27.5:1.9) {};
    \node at (A) [below] {\(\geom{A}\)};
    \node (B) at (201.5:2.4) {};
    \node at (B) [below] {\(\geom{B}\)};
    \node (C) at (70.5:1.9) {};
    \node at ([shift=({-0.1,0.15})]C) {\(\geom{C}\)};
    \draw (B.center) -- ++ (43:4.5) node (V) {};
    \node at ([shift=({0,0.2})]V) {\(\geom{V}\)};
    \draw (B.center) -- (A.center) -- (C.center);
    \node (X) at ([xshift=0.8cm]B) {};
    \node at (X) [below] {\(\geom{X}\)};
    \draw (X.center) -- ++ (43:3.7) node [right] {\(\geom{Y}\)};
    \draw (B.center) -- ++ (21.5:4.2) node (T) {} --++ (21.5:0.35) node [right] {\(\geom{I}\)};
    \node at ([shift=({0.05,-0.2})]T) {\(\geom{T}\)};
    \draw (V.center) -- (T.center);
    \draw (0,0) node {} -- ++ (151:2) node [left] {\(\geom{K}\)};
    \node at (0,0) [below] {\(\geom{O}\)};
    \draw (X.center) -- ++ (110:0.6) node (S) {};
    \node at ([shift=({-0.1,0.15})]S) {\(\geom{S}\)};
  \end{tikzpicture}
\end{figure}

\bigskip

\S.\,\,\,14.\quad \marginpar{Tab. II.\\ Fig. 6.} Now let the
triangle \(\geom{A\,B\,C}\) be isosceles and
\(\geom{B\,C}=\geom{B\,A}=a=c\), indeed the angle at
\(\geom{B}\), which was \(\omega\), let us now set as
\(\omega=2\theta\), so that with the line \(\geom{B\,G\,T}\)
having been drawn bisecting the side \(\geom{A\,C}\) and normal
to it, the angle will be \(\geom{C\,B\,G}=\theta\). Now let it
be set \(\geom{B\,T}=t\) and \(\geom{T\,Y}=u\) as before, with
it being \(\geom{B\,X}=x\) and \(\geom{X\,Y}=y\); then let the
parallelogram \(\geom{X\,S\,V\,Y}\) be completed, and it will be
that \(\geom{B\,S}=x\), and
\(\geom{X\,S}=2x\sin\theta=\geom{Y\,V}\). Since now
\(\geom{S\,V}=\geom{X\,Y}=y\), it will be that
\(\geom{B\,V}=x+y\), and hence
\begin{align*}
  &\geom{B\,T}=t=(x+y)\cos\theta\quad\text{and}\\
  &\geom{T\,V}=(x+y)\sin\theta.
\end{align*}
Let \(\geom{Y\,V}=2x\sin\theta\) be removed from this, and there
will remain \[\geom{T\,Y}=u=(y-x)\sin\theta.\] Wherefore since
we have
\[x+y=\frac{t}{\cos\theta}\quad\text{and}\quad
  y-x=\frac{u}{\sin\theta},\] the obliquangular coordinates will
be
\[x=\frac{t}{2\cos\theta}-\frac{u}{2\sin\theta}\quad\text{and}\quad
  y=\frac{t}{2\cos\theta}+\frac{u}{2\sin\theta},\] which
substituted values furnish this equation between the orthogonal
coordinates:
\begin{align*}
  &\frac{3tt}{4\cos\theta^2}+\frac{uu}{4\sin\theta^2}-\frac{at}{\cos\theta}=0,\quad\text{or}\\
  &uu=\frac{4at\sin\theta^2}{\cos\theta}-\frac{3tt\sin\theta^2}{\cos\theta^2},
\end{align*}
whence it is clear the applicate \(u\) vanishes in either the
case \(t=0\) or the case \(t=\frac{4}{3}a\cos\theta\), which
quantities will therefore give the principal axis of the
ellipse, which is
\(\geom{B\,I}=\frac{4}{3}a\cos\theta\). Wherefore since
\(\geom{B\,G}=a\cos\theta\), it will be that
\(\geom{B\,I}=\frac{4}{3}\geom{B\,G}\), and thus the centre will
fall at \(\geom{O}\), with it being
\(\geom{B\,O}=\geom{O\,I}=\frac{2}{3}\geom{B\,G}\). But if now
we were to take \(t=\geom{B\,O}=\frac{2}{3}a\cos\theta\), the
value of \(u\) will give the conjugate semiaxis,\footnote{The
  position of \(\geom{K}\) in Tab. II, Fig. 6 looks wrong. I
  think that \(\geom{O\,K}\) should be orthogonal to
  \(\geom{B\,I}\). \JE} which is \(\geom{O\,K}\), and it will be
that \(\geom{O\,K}=\frac{2a\sin\theta}{\surd 3}\), and indeed
the other semiaxis was
\(\geom{B\,O}=\geom{O\,I}=\frac{2}{3}a\cos\theta\), whence it
yields the area for this ellipse:
\[=\pi\,.\,\geom{B\,O}\,.\,\geom{O\,K} = \frac{4\pi
    aa\sin\theta\cos\theta}{3\surd 3}=\frac{2\pi
    aa\sin\omega}{3\surd 3},\] which agrees perfectly with the
general form.

\vspace{1.5cm}

\begin{center}
  \begin{tikzpicture}
    \draw[very thick] (0,0) -- (8,0);
    \draw (0,-0.1) -- (8,-0.1);
  \end{tikzpicture}
\end{center}

\newpage
\appendix
\begin{center}
  \medskip

  {\Huge\sc
    Appendices}

  \medskip
\end{center}

\section{Euler's cubic from E691}
\label{app:cubic}

In this appendix, we will discuss the minimal area problem for
ellipses circumscribing quadrilaterals (considered in E691) from
a modern point of view, and solve the problem whose
investigating Euler left to others, namely the problem of
whether there are further critical points of the area function
in the pencil of conics (i.e. roots of his cubic equation) and
what they mean geometrically.

\subsection{Euler's area formula}

Despite its complicated appearance, Euler's formula
\begin{equation}\label{eq:app_area}
  \text{area} =
  \pi\OP{sin}\omega\left(\frac{CD^2+AE^2-2BDE-F(AC-B^2)}{(AC-B^2)^{3/2}}\right)
\end{equation}
has a simple interpretation in terms of determinants. Namely, if
we write the equation
\begin{equation}\label{eq:app_conic}
  Ax^2+2Bxy+Cy^2+2Dx+2Ey+F=0
\end{equation}
of the ellipse in matrix form:
matrix equation
\(\mathbf{V}^T\mathbf{M}\mathbf{V}=0\), where
\[\mathbf{M}=
  \begin{pmatrix} A & B & D \\ B & C & E \\ D & E & F\end{pmatrix}\quad\text{and}\quad\mathbf{V}=
  \begin{pmatrix} x \\ y \\ 1\end{pmatrix}\] and write
\(\mathbf{N}=\begin{pmatrix} A & B \\ B & C\end{pmatrix}\) then
Equation \eqref{eq:app_area} becomes
\[\text{area}=-\pi\OP{sin}\omega\frac{\det(\mathbf{M})}{\det(\mathbf{N})^{3/2}}.\]
In fact, the factor \(\OP{sin}\omega\) is something of a red
herring: we will soon see that it arises as the determinant of
the matrix
\(\begin{pmatrix}\OP{cos}(\omega/2) &
  \OP{sin}(\omega/2)\\\OP{sin}(\omega/2) &
  \OP{cos}(\omega/2)\end{pmatrix}\) which changes coordinates
from orthogonal axes to obliquangular axes meeting at an angle
\(\omega\). The overall sign is also not very important: it can
be changed by multiplying \(\mathbf{M}\) by \(-1\). Let us
therefore focus on the quantity
\(\frac{\det(\mathbf{M})}{\det(\mathbf{N})^{3/2}}\).

Let \(\mathbf{v}=\begin{pmatrix} x \\
  y \end{pmatrix}\) in standard Cartesian coordinates. Consider
a new affine coordinate system \((x',y')\) related to this one
by
\[\mathbf{v}=\mathbf{sv}'+\mathbf{t}\] where
\(\mathbf{s}=\begin{pmatrix} s_{11} & s_{12} \\ s_{21} &
  s_{22}\end{pmatrix}\) is a nonsingular 2-by-2 matrix and
\(\mathbf{t}=\begin{pmatrix} t_1 \\ t_2\end{pmatrix}\) is a
vector then this sends
\(\mathbf{V}'\) to \(\mathbf{V}=\mathbf{SV}'\) where
\[\mathbf{S}=\begin{pmatrix}s_{11} & s_{12} & t_1 \\
    s_{21} & s_{22} & t_2 \\ 0 & 0 &
    1\end{pmatrix}.\] The equation of the ellipse in the new
coordinates is \((\mathbf{V}')^T\mathbf{M}'\mathbf{V}'\) where
\(\mathbf{M}'=\mathbf{S}^T\mathbf{MS}\) since
\[\mathbf{V}^T\mathbf{MV}=(\mathbf{V}')^T\mathbf{S}^T\mathbf{MSV}'.\]
The top-left 2-by-2 block \(\mathbf{N}\) becomes
\(\mathbf{N}'=\mathbf{s}^T\mathbf{Ns}\). Observe that
\begin{align*}
  \det(\mathbf{M}')&=\det(\mathbf{S}^T\mathbf{MS})=\det(\mathbf{S}^T)\det(\mathbf{M})\det(\mathbf{S})=\det(\mathbf{S})^2\det(\mathbf{M})\\
  \det(\mathbf{N}')&=\det(\mathbf{s}^T\mathbf{Ns})=\det(\mathbf{s}^T)\det(\mathbf{N})\det(\mathbf{s})=\det(\mathbf{s})^2\det(\mathbf{N})
\end{align*}
and since \(\det(\mathbf{S})=\det(\mathbf{s})\), both
\(\det(\mathbf{M})\) and \(\det(\mathbf{N})\) transform under
the affine transformation by a factor of
\(\det(\mathbf{s})^2\). If instead we rescale the equation of
the ellipse by a factor of \(\lambda>0\) then \(\det(\mathbf{M})\)
and \(\det(\mathbf{N})\) scale as
\(\det(\lambda\mathbf{M})=\lambda^3\det(\mathbf{M})\) and
\(\det(\lambda\mathbf{N})=\lambda^2\det(\mathbf{N})\)
respectively. The quantity
\[\det(\mathbf{M})^m\det(\mathbf{N})^n\] therefore scales like
\(\det(\mathbf{s})^{2(m+n)}\lambda^{3m+2n}\) under a combination
of affine transformations and positive rescalings of the
equation.

The area of the ellipse is unaffected by rescaling the equation
of the ellipse, and scales like \(\det(\mathbf{s})^{-1}\) under
affine transformations\footnote{These are passive coordinate
  transformations, so this is \(\det(\mathbf{s})^{-1}\) rather
  than \(\det(\mathbf{s})\).} so it behaves exactly like the
quantity \(\det(\mathbf{M})/\det(\mathbf{N})^{3/2}\), which we
obtain by solving \[m+n=-1/2,\qquad 3m+2n=0\] to get \(m=1\),
\(n=-3/2\). This means that the ratio
\[\text{area}\,:\,\det(\mathbf{M})/\det(\mathbf{N})^{3/2}\] is
invariant under affine transformations. Since affine
transformations act transitively on ellipses, we can use an
affine change of coordinates to make our ellipse into the unit
circle, for which the value of this ratio is \(\pi\), and,
taking
\[\mathbf{M}=\begin{pmatrix}-1 & 0 & 0 \\ 0 & -1 & 0 \\ 0 & 0 &
    1\end{pmatrix},\] we see that
\[\text{area}=\pi\det(\mathbf{M})/\det(\mathbf{N})^{3/2}.\]
Euler's factor of \(-\OP{sin}\omega\) arises because (a) he is
taking the matrix for the circle to be \(-\mathbf{M}\) and (b)
he is working in a different coordinate system again,
\(\mathbf{v}''\), where the axes meet with an angle \(\omega\);
this is related to the rectangular coordinates by the affine
transformation
\[\mathbf{v}''=\begin{pmatrix}\OP{cos}(\omega/2) &
  \OP{sin}(\omega/2)\\\OP{sin}(\omega/2) &
  \OP{cos}(\omega/2)\end{pmatrix}\mathbf{v},\] and the
determinant of this matrix
is \[\OP{cos}^2(\omega/2)-\OP{sin}^2(\omega/2)=\OP{sin}\omega.\] 

\subsection{Euler's cubic}

Euler's calculation of the differential of the area function
yielded
\begin{equation}\label{eq:app_cubic}\frac{d(\text{area})}{dB}=\frac{FB^3-4DEB^2+(3CD^2+3AE^2-ACF)B-2ACDE}{(AC-B^2)^{5/2}}\end{equation}
which we can now write as
\begin{equation}\label{eq:app_matrix_cubic}\frac{d(\text{area})}{dB}=\frac{1}{\det(\mathbf{N})^{5/2}}\left(\det(\mathbf{N})\frac{d\det(\mathbf{M})}{dB}-\frac{3}{2}\det(\mathbf{M})\frac{d\det(\mathbf{N})}{dB}\right).\end{equation}
Euler's comment in E691 was an encouragement to consider the
nature of the critical points when the cubic in the numerator
happens to have three real roots. In fact, this {\em always}
happens. To see why, recall that in Euler's coordinate system
the four points are at \((a,0)\), \((b,0)\), \((0,c)\) and
\((0,d)\) and the coefficients in the equation of the ellipse
are
\begin{gather*}
  A=cd;\quad C=ab;\quad 2D=-cd(a+b);\\
  2E=-ab(c+d)\quad\text{and}\quad F=abcd.
\end{gather*}
By relabelling points, we can assume that the intersection point
\(\geom{O}\) between the axes \(\geom{A\,B}\) and
\(\geom{C\,D}\) is inside the ellipse, so \(a\) and \(b\) have
opposite signs and \(c\) and \(d\) have opposite signs. Scaling
the axes independently, we can assume \(ab=-1\) and \(cd=-1\),
so \(A=C=-1\), \(F=1\) and \(D=(a+b)/2\), \(E=(c+d)/2\). If we
set \(\sigma=(a+b)^2\geq 0\) and \(\tau=(c+d)^2\geq 0\), then the
cubic becomes
\[B^3-\sqrt{\sigma\tau}
  B^2+\left(\frac{3}{4}(\sigma+\tau)-1\right)B -
  \frac{1}{2}\sqrt{\sigma\tau}.\] Using Sage, we can compute the
discriminant of this cubic:
\begin{gather*}\frac{1}{16}\left(9\sigma^3\tau-14\sigma^2\tau^2+9\sigma\tau^3\right)\\
  +\frac{3}{16}\left(9\sigma^3-\sigma^2\tau-\sigma\tau^2+9\tau^3\right)\\+\frac{1}{4}\left(27\sigma^2-5\sigma\tau+27\tau^2\right)\\+9\sigma+9\tau+4\end{gather*}
We can group these terms as follows:
\begin{align*}
  9\sigma^3\tau-14\sigma^2\tau^2+9\sigma\tau^3=\sigma\tau(7(\sigma-\tau)^2+2\sigma+2\tau)&\geq
  0\\
  9\sigma^3-\sigma^2\tau-\sigma\tau^2+9\tau^3=8\sigma^3+8\tau^3+(\sigma+\tau)(\sigma-\tau)^2&\geq
  0\\
  27\sigma^2-5\sigma\tau+27\tau^2=\frac{5}{2}(\sigma-\tau)^2+\frac{49}{2}(\sigma^2+\tau^2)&\geq 0
\end{align*}
so all of the terms are non-negative, and the constant term is
\(4\), so the discriminant is positive. This means that Euler's
cubic always has three distinct real roots. In what follows, we
will see that precisely one of them lies in the domain
\(-\sqrt{AC}<B<\sqrt{BC}\) and explain the nature of the other
two, depending on the shape of the quadrilateral \(\geom{A\,B\,C\,D}\).

\subsection{Hyperbolas}

When \(B^2=AC\), the conic in the pencil is either a parabola or
a pair of parallel lines; in other words it is a conic with
precisely one asymptote (with multiplicity \(2\)). At such a
value of \(B\), Euler's formula has a pole. However, when
\(B^2>AC\), the only issue with the area formula is that it
takes on imaginary values. If we replace the denominator
\((AC-B^2)^{3/2}\) by \((B^2-AC)^{3/2}\) we obtain a real-valued
function of \(B\) for the domain \(B^2>AC\). Let us write
\(\mathcal{A}\) for the function
\[\mathcal{A}(B)=\begin{cases}
\frac{\det(\mathbf{M})}{(AC-B^2)^{3/2}} &\text{ if
}B^2<AC\smallskip\\
\frac{\det(\mathbf{M})}{(B^2-AC)^{3/2}} &\text{ if }B^2>AC
\end{cases}\] Crucially, Euler's calculation of critical points
is unaffected, since we are just rescaling his function by a
piecewise constant scale factor (equal to \(1\) on \(B^2<AC)\)
and \(-i\) on \(B^2>AC\)).

What does this new quantity \(\mathcal{A}\) mean for hyperbolas?
The argument we gave before shows that it still transforms like
an area under affine transformations. Moreover, affine
transformations act transitively on hyperbolas. Therefore if we
can interpret \(\mathcal{A}\) geometrically as an area for a
single hyperbola using constructions that only involve geometric
features which are preserved by affine transformations (like
parallelism, tangency, or asymptotes) then this geometric
interpretation will hold for any hyperbola.

\begin{figure}[htb]
  \centering
  \begin{tikzpicture}
    \draw (-2,0) -- (4.5,0);
    \draw (0,-3) -- (0,3);
    \draw[domain=0.35:4.5, smooth, variable=\x,very thick] plot ({\x},{1/\x});
    \draw[domain=-2:-0.35, smooth, variable=\x,very thick] plot ({\x},{1/\x});
    \filldraw[fill=lightgray,draw=black] (0,0) -- (0,1) -- (4,0) -- cycle;
    \draw [decorate,decoration={brace,amplitude=5pt,mirror,raise=0.5ex}] (0,1) -- (0,0) node[midway,xshift=-2em]{\(2/x\)};
    \draw [decorate,decoration={brace,amplitude=5pt,mirror,raise=0.5ex}] (0,0) -- (4,0) node[midway,yshift=-1.5em]{\(2x\)};
    \node at (2,1/2) {\(\bullet\)};
    \node at (2,1/2) [above right] {\((x,1/x)\)};
  \end{tikzpicture}
  \caption{The area invariant for a hyperbola.}
  \label{fig:app_hyperbola_triangle}
\end{figure}
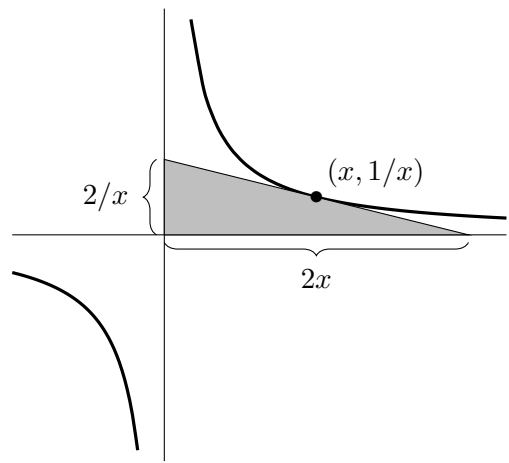

Consider the hyperbola \(xy=1\) and pick a point
\(p=(x,1/x)\). Look at the triangle formed by the tangent line
at \(p\) and the asymptotes (Figure
\ref{fig:app_hyperbola_triangle}); this has area \(1/2\),
independently of the point \(p\). If we make an affine
transformation, since asymptotes are sent to asymptotes and
tangent lines to tangent lines, and since all areas scale by the
same factor, every hyperbola has this property: that the area of
the triangle formed by the asymptotes and any tangent line is
independent of the choice of tangent line. Let us write
\(\mathcal{T}\) for this area invariant of a hyperbola. Since
all these geometric features are preserved by affine
transformations, and since both \(\mathcal{A}\) and
\(\mathcal{T}\) scale as areas under affine transformations, we
deduce that the ratio \(\mathcal{A}\,:\,2\mathcal{T}\) is
unchanged by affine transformations. For the hyperbola \(xy=1\),
\[\mathbf{M}=\begin{pmatrix} 0 & -\frac{1}{2} & 0 \\ -\frac{1}{2}
    & 0 & 0 \\ 0 & 0 & 1\end{pmatrix}\] so
\[\mathcal{A}=\frac{\det(M)}{(-\det(N))^{3/2}}=\frac{1/4}{1/8}=2,\]
so we conclude that
\[\mathcal{T}=\frac{1}{4}\mathcal{A}\] for any hyperbola. So if
there are roots of Euler's cubic \eqref{eq:app_cubic} in the
domain \(B^2>AC\), they will correspond to critical points of
this area invariant \(\mathcal{T}\).

Note that \(\mathcal{A}\) may vanish or even become negative: we
should think of it as a signed area. For example, in the family
\(xy=t\), we have \(\mathcal{A}=\frac{t/4}{1/8}=2t\). When
\(t=0\), the hyperbola degenerates to a union of two lines, and
\(\mathcal{A}\) and \(\mathcal{T}\) vanish. When \(t\) becomes
negative, it is natural to consider both \(\mathcal{A}\) and
\(\mathcal{T}\) to be negative. If instead we took the absolute
value, we would get \(\mathcal{A}=2|t|\), which is not
differentiable at \(t=0\). The choice of sign is fixed by taking
\((-\det(N))^{3/2}\) to be positive.

\subsection{Four cases: the moduli space of quadrilaterals}

Making an affine transformation simply rescales \(\mathcal{A}\)
by a constant factor \(\det(\mathbf{s})^{-1}\), so we have the
freedom to move our quadrilateral around by an affine
transformation. How many conelliptic quadrilaterals
\(\geom{A\,B\,C\,D}\) are there up to affine transformations?

\begin{figure}[htb]
  \centering
  \begin{tikzfadingfrompicture}[name=fade left]
        \shade[left color=transparent!0, right color=transparent!100] (0,0) rectangle (2,2);
      \end{tikzfadingfrompicture}
  \begin{tikzpicture}
    \filldraw[path fading=fade left,fill=lightgray,draw=none] (6,6) -- (1,1) -- (6,1) -- cycle;
    \draw[dashed] (5,1) -- (6,1);
    \draw[dashed] (5,5) -- (6,6);    
    \draw[->] (0,0) -- (5,0) node [right] {\(a\)};
    \draw[->] (0,0) -- (0,5) node [left] {\(c\)};
    \draw (1,1) -- (5,5) node [midway,sloped] {trapezia};
    \draw (1,1) -- (5,1);
    \node at (3,1) {kites};
    \node[fill=white] at (1,1) [left] {parallelograms};
    \node at (1,1) {\(\bullet\)};
    \node at (4,2) {irregular};
    \begin{scope}[shift={(-0.2,2)},scale=0.5]
      \filldraw[fill=white,draw=black] (-1,0) node {\(\bullet\)} -- (0,-1) node {\(\bullet\)} -- (1,0) node {\(\bullet\)} -- (0,1) node {\(\bullet\)} -- cycle;
      \draw (-1.5,0) -- (1.5,0);
      \draw (0,-1.5) -- (0,1.5);
    \end{scope}
    \begin{scope}[shift={(3,-0.2)},scale=0.5]
      \filldraw[fill=white,draw=black] (-1,0) node {\(\bullet\)} -- (0,-1) node {\(\bullet\)} -- (2,0) node {\(\bullet\)} -- (0,1) node {\(\bullet\)} -- cycle;
      \draw (-1.5,0) -- (2.5,0);
      \draw (0,-1.5) -- (0,1.5);
    \end{scope}
    \begin{scope}[shift={(2,4)},scale=0.5]
      \filldraw[fill=white,draw=black] (-1,0) node {\(\bullet\)} -- (0,-1) node {\(\bullet\)} -- (2,0) node {\(\bullet\)} -- (0,2) node {\(\bullet\)} -- cycle;
      \draw (-1.5,0) -- (2.5,0);
      \draw (0,-1.5) -- (0,2.5);
    \end{scope}
    \begin{scope}[shift={(5,3)},scale=0.5]
      \filldraw[fill=white,draw=black] (-1,0) node {\(\bullet\)} -- (0,-1) node {\(\bullet\)} -- (4,0) node {\(\bullet\)} -- (0,2) node {\(\bullet\)} -- cycle;
      \draw (-1.5,0) -- (2.5,0);
      \draw (0,-1.5) -- (0,1.5);
    \end{scope}
  \end{tikzpicture}
  \caption{The moduli space of quadrilaterals up to affine transformations.}
  \label{fig:moduli}
\end{figure}
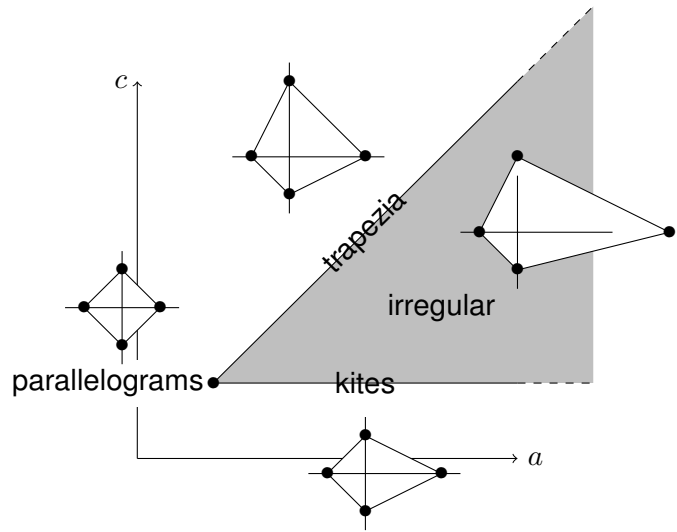

Like Euler, we can use an affine transformation to make
\(\geom{A\,B}\) and \(\geom{C\,D}\) our axes. By possibly
relabelling points, we can assume that the point \(\geom{O}\)
where these lines intersect is inside the quadrilateral, so that
\[\geom{A}=(a,0),\quad \geom{B}=(b,0),\quad
  \geom{C}=(0,c),\quad \geom{D}=(0,d),\] where \(a\) and \(b\)
have opposite signs and \(c\) and \(d\) have opposite signs. By
reflecting in the axes, we can assume that \[0<-b\leq
  a\quad\text{and}\quad 0<-d\leq c.\] By rescaling the axes
independently, we can assume that \(b=d=-1\). By swapping
\(\geom{A\,B}\) and \(\geom{C\,D}\) if necessary, we can further
assume that \(a\geq c\). We have now used up all of the freedom
in choosing our coordinates, and we see that the quadrilateral
is determined up to affine transformations by the two numbers
\(a\) and \(c\), which live in the wedge
\[\{(a,c)\in\RR^2\,:\,1<c\leq a\}.\] This is the {\em moduli
  space of conelliptic quadrilaterals} up to affine
transformations (Figure \ref{fig:moduli}).

\begin{itemize}
\item The vertex \((1,1)\) of the wedge corresponds to a square
  \(a=c=1\), \(b=d=-1\). The images of the square under affine
  transformations are parallelograms, and all parallelograms are
  affine-equivalent to a square.
\item The points along the diagonal have \(a=c\) which
  correspond to {\em trapezia} in the original sense of Proclus:
  one pair of opposite sides (\(\geom{A\,C}\) and
  \(\geom{B\,D}\)) is parallel.
\item The points along the bottom edge have \(c=-d\), which correspond
  to {\em kites}.
\item The remaining quadrilaterals in the interior of the wedge
  are irregular in the sense that they have no affine
  symmetries.
\end{itemize}
The kites and trapezia have an affine symmetry group of order
\(2\) comprising a reflection. The square has an affine symmetry
group of order \(8\) (the dihedral group). The moduli space is
really an orbifold; these edges are part of the orbifold locus
with stabiliser \(\ZZ/2\) and the vertex is an orbifold point
with stabiliser \(D_8\) (this orbifold is just the quotient of
\(\RR^2\) by the symmetries of the square).

In fact, a better system of coordinates to use on the moduli
space would
be \[\sigma=\frac{(a+b)^2}{ab}\quad\text{and}\quad\tau=\frac{(c+d)^2}{cd}\]
since these are actually affine invariants.

There are therefore four broad ``types'' of quadrilateral up to
affine transformations: parallelograms, kites, trapezia and
irregular quadrilaterals. The parallelograms form a single orbit
of the group of affine transformations, whereas the kites,
trapezia and irregular quadrilaterals are all comprised of
infinitely many orbits distinguished by the ``moduli''
\(\sigma\) and \(\tau\) (so, for example, it's possible to have
two kites which are not related by any affine
transformation). This classification of quadrilaterals will help
us to organise the subsequent case analysis.

\subsection{Degenerate and limiting conics}

{\bf Degenerate conics.} Recall that a conic is called {\em
  degenerate} if \(\det(\mathbf{M})=0\). A degenerate conic is a
pair of lines. Since our points are conelliptic (and hence no
three are collinear), the only degenerate conics in the pencil
are the pairs of lines
\[\geom{A\,B}\cup\geom{C\,D},\quad
  \geom{A\,C}\cup\geom{B\,D},\quad\text{and }
  \geom{A\,D}\cup\geom{B\,C}.\] The first of these consists of
the axes, and can be thought of as the conic \(B=\infty\) in the
pencil. The other two appear at the roots \(B\) of the
quadratic \[\det(\mathbf{M})=FB^2-2DEB+CD^2+AE^2-FAC,\] which
are (after some manipulation):
\[\frac{DE}{F}\pm\sqrt{\left(\frac{D^2}{F} -
      A\right)\left(\frac{E^2}{F} - C\right)}.\]

\noindent {\bf The limiting conics.} The poles of the function
\(\mathcal{A}\) occur at \(B=\pm \sqrt{AC}\), that is when
\(\det(\mathbf{N})=0\). The corresponding conics, which we have
been calling the ``limiting conics'' of the pencil, are
precisely the conics in the pencil with just one asymptote. To see
this, projectivise the equation of the conic by introducing a
third variable \(z\):
\[Ax^2+2Bxy+Cy^2+2Dxz+2Eyz+Fz^2=0.\] The slopes of the
asymptotes are the roots of the quadratic form
\(Ax^2+2Bxy+Cy^2\) obtained by setting \(z=0\), that is
\[\frac{-B\pm\surd(B^2-AC)}{C}.\] We see that the asymptotes
are parallel (both having slope \(-B/C\)) if and only if
\(B^2=AC\). A conic with a single (real) asymptote is either a
parabola or a pair of parallel lines (ellipses have no real
asymptotes and hyperbolas or non-parallel lines have two).
Putting these observations together with those about degenerate
conics, we deduce:

\bigskip

\noindent {\bf Lemma 1.} {\em A limiting conic \(B=+\sqrt{AC}\)
  (respectively \(B=-\sqrt{AC}\)) is either a parabola (if
  \(\det(\mathbf{M})\neq 0)\)) or a pair of parallel lines (if
  \(\det(\mathbf{M})=0\)). The latter can only happen if the
  lines \(\geom{A\,C}\) and \(\geom{B\,D}\) are parallel
  (respectively if \(\geom{A\,D}\) and \(\geom{B\,C}\) are
  parallel).}

\bigskip

\noindent This implies:

\bigskip

\noindent {\bf Proposition 2.} {\em Fix a conelliptic quadrilateral. The limiting
conics have the types shown in the following table:}

\begin{center}
  \begin{tabular}{c|c|c}
    Type & \(B=-\sqrt{AC}\) & \(B=+\sqrt{AC}\) \\\hline
         &&\\
    parallelogram & parallel lines & parallel lines\\
         & \(\geom{A\,D}\cup\geom{B\,C}\) &
                                            \(\geom{A\,C}\cup\geom{B\,D}\)\\
         &&\\
    trapezium & parallel lines & parabola\\
         & \(\geom{A\,D}\cup\geom{B\,C}\)&\\
         &&\\
    kite or irregular & parabola & parabola\\
  \end{tabular}
\end{center}

\noindent The apparent asymmetry in the case of a trapezium is
simply a manifestation of our convention that \(c\leq a\).

\subsection{The graph of the area function}

We now consider the shape of the graph of \(\mathcal{A}\) for
each of the four types of quadrilateral. The results are
illustrated in Figure \ref{fig:pencils_areas}. Recall that with
our choices \(0<-b=1\leq a\) and \(0<-d=1\leq c\), so:
\begin{gather*}
  A=cd<0,\quad C=ab<0,\quad F=AC>0,\\
  D=-\frac{1}{2}cd(a+b)\geq 0,\quad E=-\frac{1}{2}ab(c+d)\geq 0.
\end{gather*}

\bigskip

\noindent {\bf Lemma 3.} {\em The area function \(\mathcal{A}\)
  for ellipses, attains its minimum on the domain
  \(-\sqrt{AC}<B<\sqrt{AC}\).}
\begin{proof}  
  The numerator of \(\mathcal{A}\) is
  \(F(AC-B^2)-CD^2-AE^2+2BDE\). We have
  \begin{align*}
    -CD^2-AE^2+2BDE&=\left(-CD^2-AE^2+2DE\sqrt{AC}\right)-2DE\left(\sqrt{AC}-B\right)\\
                  &=\left(D\sqrt{-C}+ E\sqrt{-A}\right)^2-2DE\left(\sqrt{AC}-B\right).
  \end{align*}
  Let us write \(L\) for the (non-negative) quantity
  \(\left(D\sqrt{-C}+E\sqrt{-A}\right)^2\). Since
  \(D,E\geq 0\), this is positive unless \(D=E=0\). Then the
  numerator can be written as
  \[L+\left[\sqrt{AC}-B\right]\cdot\left[F\left(\sqrt{AC}+B\right)-2DE\right]\]
  The area function is therefore
  \[\mathcal{A}=\frac{L}{(AC-B^2)^{3/2}}+\frac{F\left(\sqrt{AC}+B\right)-2DE}{(\sqrt{AC}+B)^{3/2}\sqrt{\sqrt{AC}-B)}}.\]
  If \(L>0\) then the term
  \(\frac{L}{(AC-B^2)^{3/2}}\sim\mathcal{O}((\sqrt{AC}-B)^{-3/2})\)
  dominates as \(B\) approaches \(\sqrt{AC}\), and tends to
  \(+\infty\). If \(L=0\) then \(D=E=0\) so
  \(\mathcal{A}=\frac{F}{(AC-B^2)^{1/2}}\), which also tends to
  \(+\infty\) as \(B\to\sqrt{AC}\). A similar argument shows
  that \(\mathcal{A}\to +\infty\) as \(B\to-\sqrt{AC}\). Since
  the restriction of \(\mathcal{A}\) to any closed subinterval
  in \((-\sqrt{AC},\sqrt{AC})\) attains its minimum, and since
  \(\mathcal{A}\) can be bounded from below outside a
  sufficiently large closed subinterval (because it tends to
  \(+\infty\)), this shows that \(\mathcal{A}\) attains its
  minimum on the interval \((-\sqrt{AC},\sqrt{AC})\).
\end{proof}

\noindent {\bf Lemma 4.} {\em As \(B\to\pm\infty\), the function
  \(\mathcal{A}(B)\) tends to zero from above.}
\begin{proof}
  The dominant terms in the numerator and denominator of
  \(\mathcal{A}\) are \(FB^2\) and \(|B|^3\) respectively (since
  we are choosing \((B^2-AC)^{3/2}\) to be positive). With our
  choice of coordinates, since \(a\) and \(b\) have the same
  sign and \(c\) and \(d\) have the same sign,
  \(F=abcd>0\). Therefore \(\mathcal{A}(B)\sim F/B\) which is
  positive and tends to zero.
\end{proof}

\noindent {\bf Lemma 5.} {\em If \(\det(\mathbf{M})\neq 0\) when
  \(B=\sqrt{AC}\) (respectively \(B=-\sqrt{AC}\)) then, as
  \(B\to\sqrt{AC}\) from above (respectively \(B\to-\sqrt{AC}\)
  from below) the function \(\mathcal{A}(B)\) tends to
  \(-\infty\).}
\begin{proof}
  As \(B\to\pm\sqrt{AC}\), we have
  \[\det(\mathbf{M})\to CD^2+AE^2\mp
    2DE\sqrt{AC}=-(D\surd(-C)\mp E\surd(-A))^2\leq 0.\] If we
  assume \(\det(\mathbf{M})\neq 0\) at \(B=\pm\sqrt{AC}\) then
  this shows \(\det(\mathbf{M})\) is strictly negative at
  \(B=\pm\sqrt{AC}\). The denominator is \((B^2-AC)^{3/2}\)
  which we are assuming tends to zero from above. Therefore
  \(\mathcal{A}=\det(\mathbf{M})/(B^2-AC)^{3/2}\) is negative
  and tends to infinity.
\end{proof}

\noindent {\bf Critical hyperbolas from parabolic limiting
  conics.} By Lemma 1, we can have \(\det(\mathbf{M})\neq 0\)
when \(B=\pm\sqrt{AC}\) if and only if the corresponding
limiting conic is a parabola. Suppose then that the limiting
conic \(B=\sqrt{AC}\) is a parabola (which happens for anything
other than a parallelogram). Then the function \(\mathcal{A}\)
restricted to the domain \(\sqrt{AC}<B<\infty)\) must go from
having large negative values at one end to having positive but
decaying values at the other. This means that the derivative
\(\frac{d\mathcal{A}}{dB}\) must start off positive and end up
negative, so there must be a critical point of \(\mathcal{A}\)
in this domain. See Figure \ref{fig:pencils_areas}(b--d) for an
illustration. For kites and irregular quadrilaterals, the same
argument guarantees the presence of another critical point in
the domain \(-\infty<B<-\sqrt{AC}\).

\bigskip

\noindent {\bf Kites.} In fact, for kites we can say a little
more about these critical hyperbolas. We have \(c=-d\), so
\(E=0\) and Euler's cubic becomes
\[-abc^2B^3+ab\left(\frac{3}{4}c^4(a+b)^2-abc^4\right).\] This
has a root at \(B=0\) (the minimal ellipse) and two critical
hyperbolas at \[B=\pm \frac{c}{2}\sqrt{3a^2+3b^2+2ab}.\] Since
\(a^2+b^2\geq 2ab\), we have
\[\frac{c}{2}\sqrt{3a^2+3b^2+2ab}=\frac{c}{2}\sqrt{2a^2+b^2+(a+b)^2}\geq
  c\sqrt{|ab|}=\sqrt{AC},\] which confirms that these critical
points \(B\) are indeed in the range \(B^2>AC\).

All that remains is to understand the positions of the two
remaining critical points in the case of a parallelogram and a
trapezium. But it follows from Equation
\eqref{eq:app_matrix_cubic} that if
\(\det(\mathbf{M})=\det(\mathbf{N})=0\) for some \(B\) then
\(B\) is a root of Euler's cubic. In other words, the value
\(B=\pm\sqrt{AC}\) is a root if and only if the corresponding
limiting conic is a pair of parallel lines. The parallelogram
has two such limiting conics, the trapezium has one. In fact, we
can completely solve Euler's cubic in both cases.

\bigskip

\noindent {\bf Parallelograms.} For a parallelogram, we have
\(a=-b\) and \(c=-d\), so \(D=E=0\) and the area function
reduces to \[\mathcal{A}=\frac{1}{\sqrt{|B^2-AC|}}.\] Euler's
cubic becomes
\[FB^3=ACFB,\] which has the solutions \(B=0\) (which Euler
noted) and \(B=\pm\sqrt{BC}\). So the roots of Euler's cubic
correspond either to the minimal ellipse or else to the limiting
conics, where the function \(\mathcal{A}\) is not defined.

\bigskip

\noindent {\bf Trapezia.} For a trapezium, we have
\[a=c,\quad b=d,\quad A=C=ab<0,\quad F=a^2b^2,\quad D=E=-\frac{1}{2}ab(a+b)>0.\]
The cubic becomes
\[a^2b^2B^3-a^2b^2(a+b)^2B^2+\left(\frac{3}{2}a^3b^3(a+b)^2-a^4b^4\right)B-\frac{1}{2}a^4b^4(a+b)^2.\]
This factorises as
\[a^2b^2(B-ab)\left(B^2-(a^2+ab+b^2)B+\frac{ab}{2}(a+b)^2\right).\]
The root \(B=ab=-\sqrt{AC}\) is the final one we have not yet
found: it corresponds to the limiting conic given by the
parallel lines \(\geom{A\,D}\cup\geom{B\,C}\). The other two
roots can be found explicitly by solving the quadratic:
\[B=a^2+ab+b^2\pm\sqrt{a^4-a^2b^2+b^4}.\]

\section{Convergence in E563}
\label{app:convergence}

In this appendix, for completeness, I will examine the questions
of convergence of the various infinite series and products which
Euler discusses in the paper.

At the end of \S.4, Euler gives the following series:
\[
  \frac{\pi c}{2\surd 2}\left(1-\frac{1.\,1}{4.\,4}\cdot n^2-\frac{1.\,1.\,3.\,5}{4.\,4.\,8.\,8}\cdot n^4-\frac{1.\,1.\,3.\,5.\,7.\,9}{4.\,4.\,8.\,8.\,12.\,12}\cdot n^6-\text{etc.}\right)
\]
and tells us:
\begin{quote}
  ``which series always converges to a limit, however different
  the axes of the ellipse will have been from one another,
  because \(n\) is always less than unity, and moreover the
  numerical coefficients are vigorously decreasing.''
\end{quote}
Indeed, convergence is immediate from D'Alembert's ratio test if
\(n<1\): the ratios of consecutive terms are
\[\frac{3.\,5}{8.\,8}\,n^2,\quad \frac{7.\,9}{12.\,12}\,n^2,\quad
  \frac{11.\,13}{16.\,16}\,n^2,\quad\text{etc.}\] which converge to
zero as long as \(n<1\). However, it seems odd to say that the
numerical coefficients decrease vigorously, since the sequence
\[\frac{3.\,5}{8.\,8},\quad \frac{7.\,9}{12.\,12},\quad
\frac{11.\,13}{16.\,16},\quad\text{etc.}\] has general term
\(\frac{16n^2-1}{16(n+1)^2}\) which converges to \(1\) as \(n\to
\infty\), so the numerical coefficients
\[\frac{1.\,1.\,3.\,5.\,7.\,9.\,11.\,13...}{4.\,4.\,8.\,8.\,12.\,12.\,16.\,16...}\]
decrease more and more gradually.

In \S.5, Euler handles the case when \(n=1\). Although he does
not discuss convergence, we can see with hindsight that this
series converges when \(n=1\) thanks to Raabe's test. Raabe's
test asserts that \(\sum a_k\) converges if the limit
\(\lim_{k\to \infty}k(a_k/a_{k+1}-1)\) exists and is strictly
greater than \(1\) {\cite[\S 12.2(2)]{Bromwich}}. In our
case,
\[\lim_{k\to\infty} k\left(\frac{a_k}{a_{k+1}} - 1\right)=
  \lim_{k\to\infty}k\left(\frac{16(k+1)^2}{16k^2-1}-1\right)=\frac{32k^2+17k}{16k^2-1}
  = 2.\] The convergence of the infinite product for
\(\frac{\pi}{2\surd 2}\) in \S.5 can be established by
rewriting it as
\[\frac{4.\,4}{3.\,5}\cdot\frac{8.\,8}{7.\,9}\cdot\frac{11.\,13}{12.\,12}\cdots
  = \frac{16}{15} \cdot \frac{64}{63} \cdot \frac{144}{143} =
  \left(1+\frac{1}{15}\right) \left(1+\frac{1}{63}\right)
  \left(1+\frac{1}{143}\right)\cdots\] and {\cite[\S
  39]{Bromwich}} tells us that this converges if and only if the
sum
\[\frac{1}{15}+\frac{1}{63}+\frac{1}{143}+\cdots=\sum_{k=1}^\infty
  \frac{1}{4k^2-1}\] converges. But this sum is bounded above by
\(\sum 1/k^2=\pi^2/6\), and so converges.


\begin{thebibliography}{9}
\bibitem{Bachtold} {\sc B\"{a}chtold, M.} (2018). Answer to ``Who invented the Leibniz notation \(\frac{d^2y}{dx^2}\) for the *second* derivative'', History of Science and Mathematics Stackexchange. Oct. 26, 2018. \verb|https://hsm.stackexchange.com/a/7857/|

\bibitem{Bernoulli} {\sc Bernoulli, D.} (1733). Letter to L. Euler
  22. September 1733. (E864), Published in {\sc G. Enestr\"{o}m}, Der Briefwechsel zwischen Leonhard Euler und Daniel
    Bernoulli (The exchange of letters between Leonhard Euler
  and Daniel Bernoulli), {\em Bibliotheca Mathematica} (1906-1907), 3. Folge 7. Band, pp. 134--153.

\bibitem{Bromwich} {\sc Bromwich, T. J. I.'a.}. (1908). An introduction to the theory of infinite series. MacMillan, London.

\bibitem{CajoriNotation} {\sc Cajori, F.} (1929). A history of mathematical notations, Volume II. (1929) Open Court Publishing, Chicago.
  
\bibitem{CapobiancoEneaFerraro} {\sc Capobianco, G.}, {\sc
    Enea, M. R.} and {\sc Ferraro, G.} (2017). Geometry and analysis in
  Euler's integral calculus. 
  {\em Arch. Hist. Exact. Sci.} 71:1--38.

\bibitem{Dorrie} {\sc D\"{o}rrie, J.} (1965). 100 Great Problems of
  Elementary Mathematics: Their History and Solution,
  transl. D. Antin, Dover.

\bibitem{E52} {\sc Euler, L.} (1741). Solvtio problematis
    rectificationem ellipsis requirentivm. (E52, Solution of a
  problem requiring the rectification of an ellipse.) {\em Commentarii Academiae Scientiarum Petropolitanae (pro Anno
  1736),} Volume VIII, pp. 86--98. Appeared in 1741.
  Reprinted in Opera Omnia, Series 1, Volume 20, pp. 8--20
  (editor Adolf Krazer). Original text available online at
  \verb|https://scholarlycommons.pacific.edu/euler/|

\bibitem{E154} --- (1750). Animadversiones in
    rectificationem ellipsis. (E154, Observations on the rectification of the ellipse.)
  Opuscula Varii Argumenti, Volume 2, pp. 121--166.
  Reprinted in Opera Omnia, Series 1, Volume 20, pp. 21--55
  (editor Adolf Krazer). Original text available online at
  \verb|https://scholarlycommons.pacific.edu/euler/|

\bibitem{E448} --- (1774). Nova series infinita maxime
    convergens perimetrvm ellipsis exprimens. (E448, A new rapidly
  converging infinite series expressing the perimeter of the
  ellipse.)  {\em Novi Commentarii Academiae Scientiarum 
  Petropolitanae (pro Anno 1773),} Volume XVIII, pp. 71--84. Appeared
  in 1774. Reprinted in Opera Omnia, Series 1, Volume 20,
  pp. 357--370 (editor Adolf Krazer). Original text available
  online at \verb|https://scholarlycommons.pacific.edu/euler/|

\bibitem{E563} --- (1784). De ellipsi minima dato
    parallelogrammo rectangvlo circvmscribenda. (E563, On the minimal
  ellipse circumscribing a given right-angled parallelogram.)
  {\em Acta Academiae Scientiarum Imperialis Petropolitanae (pro Anno
  1780),} Volume IV, pp. 3--17. Appeared in 1784.
  Reprinted in Opera Omnia, Series 1, Volume 28, pp. 322--335
  (editor Andreas Speiser). Original text available online at
  \verb|https://scholarlycommons.pacific.edu/euler/|

\bibitem{E691} --- (1795). Problem geometricvm, qvo inter
    omnes ellipses, qvae per data qvatvor pvncta tradvci
    possvnt, ea qvaeritvr, qvae habet aream minimam (E691, A
  geometric problem, in which amongst all ellipses which can be
  drawn through four given points, that one is sought which has
  minimal area.) {\em Nova Acta Academiae Scientiarum Imperialis
  Petropolitanae (pro Anno 1791),} Volume IX, pp. 132--145. Appeared
  in 1795. Reprinted in Opera Omnia, Series 1, Volume 29,
  pp. 226--238 (editor Andreas Speiser). Original text available
  online at \verb|https://scholarlycommons.pacific.edu/euler/|

\bibitem{E692} --- (1795). Solvtio problematis maxime
    cvriosi, qvo inter omnes ellipses, qvae circa datvm
    triangvlvm circvmscribi possvnt, ea qvaeritvr, cvivs area
    sit omnivm minima. (E692, Solution to a problem of greatest
  curiosity in which, amongst all ellipses which can be
  circumscribed around a given triangle, that one is sought,
  whose area is smallest of all.) {\em Nova Acta Academiae Scientiarum
  Imperialis Petropolitanae (pro Anno 1791),} Volume IX,
  pp. 146--153. Appeared in 1795. Reprinted in Opera Omnia,
  Series 1, Volume 29, pp. 239--246 (editor Andreas
  Speiser). Original text available online at
  \verb|https://scholarlycommons.pacific.edu/euler/|

\bibitem{Gergonne} {\sc Gergonne, J. D.} (1826-7). Probl\`{e}mes de g\'{e}om\'{e}trie. (Geometry problems.)
  {\em Annales de Math\'{e}matiques Pures et Appliqu\'{e}s.} Volume
  XVII, p. 284. Original available online at
  \verb|https://www.numdam.org/item/AMPA_1826-1827__17__283_0/|

\bibitem{John} {\sc John, F.} (1948). Extremum problems with
    inequalities as subsidiary conditions. In: Studies and Essays
  Presented to R. Courant on his 60th Birthday
  187--204. Interscience Publishers, Inc., New York, N. Y.
  
\bibitem{Sommerville} {\sc Sommerville, D. M. Y.} (1933). Analytical conics, G. Bell and Sons, London. Third edition.
  
\bibitem{Steiner1} {\sc Steiner, J.} (1827). Aufl\"{o}sung einer geometrischen
    Aufgabe aus Gergonne's Annales de
    Math\'{e}m. t. XVII. p.284. (Solution to a geometrical
  exercise from Gergonne's Annales de
  Math\'{e}m. vol. XVII. p.284.) {\em Crelle's Journal}, Band II,
  pp. 64--65. See J. Steiner, Opera Vol. I, p.123--124,
  available online at
  \verb|https://archive.org/details/jacobsteinersge04steigoog|.

\bibitem{Steiner2} {\sc Steiner, J.} (1828-9). D\'{e}veloppement d'une s\'{e}rie de
    th\'{e}or\`{e}mes relatifs aux sections coniques.
  (Development of a series of relating to conic sections.)
  {\em Annales de Math\'{e}matiques Pures et Appliqu\'{e}s,} Volume
  XIX, pp. 37--64. Original available online at
  \verb|https://www.numdam.org/item/AMPA_1828-1829__19__37_0/|
  
\bibitem{Wallis} {\sc Wallis, J.} (1656). Arithmetica infinitorvm
    sive nova methodus inquirendi in curvilineorum quadraturam
    aliaque difficiliora Matheseos problemata (Arithmetic of
  infinities or new methods of inquiring into the quadrature of
  curvilinear figures and other more difficult mathematics problems.)
  Oxford. Original text avaiable online at
  \verb|https://archive.org/details/ArithmeticaInfinitorum/|
\end{thebibliography}
\end{document}